\documentclass[preprint,12pt]{elsarticle}
\usepackage[english]{babel}
\usepackage[utf8]{inputenc}
\usepackage[T1]{fontenc}

\usepackage{xcolor}
\usepackage{graphicx}
\usepackage{caption}
\usepackage{layout}
\usepackage{hhline}
\usepackage{todonotes}
\usepackage{multirow}
\usepackage{amsmath,amsfonts}

\newtheorem{theorem}{Theorem}

\newtheorem{assumption}[theorem]{Assumption}

\newcommand{\codename}[1]{\textsl{#1}}

\definecolor{chl}{RGB}{35,36,142}  
\newcommand{\hl}[1]{{\color{black}{#1}}}

\begin{document}

\begin{frontmatter}

\title{Coupling parallel adaptive mesh refinement with \\
 a nonoverlapping domain decomposition solver}

\author[prague,munich]{Pavel~K{\r u}s}
\ead{kus@math.cas.cz}
\author[prague,manchester]{Jakub~{\v S}{\'\i}stek}
\ead{sistek@math.cas.cz}

\address[prague]{Institute of Mathematics of the Czech Academy of Sciences, \\{\v Z}itn{\' a} 25, 115 67 Prague, Czech Republic}
\address[munich]{Max Planck Computing and Data Facility, Max Planck Institute,\\ Gießenbachstraße 2, 85748 Garching bei M\"unchen, Germany}
\address[manchester]{School of Mathematics, The University of Manchester, \\ Manchester, M13 9PL, United Kingdom}

\begin{abstract}
We study the effect of adaptive mesh refinement on a parallel domain decomposition solver of a linear system of algebraic equations.
These concepts need to be combined within a parallel adaptive finite element software.
A prototype implementation is presented for this purpose.
It uses adaptive mesh refinement with one level of hanging nodes.
Two and three-level versions of the Balancing Domain Decomposition based on Constraints (BDDC) method are used to solve the arising system of algebraic equations.
The basic concepts are recalled and components necessary for the combination are studied in detail.
Of particular interest is the effect of disconnected subdomains, a typical output of the employed mesh partitioning based on space-filling curves,
on the convergence and solution time of the BDDC method. 
It is demonstrated using a large set of experiments that while both refined meshes and disconnected subdomains have a negative effect on the convergence of BDDC, the number of iterations remains acceptable.
In addition, scalability of the three-level BDDC solver remains good on up to a few thousands of processor cores.
The largest presented problem using adaptive mesh refinement has over 10$^9$ unknowns and is solved on 2048 cores.
\end{abstract}


\end{frontmatter}


\section{Introduction}
\label{sec:Introduction}
Adaptive mesh refinement (AMR) is a well-established technique in the framework of the finite element method (FEM). 
Its use in combination with massively parallel calculations is, however, still rather limited. 
This should not be surprising, since developing an efficient parallel code employing a general AMR is highly nontrivial and
requires meeting several contradicting requirements. 
In particular, keeping the load of CPU cores balanced is the main issue during the parallel process. 
Other issues involve management of hanging nodes (or enforced refinements) and a choice of a suitable algebraic solver.

In this contribution, we present an implementation of a massively parallel AMR code.
It is coupled to the parallel linear system solver based on nonoverlapping domain decomposition method. 
The impact of AMR on the solver is evaluated.
A subsequent goal of the paper is convincing the reader that hanging nodes treatment and its implementation can be relatively easy and straightforward, 
as opposed to the general feeling that nonconforming meshes are best to be avoided. 
The simplicity is achieved by allowing only so called first-level hanging nodes and equal order elements.

The goal of an AMR strategy is to improve the approximate solution in those parts 
of the domain, where the behaviour of the exact solution is complex, such as close to singularities or within boundary and internal layers. 
The structure of the mesh becomes complicated while the degrees of freedom are efficiently distributed and their number is kept low.
On the other hand, 
relatively simple, often structured, globally fine meshes are typically preferred in massively parallel codes,
resulting in very large number of degrees of freedom.

An integral part of parallel computations is partitioning the computational mesh into subdomains. 
The strategy employed in our work leads to parts with approximately equal number of elements.
Although this requirement may seem rather simple, 
partitioning adaptively refined meshes in a scalable way on thousands of CPU cores is still a challenging task. 

Perhaps the simplest strategy suitable for structured meshes is dividing them into geometrically regular subdomains.
However, it cannot be used for adaptively refined subdomains, 
where the number of elements in different regions can differ largely, leading to huge load imbalance. 
A widely adopted alternative is translating the problem of dividing the mesh into partitioning the graph of the mesh. 
The graph partitioning is then performed by standard libraries, such as \codename{METIS}~\cite{Karypis:1998-FHQ} or \codename{ParMETIS}~\cite{Karypis-1999-PMS}.
While the latter is meant for parallel repartitioning, the approach is not scalable to thousands of cores as required by massively parallel adaptive simulations. 

A considerably simpler approach is based on partitioning space-filling curves, 
as it is done in the \codename{p4est} library~\cite{burstedde_p4est:_2011, isaac_recursive}, one of the building blocks of our implementation. 
Its scalability has been tested in~\cite{burstedde_p4est:_2011} up to hundreds of thousands of cores, 
and its use for development of a parallel FEM library \codename{deal.II} is described in~\cite{bangerth_algorithms_2012}. 
\hl{The authors use a globally assembled system matrix distributed by rows.}
The main difference of our work seems to be the use of the \emph{subassembled} system matrix, i.e.\ matrix only assembled subdomain-wise,
as it is common to nonoverlapping domain decomposition methods.
\hl{An interface degree of freedom is present in two or more subdomains, and its assembly is not finalised over the interface.
The use of this matrix format naturally avoids partitioning of matrix rows of the fully assembled matrix and thus circumvents the issue with
the assembly at hanging nodes in the vicinity of interface described in~\cite{burstedde_p4est:_2011}.}

On the other hand, properties of subdomains play a more important role in our approach.
The drawback of the use of space-filling curves for mesh partitioning is the poor shape of subdomains. 
Moreover, disconnected subdomains composed of several independent components are common. 
One of our goals is to investigate how this type of subdomains affects the performance of the nonoverlapping domain decomposition solver.

For the solution of the system of equations, we use the Balancing Domain Decomposition based on Constraints (BDDC) method~\cite{Dohrmann-2003-PSC} 
and its extension to multiple levels, the \emph{Multilevel BDDC}~\cite{Mandel-2008-MMB,Tu-2007-TBT3D}.
The potential of the multilevel method to scale to 500 thousand cores was recently demonstrated in~\cite{Badia-2016-MBD}.
Another parallel implementation of multilevel BDDC is available in our open-source \codename{BDDCML} library~\cite{Sousedik-2013-AMB}, 
the second building block of our implementation.

A crucial role in BDDC is played by constraints which enforce continuity of suitably defined \emph{coarse degrees of freedom} across subdomains. 
Typically, these are point values at selected nodes called \emph{corners} or 
averages over faces and/or edges among neighbouring subdomains.
If a sufficient number of constraints is properly selected, the local Neumann problems on subdomains become uniquely solvable.

The strategy for handling disconnected subdomains employed in \codename{BDDCML} is described. 
In fact, it is quite simple: 
a local graph of computational mesh is created for each subdomain and its continuity is analyzed. 
If more graph components are detected, each of them is treated separately during the selection of constraints.
While this may lead to a certain load imbalance due to a sudden increase of number of constraints for disconnected subdomains, 
the principal amount of work, which is given by the size of each subdomain, remains unchanged.
This effect is studied in detail in our paper.


The rest of the paper is organized as follows. 
In Section~\ref{sec:Adaptive}, we review adaptive mesh refinement strategies and describe a simple approach to dealing with hanging nodes. 
The BDDC method is recalled in Section~\ref{sec:Balancing} with the emphasis on accommodating disconnected subdomains within the method.
Section~\ref{sec:Coupling} is devoted to combination of these concepts and discussion of specific issues arising in coupling AMR with a parallel BDDC solver.
Finally, numerical results testing the developed implementation are presented in Section~\ref{sec:Numerical},
and conclusions are drawn in Section~\ref{sec:Conclusions}.

\section{Adaptive meshes and hanging nodes}
\label{sec:Adaptive}

\begin{figure}[tbh]
\centering
 \includegraphics[width = 0.45\textwidth]{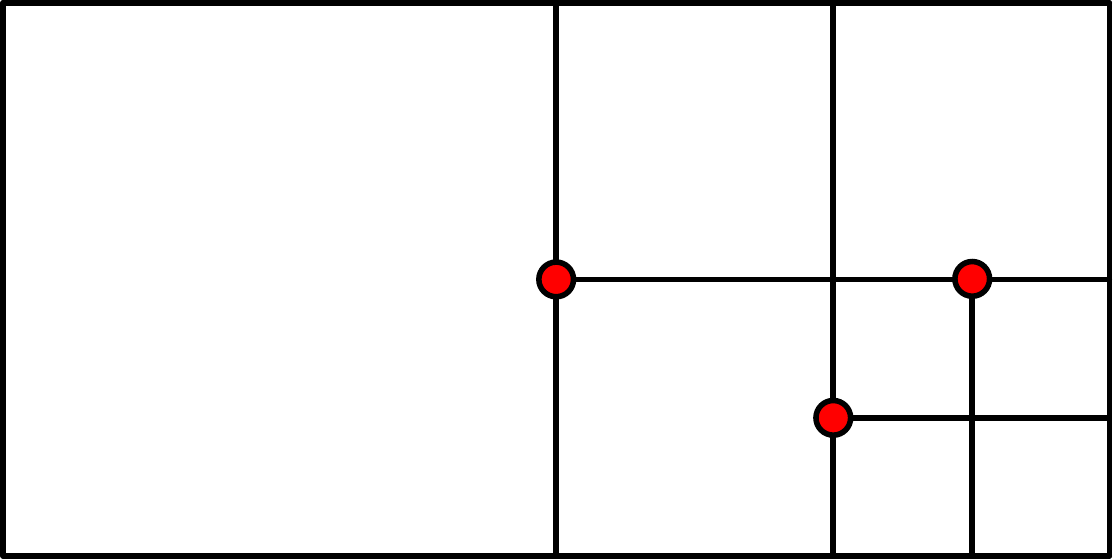}
 \hskip 5mm
 \includegraphics[width = 0.45\textwidth]{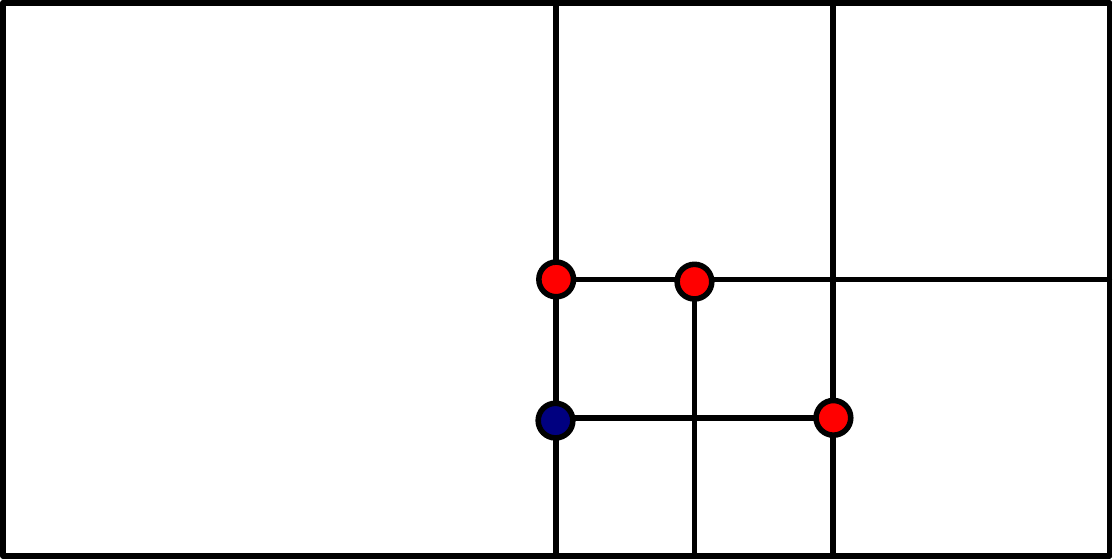} 
 \caption{Mesh satisfying the 2:1 regularity of Assumption~\ref{ass:1-2} (left) contains only \emph{first-level} hanging nodes (red dots). 
Mesh not satisfying the 2:1 regularity (right) contains also \emph{higher-level} hanging nodes (blue dot).
To comply with the assumption, the refined element would trigger an enforced refinement of the large (left) element in this case.
}
 \label{fig:1to2}
\end{figure}

The purpose of mesh adaptivity is to ensure precise resolution of details in troublesome areas of the domain 
while keeping the overall size of the resulting system within reasonable bounds.
This often cannot be achieved by uniform refinements. 
The price to pay are several complications of the discretisation algorithm which have to be addressed. 
This is particularly true in the case of adaptivity in parallel setting as will be further elaborated in Section~\ref{sec:Coupling}. 
There is, however, one particular issue, which is rather challenging even in the case of serial calculation: 
a proper treatment of the so called \emph{hanging nodes}. 
These nodes appear when several smaller elements are adjacent to an edge (in 2D) or a face (in 3D) of a larger element, 
see Fig.~\ref{fig:1to2}.
It is especially demanding when higher-order elements are used.
Various techniques have been developed and used in different settings in the case of serial calculations. 
However, one has to be careful when dealing with adaptivity in parallel 
where hanging nodes might appear at subdomain interface. 

Hanging nodes present a complication of the numerical scheme, and many authors proposed various techniques to avoid them 
from the beginning by introducing extra refinements. 
These \emph{enforced} refinements are not necessary for a better precision of the solution but have the only purpose of keeping the mesh regular, 
i.e.\ face-to-face.
Let us just mention the red--green algorithm, although many others exist.
Nevertheless, these algorithms bring other disadvantages, and it seems that most approaches towards mesh adaptivity involve hanging nodes nowadays.
The presence of hanging nodes is quite simply manageable when a discontinuous approximation, e.g.\ the discontinuous Galerkin method, is used. 
This can be employed for development of hybrid continuous-discontinuous methods as in \cite{badia_adaptive_2013}, 
where discontinuous approximation is used at hanging nodes. 

If continuous approximation is used, the use of  
\emph{arbitrary-level} hanging nodes (see e.g. \cite{solin_arbitrary-level_2008} for 2-D and \cite{kus_arbitrary-level_2014} for 3-D results or recent 
works \cite{Zander-2016}  or  \cite{DiStolfo-2016} for alternative approaches) becomes technically somewhat difficult. 
The \emph{level} here corresponds to the number of subsequent refinements at one side of an element face that were needed for creation of the hanging node,
see Fig.~\ref{fig:1to2} for first and second level hanging nodes.
The advantage of the arbitrary-level approach is that there are no additional refinements enforced by mesh regularity reasons. It leads, 
however, to a substantial algorithmic complexity due to the non-local character of constraints on continuity at 
higher-level hanging nodes and their propagation through the mesh \cite{kus_arbitrary-level_2014}.
This holds especially in 3D. From this reason, most authors resort to the use of first-level hanging nodes only (see, e.g., \cite{Demkowicz-book2}), 
which seems to be the most feasible solution and which is used in our current work as well.

In the rest of this section, we describe in more detail the way hanging nodes are treated in our solver.
The ideas are not new, but we recall them for the sake of completeness and in a way suitable for a subsequent coupling with nonoverlapping domain decomposition.
We proceed by imposing two important restrictions, formulated as assumptions on the mesh.

\begin{assumption}[2:1 mesh regularity]
\label{ass:1-2}
Only first-level hanging nodes are present in the adaptive meshes. 
\end{assumption}
This important assumption means that no more than two (four) other elements can be adjacent to an element edge (face) in 2D (3D), respectively, see Fig.~\ref{fig:1to2}. 
Consequently, a limited number of enforced refinements can appear in order to fulfill this assumption.
This restriction is used by a majority of authors as a reasonable trade-off between performance gain and complexity of implementation. 

\begin{assumption}[Equal order of elements]
\label{ass:equal_order}
Finite elements have a uniform polynomial order. 
\end{assumption}
This restriction rules out the superior $hp$-adaptivity and leaves us with the $h$-adaptivity, 
however including powerful higher-order approximations.

\begin{figure}[tbh]
\centering
 \includegraphics[width = 0.495\textwidth]{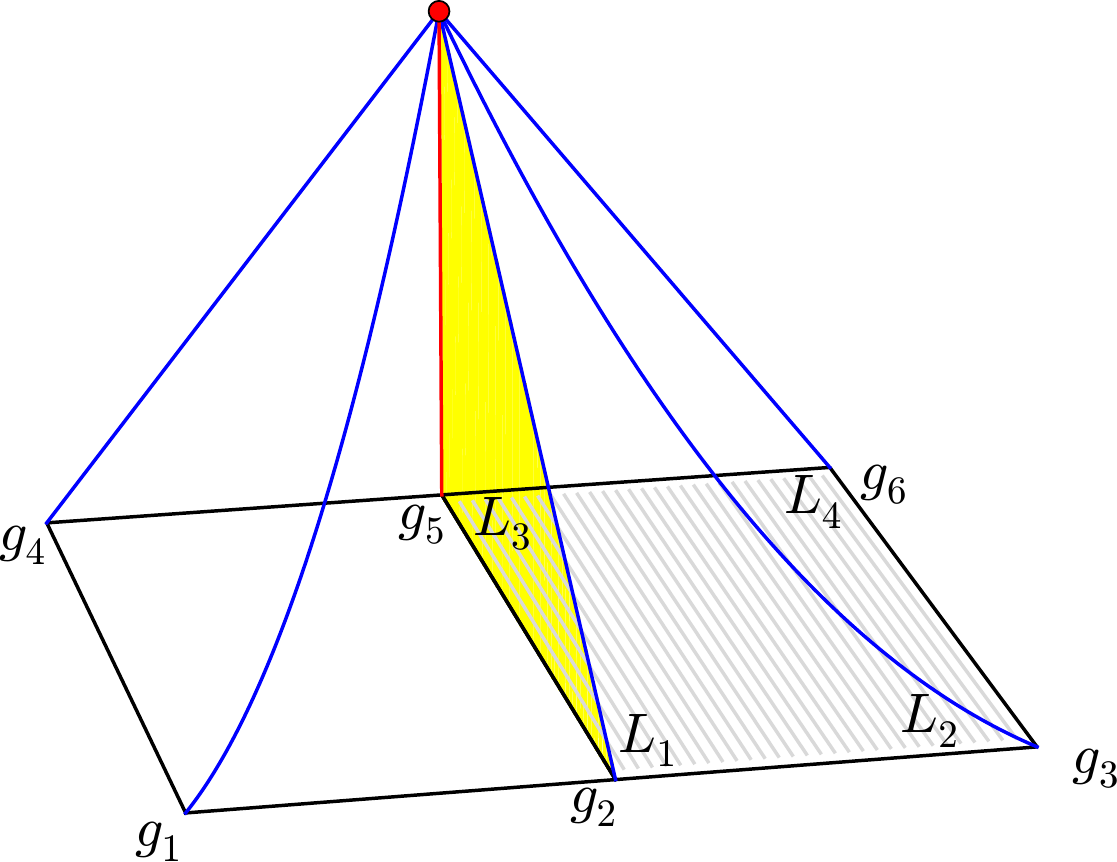}
 \includegraphics[width = 0.495\textwidth]{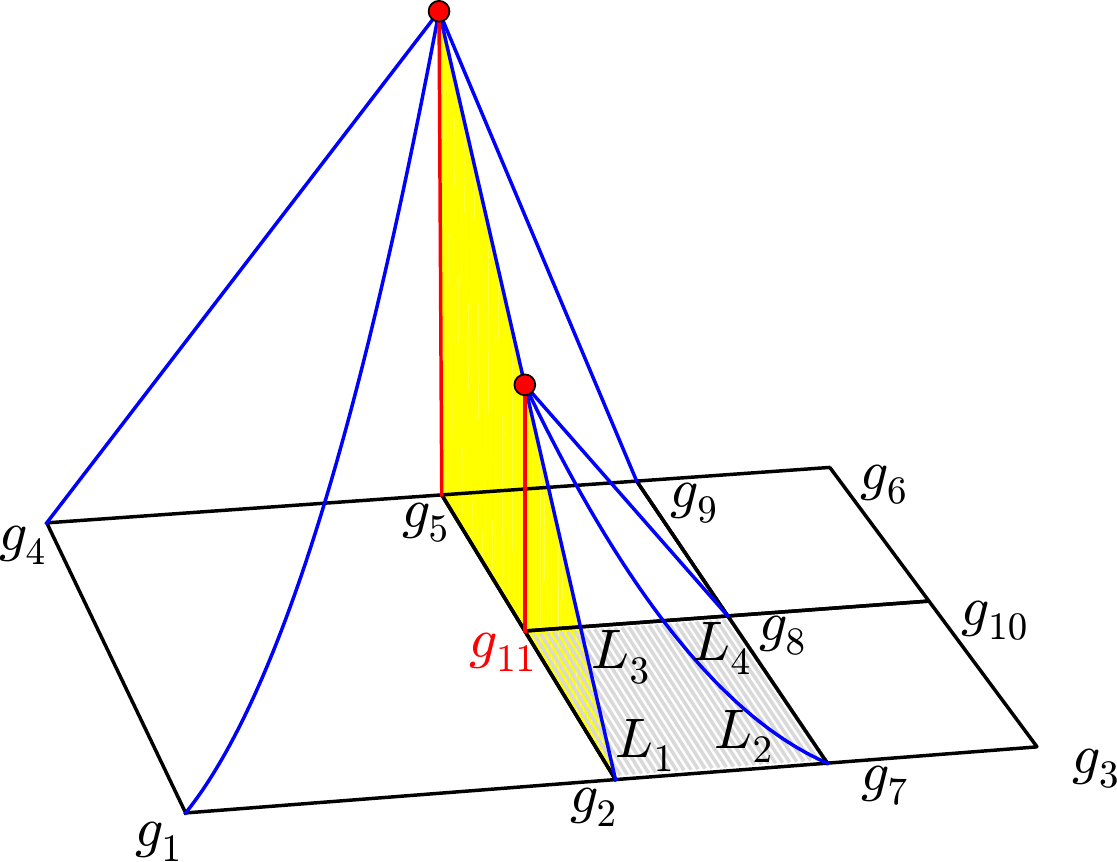} 
 \caption{\label{fig:basis_function_glob_constrain}Linear basis function in regular (left) and adapted (right) meshes. The figure represents the first approach, where 
 a global degree of freedom $g_{11}$ is assigned to the hanging node and needs to be handled after the global assembly.}
\end{figure}

\begin{assumption}[Isotropic refinements]
\label{ass:isotropic}
Elements are only refined into 4 (2D) or 8 (3D) smaller elements in a most natural way by bisecting all element edges.
\end{assumption}
Although anisotropic refinements might be useful e.g.\ near to boundary layers, the potential performance gain is rather limited.

\subsection{Finite element assembly}

We will illustrate the finite element assembly on an example of a simple mesh of (initially) two bilinear elements, see Fig.~\ref{fig:basis_function_glob_constrain} (left). 
If we do not consider boundary conditions, as if the two-element patch was inside the computational domain, 
there are 6 degrees of freedom denoted $g_1, g_2, \ldots, g_6$. 
In standard finite element codes, the global stiffness matrix is built element-wise, 
first constructing local stiffness matrix of each element (4$\times$4 in our example) followed by the \emph{assembly}, 
i.e.~distributing the values into appropriate places in the global matrix (6$\times$6 in the example). 
The local stiffness matrix is obtained using shape functions local to the element, which correspond to local degrees of freedom, such as
$L_1,\ldots, L_4$ indicated on the hatched element $\Omega _{j}^{e}$ in Fig.~\ref{fig:basis_function_glob_constrain} (left).
Denoting shape functions by the same symbols as corresponding degrees of freedom, 
the local stiffness matrix on element $\Omega _{j}^{e}$ can be written as 
\begin{equation}
\label{eq:local_assembly}
A_{j}^{e}=
\begin{bmatrix}
    a(L_1, L_1) & \dots & a(L_1, L_4) \\
    \vdots & \ddots & \vdots\\
    a(L_4, L_1) & \dots & a(L_4, L_4) \\
\end{bmatrix}
\end{equation}
where $a(\cdot,\cdot)$ is the bilinear form in a weak formulation, e.g.~of the Poisson problem used in our computations in Section~\ref{sec:Numerical}.
The global assembly procedure can be formally written as 
\begin{equation}
\label{eq:global_assembly}
A = \sum_{j=1}^{N_e} R_{j}^{eT} A_{j}^{e} R_{j}^{e},
\end{equation}
where $A$ is the global stiffness matrix, $N_e$ is the number of elements in the mesh, 
and $R_{j}^{e}$ is the 0--1 matrix restricting global vector of unknowns to those present on the $j$-th element.
In our example, $R_{j}^{e}$ for the hatched element is 
\begin{equation}
R_{j}^{e}=
\begin{bmatrix}
    0 & 1 & 0 & 0 & 0 & 0 \\
    0 & 0 & 1 & 0 & 0 & 0 \\
    0 & 0 & 0 & 0 & 1 & 0 \\
    0 & 0 & 0 & 0 & 0 & 1 \\
\end{bmatrix}
\end{equation}
In an actual implementation, matrix $R_{j}^{e}$ is implemented as a vector of global indices 
expressing the local-to-global correspondence of degrees of freedom, in our case 
\begin{equation}
\label{loc_glob_global_constrain}
(L_1, L_2, L_3, L_4) \leftrightarrow (g_2, g_3, g_5, g_6). 
\end{equation}

Similarly, the global right-hand side vector $f$ is assembled from the finite element contributions $f^e$ as
\begin{equation}
\label{eq:global_rhs_assembly}
f = \sum_{j=1}^{N_e} R_{j}^{eT} f_{j}^{e}.
\end{equation}
Finally, the global system of equations is 
\begin{equation}
\label{eq:global_system}
 Au = f,
\end{equation}
where $u$ is the global vector of unknown coefficients of finite element basis functions.
In general, application of suitable boundary conditions makes system (\ref{eq:global_system}) uniquely solvable.

\subsection{Continuity of the solution at a hanging node}
\label{sec:global_contraints}

When hanging nodes are present in the discretisation, the situation gets slightly more complicated. 
As can be seen in Fig.~\ref{fig:basis_function_glob_constrain} (right), we can no longer allow an arbitrary value 
in $g_{11}$ if we want to keep the solution continuous, which is the case for the considered conforming finite element method.
There are generally two possibilities how this can be achieved, which are both described in this section. 
Although the approaches are equivalent mathematically, their algorithms differ.

In the \emph{global approach}, degrees of freedom are kept in the hanging nodes, and the same local (\ref{eq:local_assembly}) and global 
(\ref{eq:global_assembly}) matrices are
used to create the global system. 
In order to obtain continuous results, additional constraints have to be satisfied by the degrees of freedom. 
In our example of Fig.~\ref{fig:basis_function_glob_constrain},
\begin{equation}
\label{eq:correspondence_global}
(L_1, L_2, L_3, L_4) \leftrightarrow (g_2, g_7, {\mathbf{g_{11}}}, g_8),
\end{equation}
and the continuity leads to the constraint
\begin{equation}
\label{node_constraint}
  g_{11} = \frac{g_2 + g_5}{2}.  
\end{equation}
In a more complicated mesh with more hanging nodes, we would have to add a number of such constraints, 
which can be written in the form of a matrix equation
\begin{equation}
\label{eq:global_constraints}
Gu = 0.
\end{equation}
Matrix $G$ is typically rectangular with much less rows than columns. 
In our example, (\ref{node_constraint}) translates to 
\begin{equation}
G = \begin{bmatrix} 0, \frac{1}{2}, 0, 0, \frac{1}{2}, 0, 0, 0, 0, 0, -1 \end{bmatrix}.
\end{equation}
System (\ref{eq:global_system}) now has to be solved subject to the constraint (\ref{eq:global_constraints}),
or in other words, satisfying $u \in \mathrm{null}\,G$.
Various techniques are possible, see e.g. \cite{bangerth_data_2009}
for more details of this approach for $hp$-adapted meshes.

\begin{figure}[tbh]
\centering
 \includegraphics[width = 0.495\textwidth]{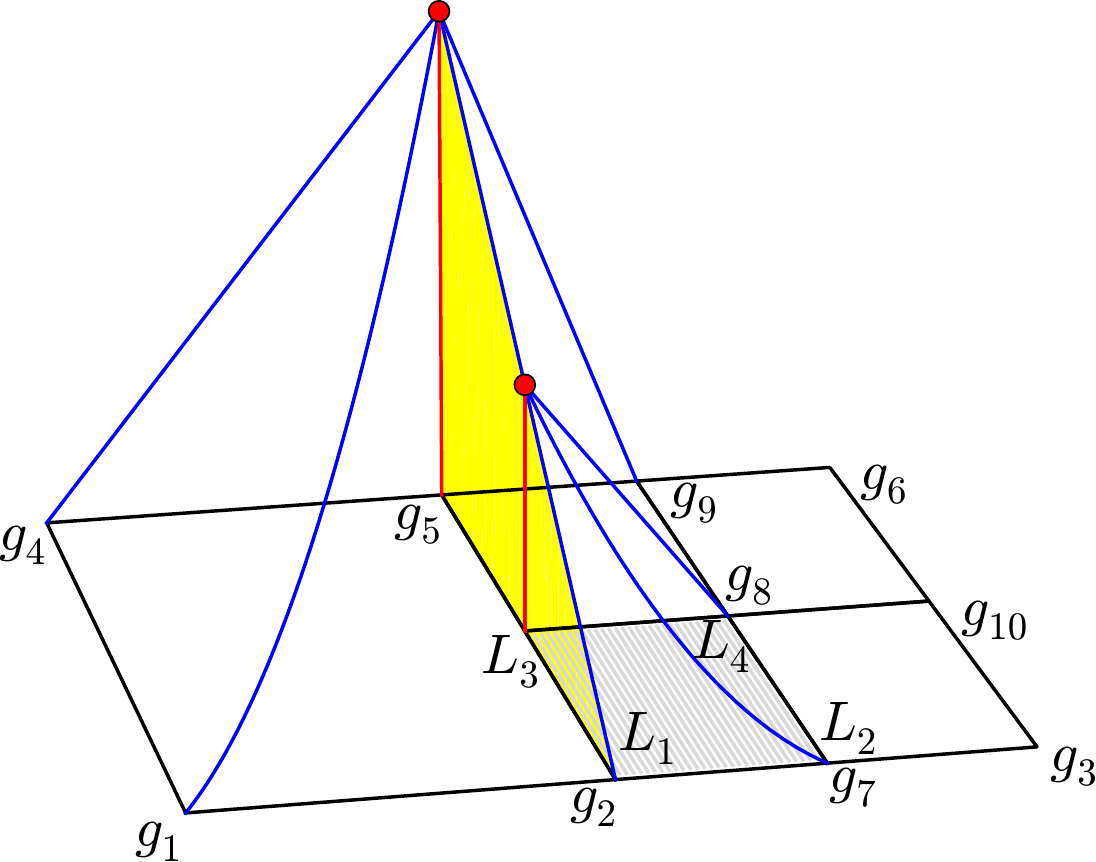}
 \includegraphics[width = 0.49\textwidth]{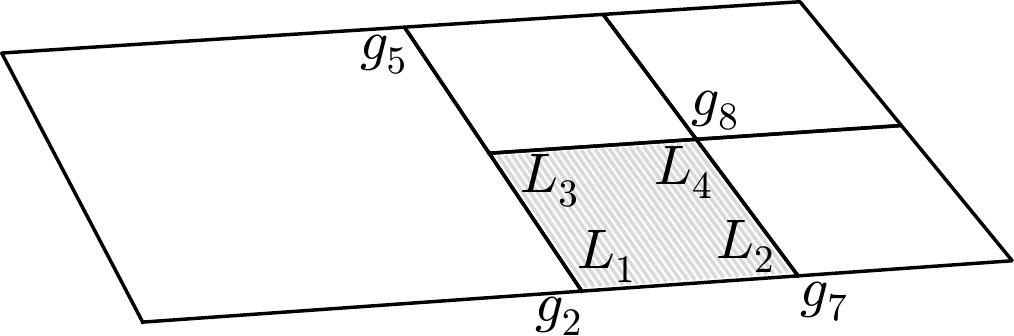} 
 \caption{\label{fig:basis_function_loc_constrain}The second approach to handling hanging nodes --- no global degree of freedom is assigned to the hanging node. There is no global degree of freedom associated with $L_3$. Instead refined elements contribute to non-local global nodes during the assembly.}
\end{figure}

The second, \emph{local approach}, is more suitable for our situation, and it is employed in our implementation. 
It relies on elimination of hanging nodes on the element level before the global matrix assembly. 
The important difference from the global approach is that we do not assign global degrees of freedom to hanging nodes. 
Instead, we directly apply the constraints to create continuous basis functions, 
which in turn give rise to a continuous approximate solution.

Due to Assumptions~\ref{ass:1-2}--\ref{ass:isotropic} 
one can keep simple relationship between local and global degrees of freedom similar to (\ref{loc_glob_global_constrain}). 
In particular, each element contributes to as many global degrees of freedom as is the number of local degrees of freedom. 
In our example, correspondence~(\ref{eq:correspondence_global}) will be replaced by
\begin{equation}
\label{loc_glob_local_constrain}
(L_1, L_2, L_3, L_4) \leftrightarrow  (g_2, g_7, {\mathbf{g_5}}, g_8).
\end{equation}
Note that instead of the global degree of freedom $g_{11}$, which is not considered now, the element will contribute to $g_5$. 
The reason can be seen from Fig.~\ref{fig:basis_function_loc_constrain}. 
We eliminate $g_{11}$ directly during the assembly by distributing contributions from local degree of 
freedom $L_3$ between $g_2$ and $g_5$. 
This can be also understood as applying first the local transform
\begin{equation}
\bar{A}_{j}^e = T_j^T A_j^e T_j
\end{equation}
given for the example in Fig.~\ref{fig:basis_function_loc_constrain} by the transition matrix
\begin{equation}
\label{eq:transform}
T_j=
\begin{bmatrix}
    1   & 0 & 0   & 0 \\
    0   & 1 & 0   & 0 \\
    1/2 & 0 & 1/2 & 0 \\
    0   & 0 & 0   & 1 \\
\end{bmatrix}.
\end{equation}
Here $\bar{A}_{j}^e$ is a modified local stiffness matrix. 
The transform is followed by an assembly procedure 
\begin{equation}
A = \sum_{j=1}^{N_e} \bar{R}_{j}^{eT} \bar{A}_{j}^{e} \bar{R}_{j}^{e},
\end{equation}
where $\bar{R}_{j}^{e}$ is now adjusted according to (\ref{loc_glob_local_constrain}).

The matrix $T_j$ is of square shape and the same dimension as $A_j^e$ due to Assumptions~\ref{ass:1-2}--\ref{ass:isotropic}.
Both $T_j$ and $\bar{R}_{j}^{e}$ would get more complicated for more general meshes with different polynomial orders 
and higher-level hanging nodes. 
Such cases are elaborated e.g.\ in \cite{kus_arbitrary-level_2014,solin_arbitrary-level_2008,DiStolfo-2016}.

\subsection{Higher order basis functions}
\begin{figure}[tbh]
\centering
 \includegraphics[width = 0.495\textwidth]{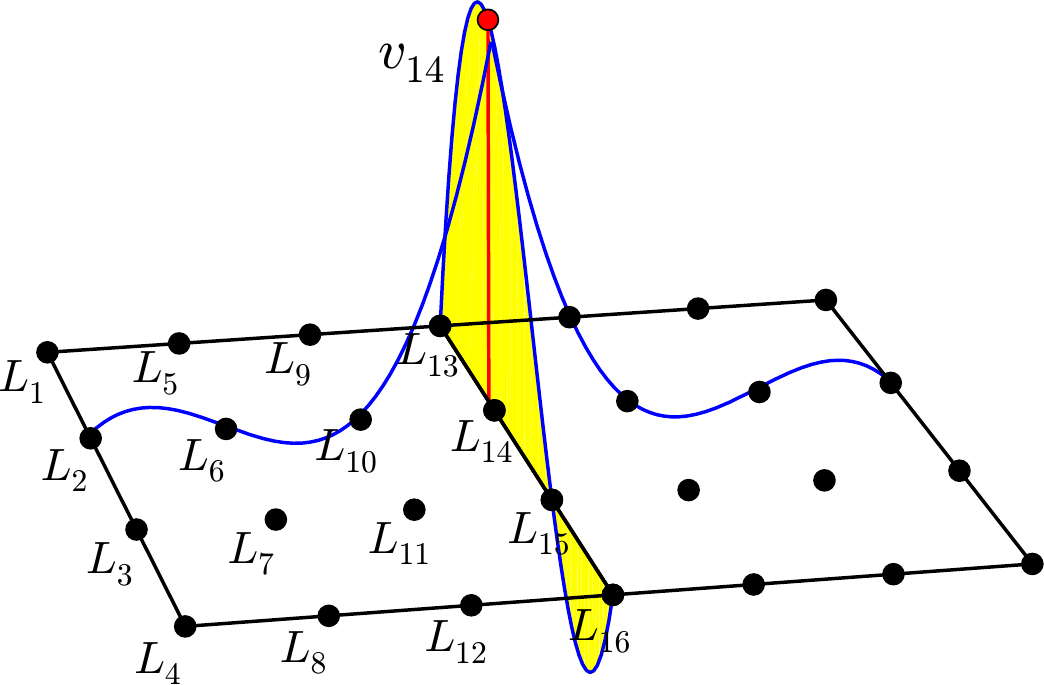}
 \includegraphics[width = 0.495\textwidth]{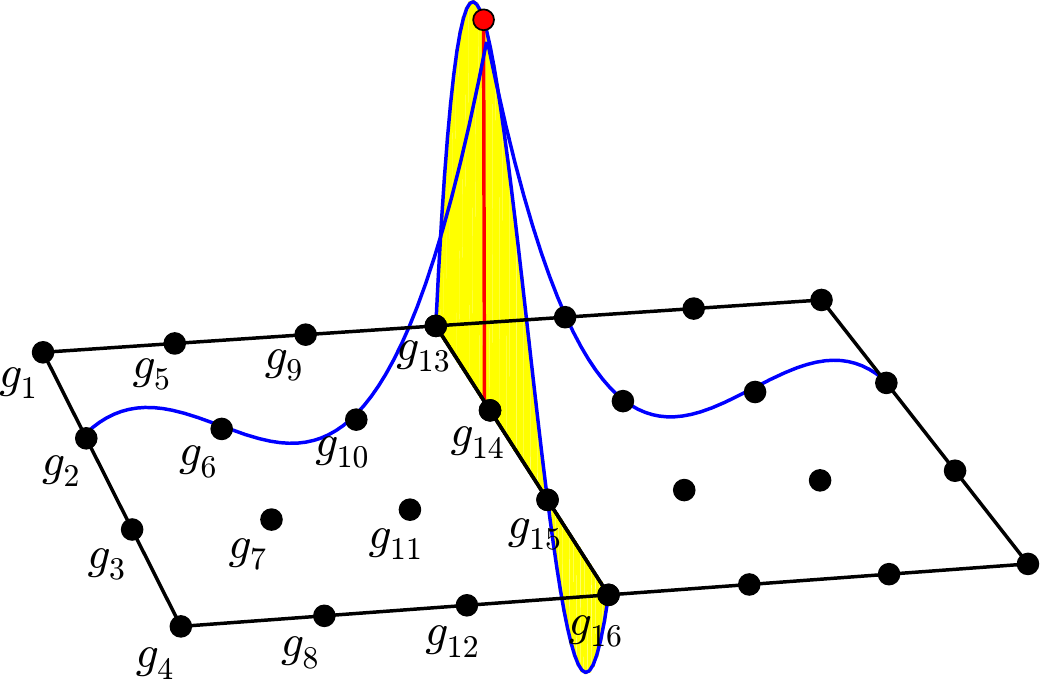} \\
 \includegraphics[width = 0.495\textwidth]{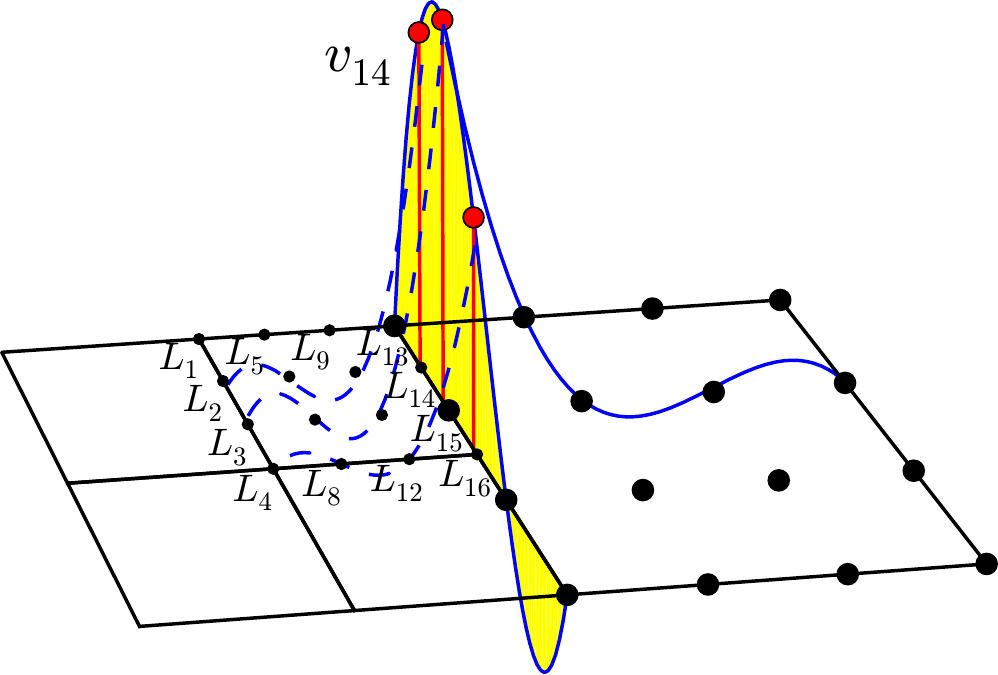}
 \includegraphics[width = 0.495\textwidth]{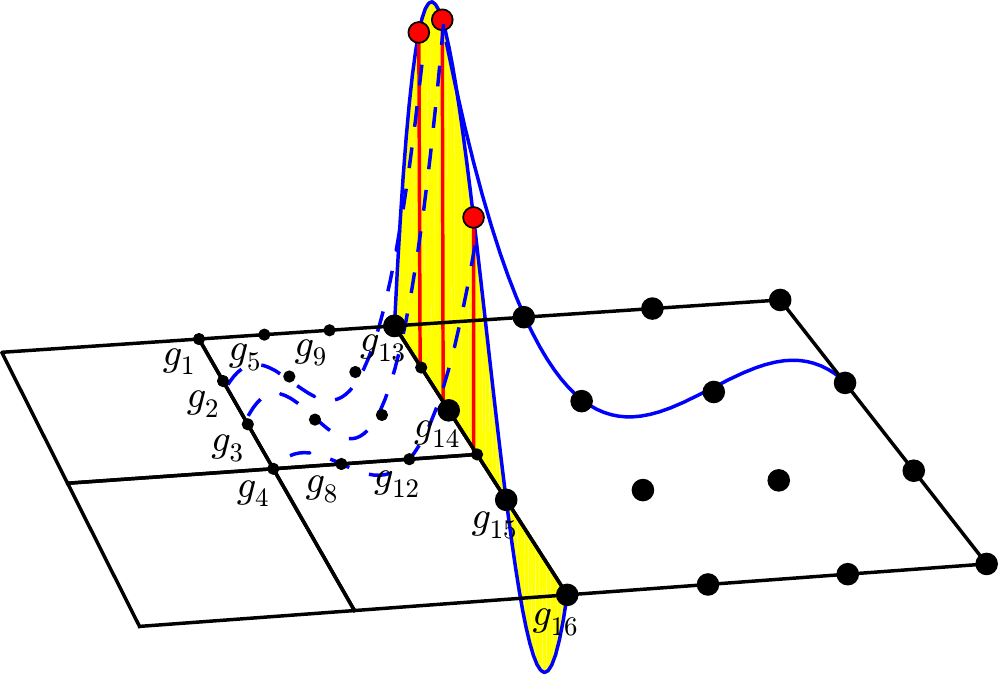} 
 \caption{\label{fig:basis_function_loc_constrain_higher}
Example of the basis function $v_{14}$ of fourth order before (top) and after (bottom) the refinement of the left element. 
Each local degree of freedom from the left element corresponds to exactly one global degree of freedom ($L_{14} \rightarrow g_{14}$) in both cases.}
\end{figure}

In the previous section, 
we explained in detail how the transition matrix $T_j$ can be constructed in the case of bilinear elements. 
In the case of higher-order basis functions, and still under Assumptions~\ref{ass:1-2}--\ref{ass:isotropic}, 
we can proceed in a very similar way. In the top row of
Fig.~\ref{fig:basis_function_loc_constrain_higher}, two elements from fourth-order regular mesh are shown together with an example 
of the basis function~$v_{14}$.
This global basis function corresponds to the degree of freedom~$g_{14}$, which lies on the edge shared by two elements, so two local basis functions 
contribute to the construction of the global basis function. 
From the left element, it is the shape function corresponding to the local degree of freedom~$L_{14}$.
The correspondence of local and global degrees of freedom for the left-hand side element is trivial: 
\begin{equation}
\label{loc_glob_local_higher}
 (L_1, L_2, \ldots, L_{16}) \leftrightarrow  (g_1, g_2, \ldots, g_{16})
\end{equation}
and the transformation matrix $T_j$ constructed for this element is the identity.

Let us now describe the situation of hanging nodes, as it is shown in the bottom row of Fig.~\ref{fig:basis_function_loc_constrain_higher}, where the investigated
element has been refined. Let us focus on one of the newly created elements as shown in the left part of the figure, where interior local degrees 
of freedom are omitted due to a lack of space. In the right-hand side of the figure we can see newly enumerated global degrees of 
freedom $g_1,g_2,\ldots,g_{12}$ which correspond to newly created nodes. Global degrees of freedom $g_{13},\ldots,g_{16}$, however, remain the same,
since the right-hand side element was not refined. Thanks to the new usage of $g_1,g_2,\ldots,g_{12}$, the relation~(\ref{loc_glob_local_higher})
will still hold meaning that the discussed element will contribute to global degrees of freedom $g_1,g_2,\ldots,g_{16}$. The transformation matrix, 
however, now has a more complicated form 
\begin{equation}
\label{eq:transform_higher}
T_j=
\begin{bmatrix}
    1 & 0 & 0 & \dots  & 0 & 0 & 0 & 0 \\
    0 & 1 & 0 & \dots  & 0 & 0 & 0 & 0 \\
    0 & 0 & 1 & \dots  & 0 & 0 & 0 & 0 \\
   \vdots & \vdots & \vdots &  \ddots & \vdots & \vdots & \vdots & \vdots \\
    
    0 & 0 & 0 & \dots  & v_{13}(L_{13}) & v_{14}(L_{13}) & v_{15}(L_{13}) & v_{16}(L_{13}) \\
    0 & 0 & 0 & \dots  &  v_{13}(L_{14}) & v_{14}(L_{14}) & v_{15}(L_{14}) & v_{16}(L_{14}) \\
    0 & 0 & 0 & \dots  &  v_{13}(L_{15}) & v_{14}(L_{15}) & v_{15}(L_{15}) & v_{16}(L_{15}) \\
    0 & 0 & 0 & \dots  & v_{13}(L_{16}) & v_{14}(L_{16}) & v_{15}(L_{16}) & v_{16}(L_{16}) \\
\end{bmatrix}
\end{equation}
The structure of the matrix is the following. Each row corresponds to a local degree of freedom~$L_i$, each column to a global basis function~$v_{k}$ 
corresponding to a global degree of freedom~$g_k$ and the matrix entry is the value of the respective function in the respective node~$v_{k}(L_i)$.
Note that in the first part of the matrix, nodes are regular and thus the basis function is 1 in the corresponding node and vanishes in all other 
nodes, leading to the identity matrix. Also note that the structure of the matrix is the same (although more complicated) as in (\ref{eq:transform}). 
Neither this matrix is constructed explicitly during the calculation.

This correspondence of (Lobatto-type) local degrees of freedom on constrained elements is well described in paper~\cite{isaac_recursive} during the 
discussion of \texttt{Lnodes} extension of the \codename{p4est} library. It can be shown that under Assumptions~\ref{ass:1-2}--\ref{ass:isotropic}, 
the matrix $T_j$ can always be constructed for both two and three dimensions and that it will be a square matrix.
Another benefit of the local approach towards hanging node elimination is that it is compatible with the
\texttt{Lnodes} extension of the \codename{p4est} library and thus global numbering of degrees of freedom of the distributed mesh can be taken directly from \codename{p4est}. 
See Section~\ref{sec:p4est} for details.

\section{Balancing Domain Decomposition based on Constraints}
\label{sec:Balancing}

The BDDC method~\cite{Dohrmann-2003-PSC} is used as a recent example of nonoverlapping domain decomposition (DD) techniques.
The starting point of these methods is a division of the computational domain $\Omega$ and the computational mesh into $N_S$ nonoverlapping subdomains,
$\Omega_i,\ i = 1,\dots,N_S$.
We use $\Omega$ and $\Omega_{i}$ also to denote the set of all degrees of freedom and those present in the $i$-th subdomain, respectively.
The \emph{interface} $\Gamma$ is formed as the set of unknowns that are present in more than one subdomain.

\subsection{Iterative substructuring}
\label{sec:Substructuring}

\hl{The standard approach employed by nonoverlapping DD methods is called \emph{iterative substructuring}, and it is described in detail, e.g.,\ in the monograph~\cite{Toselli-2005-DDM}.}
\hl{The solution} of the global system of equations arising from the finite element discretisation is split into solving a local \emph{discrete Dirichlet problem} 
on each subdomain and a global interface problem with the \hl{\emph{Schur complement with respect to the interface}} and the \emph{reduced right-hand side vector}.

In more detail,
suppose that $A_{i}$ is the local stiffness matrix assembled only from finite elements in the subdomain $\Omega _i$.
This is sometimes called \emph{subassembly}, and the 
assembly could be finalized as 
\begin{equation}
\label{eq:subassembly}
A = \sum_{i=1}^{N_S} R_i^T A_{i} R_i,
\end{equation}
where $R_{i}: \Omega \rightarrow \Omega_i$ is the 0--1 restriction matrix selecting unknowns local to subdomain $\Omega_{i}$ from the global vector of unknowns.
Similarly, let $f_{i}$ be a local right-hand side vector obtained by integration over the same set of elements and a subsequent subassembly.
Again, 
\begin{equation}
\label{eq:fsubassembly}
f = \sum_{i=1}^{N_S} R_i^T f_{i}.
\end{equation}
Note the similarity of (\ref{eq:subassembly}) with (\ref{eq:global_assembly}) and of (\ref{eq:fsubassembly}) with (\ref{eq:global_rhs_assembly}).


Let us denote the set of local interface unknowns as $\Gamma_{i} = \Gamma \cap \Omega_{i}$. 
Then $R^{\Gamma}_i:\Gamma \rightarrow \Gamma_{i}$ stands for the 0--1 matrix selecting local interface unknowns from the set of global interface unknowns. 

Splitting the local degrees of freedom into \emph{interior} (superscript $^{I}$) and \emph{interface} (superscript $^{\Gamma}$) ones implies a 2$\times$2 blocking of the 
subdomain problem on the $i$-th subdomain
\begin{equation}
\left[
\begin{array}[c]{cc}
A_{i}^{II}       & A_{i}^{I\Gamma} \\
A_{i}^{\Gamma I} & A_{i}^{\Gamma\Gamma}  
\end{array}
\right]
\left[
\begin{array}[c]{cc}
u_{i}^{I} \\
u_{i}^{\Gamma}  
\end{array}
\right]
=
\left[
\begin{array}[c]{cc}
f_{i}^{I} \\
f_{i}^{\Gamma}  
\end{array}
\right].
\end{equation}
The \emph{local Schur complement} with respect to $\Gamma _i$ is defined as
\begin{equation}
\label{eq:Slocal}
S_i = A_{i}^{\Gamma\Gamma} - A_{i}^{\Gamma I}\left(A_{i}^{II}\right)^{-1}A_{i}^{I\Gamma}.
\end{equation}

The \emph{global Schur complement} $S$ can be assembled from the local ones, similarly to (\ref{eq:subassembly}), as
\begin{equation}
\label{eq:Sassembly}
S = \sum_{i=1}^{N_S} R_i^{\Gamma T} S_{i} R^{\Gamma}_i.
\end{equation}
In a similar way, the reduced right-hand side vector is constructed as
\begin{equation}
g = \sum_{i=1}^{N_S} R_i^{\Gamma T} g_{i}, \quad 
g_i = f_i^{\Gamma} - A_{i}^{\Gamma I}\left(A_{i}^{II}\right)^{-1} f_{i}^{I}.
\end{equation}
Then the reduced global interface problem
\begin{equation}
\label{eq:Sug}
S u^{\Gamma} = g
\end{equation}
is to be solved.
Once we know the solution at the interface $u^{\Gamma}$, solution in the interior of each subdomain can be recovered from the discrete Dirichlet problems
\begin{equation}
\label{eq:discrete_dirichlet}
A_{i}^{II} u_i^{I} = f_{i}^{I} - A_{i}^{I\Gamma}R_i^{\Gamma} u^{\Gamma},\quad i=1,\dots,N_{S}.
\end{equation}

If we use for solving the reduced problem (\ref{eq:Sug}) a preconditioned Krylov subspace method, such as the Preconditioned Conjugate Gradient (PCG) method,
we need to solve the system
\begin{equation}
\label{eq:MSuMg}
M^{-1}S u^{\Gamma} = M^{-1}g,
\end{equation}
and only actions of the matrix $S$ and of the preconditioner $M^{-1}$ are needed (see e.g. \cite{Elman-2005-FEF}).
Consequently, one can circumvent the explicit construction of both $S$ and $M^{-1}$ by providing computational procedures for multiplying a vector by these matrices.
In a procedure implementing the action of $S$, the matrix is applied subdomain-wise according to (\ref{eq:Sassembly}) and (\ref{eq:Slocal}).
The inverse in (\ref{eq:Slocal}) is applied through solving a local discrete Dirichlet problem with the same matrix as in 
(\ref{eq:discrete_dirichlet}) in each iteration.
The key observation for the success of these methods is now the fact that the discrete Dirichlet problems are independent and they can be solved in an embarrassingly parallel way.

\subsection{BDDC preconditioner}
\label{sec:BDDCpreconditioner}

Apart of constructing the Schur complement $S$,
the division into subdomains is used also as the basis for constructing the preconditioner $M^{-1}$.
As is commonly done, the preconditioner will be provided by one step of a nonoverlapping domain decomposition method based on primal unknowns, i.e.\ a subset of $u$.
These evolved from the \emph{Neumann--Neumann} (NN) method~\cite{DeRoeck-1991-ATL} without an explicit coarse space,
through the \emph{Balancing Domain Decomposition} (BDD)~\cite{Mandel-1993-BDD} with the coarse space built from the exact nullspaces of local matrices,
to the BDDC method~\cite{Dohrmann-2003-PSC} allowing more general coarse spaces. 

Perhaps the main difference among these methods is in their approach to solving the local \emph{discrete Neumann problems}. 
These problems are singular for the so called \emph{floating} subdomains, 
\hl{i.e.\ those for which the local solution is not fixed by the boundary conditions of the original problem,}
in the case of Poisson equation as well as the linear elasticity.

In the BDDC method, local problems are regularized by constraints enforcing continuity of the finite element functions across the interface.
Note that these constraints are different from constraints on continuity of shape functions at hanging nodes discussed in Section~\ref{sec:global_contraints},
though similar in spirit.
Prescribing locally these constraints as homogeneous gives rise to a natural coarse problem, which is $A_i$-orthogonal to the local spaces. 
In this way, subdomain problems again become independent and embarrassingly parallel, 
while all the intersubdomain communication is extracted to the (global) coarse problem.

More specifically, from minimizing local energy functionals subject to the continuity constraints, each local problem on subdomain $\Omega _{i}$ leads to the saddle-point system matrix of the type
\begin{equation}
\label{eq:Neumann-matrix}
\left[
\begin{array}[c]{cc}
A_{i} & C_{i}^{T} \\
C_{i} &  0  
\end{array}
\right].
\end{equation}
Here $C_{i}$ is the matrix defining coarse degrees of freedom. 
\hl{In our case, these are arithmetic averages over edges and faces of subdomains for all cases presented in Section~\ref{sec:Numerical}. 
In addition, values at coarse nodes selected by the algorithm from~\cite{Sistek-2012-FSC} are considered for the linear elasticity problems.}

The use of matrix (\ref{eq:Neumann-matrix}) in BDDC is two-fold:
\begin{enumerate}
\item Constructing local \emph{coarse basis functions} $\Phi_{i}$ and the \emph{local coarse matrix} $S_{Ci} = -\Lambda_{i}$ in the set-up of the BDDC preconditioner 
by solving the following problem with multiple right-hand sides
\begin{equation}
\label{eq:Neumann-setup}
\left[
\begin{array}[c]{cc}
A_{i} & C_{i}^{T} \\
C_{i} &  0  
\end{array}
\right]  
\left[
\begin{array}[c]{c}
\Phi_{i}\\
\Lambda_{i}
\end{array}
\right] = \left[
\begin{array}[c]{c}
0\\
I
\end{array}
\right]  ,\quad i=1,\dots,N_{S}.
\end{equation}
Here $I$ is the identity matrix of the correct dimension.

\item Finding the \emph{local subdomain correction} $u_{i}$ from the local residual $r_{i}$ within each application of the BDDC preconditioner
\begin{equation}
\label{eq:Neumann-application}
\left[
\begin{array}[c]{cc}
A_{i} & C_{i}^{T} \\
C_{i} &  0  
\end{array}
\right]  \left[
\begin{array}[c]{c}
u_{i}\\
\mu_{i}
\end{array}
\right] = \left[
\begin{array}[c]{c}
r_{i}\\
0
\end{array}
\right]  ,\quad i=1,\dots,N_{S}.
\end{equation}
\end{enumerate}

The coarse problem matrix $S_C$ can be assembled from the local coarse matrices $S_{Ci}$ defined as
\begin{equation}
S_{Ci} = \Phi_{i}^{T} A_{i} \Phi_{i} = -\Lambda_{i}
\end{equation}
in a way, once again, analogous to (\ref{eq:subassembly}), 
\begin{equation}
\label{eq:coarse_assembly}
S_{C} = \sum_{i=1}^{N_{S}} R_{Ci}^T S_{Ci} R_{Ci},
\end{equation}
where $R_{Ci}$ is the restriction of the global vector of coarse unknowns to those present at the $i$-th subdomain.
However, we keep the matrix $S_{C}$ unassembled and distributed in the parallel sparse direct solver \codename{MUMPS}~\cite{Amestoy-2000-MPD} employed in \codename{BDDCML} for its factorization.

To avoid confusion here, it should be emphasized that in BDDC, the coarse problem is built algebraically in the solver based on the structure of subdomains. 
Consequently, it is not related to the \emph{coarse mesh} of the initial steps of the adaptive finite element algorithm.

Finally, let us briefly recall the complete action of the BDDC preconditioner, largely following notation of~\cite{Dohrmann-2003-PSC}.
Before doing that, let $R_{Bi}:\Omega_{i} \rightarrow \Gamma_{i}$ stand for the 0--1 matrix restricting local vector of unknowns to those at the interface.
Starting from a residual of the Krylov subspace method $r^{\Gamma}$, the preconditioned residual is obtained as $z^{\Gamma} = M_{BDDC}^{-1}\,r^{\Gamma}$ by the following algorithm:
\begin{enumerate}
\item Prepare the subdomain residual on the whole subdomain (both interior and interface unknowns)
\begin{equation} 
\label{eq:making_local_residual}
r_{i} = R_{Bi}^T W_i R_{i}^{\Gamma} r^{\Gamma}.
\end{equation} 
The diagonal matrix $W_i$ here applies weights to satisfy the partition of unity, 
and it is constructed in this contribution as the inverse of the number of subdomains containing an interface unknown. 
See~\cite{Certikova-2015-DAI} for more general options.
Note that $R_{Bi}^T$ only extends the residual at $\Gamma _i$ by zeros in subdomain interior.

\item Solve the subdomain corrections 
\begin{equation} 
\left[
\begin{array}[c]{cc}
A_{i} & C_{i}^{T} \\
C_{i} &  0  
\end{array}
\right]  \left[
\begin{array}[c]{c}
u_{i}\\
\mu_{i}
\end{array}
\right] = \left[
\begin{array}[c]{c}
r_{i}\\
0
\end{array}
\right], \quad i=1,\dots,N_{S}.
\end{equation} 

\item Prepare the coarse residual
\begin{equation} 
r_C = \sum_{i = 1}^{N_S} R_{Ci}^{T} \Phi_{i}^{T} r_{i}.
\end{equation} 

\item Solve the coarse problem
\begin{equation} 
\label{eq:coarse_problem}
S_{C} u_{C} = r_{C}.
\end{equation} 

\item Prolong the coarse correction to every subdomain
\begin{equation} 
u_{Ci} = \Phi_{i} R_{Ci} u_{C}.
\end{equation} 

\item Combine the coarse correction and the subdomain ones
\begin{equation} 
\label{eq:averaging_subdomain_solves}
z^{\Gamma} = \sum_{i=1}^{N_S} R_{i}^{\Gamma T} W_{i} R_{Bi} \left( u_{i} + u_{Ci} \right).
\end{equation} 
\end{enumerate}

As noted in \cite{Dohrmann-2003-PSC,Sousedik-2013-AMB}, the parallel solution of the coarse problem (\ref{eq:coarse_problem}) of the previous algorithm 
eventually becomes a bottleneck for the scalability of the parallel process despite the fact that it is solved by a parallel sparse direct solver \codename{MUMPS}.
To improve the situation, \codename{BDDCML} makes use of the \emph{Multilevel BDDC} method.
The idea of the 3-level method is quite simple:
instead of solving problem (\ref{eq:coarse_problem}) exactly by a direct method, we just use its approximation by one step of the BDDC method applied to this problem. 
In particular, due to the analogy of (\ref{eq:coarse_assembly}) with (\ref{eq:global_assembly}),
the coarse problem is seen as a small finite element problem in which subdomains play the role of elements. 
These are grouped into higher-level subdomains, and the preconditioner proceeds similarly to the standard (2-level) algorithm.
It is now clear that the idea can be repeated recursively, giving rise to the multilevel method.
\hl{The coarse degrees of freedom used on higher levels are analogous to those on the first level.}

Multilevel BDDC introduces error with each additional level because a preconditioner is used in place of an exact coarse solve. 
Thus, the preconditioner becomes less accurate, and the number of iterations typically grows with the number of levels. 
See also \cite{Mandel-2008-MMB} for the theoretical proof. 
However, the resulting coarse problem is much smaller and simpler for solution.
As a result, the algorithm is more scalable in a massively parallel setting \cite{Badia-2016-MBD,Sousedik-2013-AMB}.

\subsection{Accommodating disconnected subdomains in BDDC}
\label{sec:handling_disconnected}

In the BDDC theory \cite{Mandel-2003-CBD,Mandel-2005-ATP,Mandel-2008-MMB}, 
an important role is played by the following assumption. 
\begin{assumption}
\label{as:nullA_i}
Let us assume that
\begin{equation}
\mathrm{null}\,A_{i} \cap \mathrm{null}\,C_{i} = \{{\mathbf 0}\}.
\end{equation}
\end{assumption}
This assumption represents the necessary and sufficient condition for a unique solvability of systems (\ref{eq:Neumann-setup}) and (\ref{eq:Neumann-application})
with a positive semidefinite matrix $A_{i}$~\cite[Theorem 3.2]{Benzi-2005-NSS}. 
In the context of BDDC, Assumption~\ref{as:nullA_i} is satisfied if enough constraints are selected on each subdomain.

Let us now discuss the case in which a subdomain $ \Omega_{i}$ is composed of several disconnected parts $\Omega_{i}^{c,k}$, $k=1,\dots, n^{c}_{i}$, 
where $n^{c}_{i}$ is the number of such \emph{components}.
The matrix $A_{i}$ becomes block-diagonal with $n^{c}_{i}$ blocks.
Consequently, each block $A_{i}^{k}$ potentially has its own nullspace. 
In order to satisfy Assumption~\ref{as:nullA_i} in this case, it is natural to select constraints $C_{i}^{k}$
independently for each subdomain component and to compose the matrix of constraints $C_{i}$ based on these blocks. 
The structure of the matrix (\ref{eq:Neumann-matrix}) is then refined to
\begin{equation}
\label{eq:Neumann-matrix-blocks}
\left[
\begin{matrix}
A_{i}^{1} &           &            &               &  C_{i}^{1,T} &             &        &                 \\ 
          & A_{i}^{2} &            &               &              & C_{i}^{2,T} &        &                 \\
          &           & \ddots     &               &              &             & \ddots &                 \\
          &           &            & A_{i}^{n^c_i} &              &             &        & C_{i}^{n^c_i,T} \\
C_{i}^{1} &           &            &               &              &             &        &                 \\
          & C_{i}^{2} &            &               &              &             &        &                 \\
          &           & \ddots     &               &              &             &        &                 \\
          &           &            & C_{i}^{n^c_i} &              &             &        &                 \\
\end{matrix}
\right]
\end{equation}

If enough constraints are selected for each component $\Omega_{i}^{c,k}$, $\mathrm{null}\,A_{i}^k \cap \mathrm{null}\,C_{i}^{k} = \{{\mathbf 0}\}$ for each $k$.
Now let $R_i^{k}:\Omega_{i} \rightarrow \Omega_{i}^{c,k}$ be the 0--1 matrix selecting the $k$-th component of $\Omega_{i}$.
Then
$\mathrm{null}\,A_{i} = \mbox{span}\left\{(R_i^{k})^T\mathrm{null}\,A_{i}^k\right\},\ k=1,\dots, n^{c}_{i}$ and $\mathrm{null}\,C_{i} = \mbox{span}\left\{(R_i^{k})^T\mathrm{null}\,C_{i}^k\right\},\ k=1,\dots, n^{c}_{i}$. 
Here $(R_i^{k})^T$ extends the component values to the whole subdomain by zeros.
It follows that $\mathrm{null}\,A_{i} \cap \mathrm{null}\,C_{i} = \{{\mathbf 0}\}$,
and Assumption~\ref{as:nullA_i} is satisfied. 

This gives a complete picture for problems with one-dimensional nullspace of $A_{i}$ for floating subdomains 
(i.e.\ those without Dirichlet boundary conditions) such as the Poisson problem.
However another situation can still arise for problems with higher-dimensional nullspaces, such as in the case of linear elasticity, which has three and six-dimensional nullspaces for floating subdomains in two and three dimensions, respectively.
In this case, two components can be also \emph{loosely coupled}, by which we mean that there is not enough coupling between the two components to make the nullspace as simple as in the case of a connected subdomain. 
As an example for linear elasticity in 2D, this case arises if two mesh components share only one node. 
Then the two components are not mutually fixed and can rotate relatively to each other around this common node.
In general, for two components of a floating subdomain, the dimension of the nullspace is between that of the simple floating subdomain, and twice this number
for two completely disconnected mesh components.  

Our approach to resolve this issue is in a proper definition and detection of components of the subdomain mesh.
To this end, we detect components of the \emph{dual graph of the mesh}, in which mesh elements become graph vertices and an edge between two vertices is present whenever the two elements share an appropriate number of nodes depending on the type of finite elements. 
Returning once again to the case of linear elasticity, this corresponds to sharing an element edge in 2D and an element face in 3D.
Then, if a loose coupling appears within the subdomain mesh, more graph components are detected and constraints are generated independently for each of them as in the case of completely disconnected mesh parts.

Let us now look at the structure of the matrix (\ref{eq:Neumann-matrix-blocks}) from the point of view of a direct solver applied to (\ref{eq:Neumann-setup}).
\hl{The main goal of the employed load balancing is to approximately equal the number of finite elements within each subdomain.} 
If the subdomain is composed of several components, either disconnected or loosely coupled, 
the structure of the leading block $A_{i}$ simplifies and can even become block diagonal. 
For completely independent mesh parts, the system (\ref{eq:Neumann-setup}) in fact splits into several independent saddle-point systems,
although this is not explicitly utilized in our solver.
However, due to the independent detection of constraints for each component, 
number of rows of $C_{i}$ corresponding to the number of coarse degrees of freedom of the subdomain can increase significantly.
Consequently, problem (\ref{eq:Neumann-setup}) is solved for a larger number of right-hand sides and the dimension of the global coarse problem also increases for a disconnected subdomain.
These are presumably two somewhat contradictory effects on performance of a direct solver applied to (\ref{eq:Neumann-setup}),
creating some room for load imbalance. 
Consequently, the impact of disconnected subdomains on an actual performance of the parallel BDDC solver is investigated in Section~\ref{sec:Numerical}.

\section{Coupling parallel AMR with BDDC}
\label{sec:Coupling}
In the previous two sections, we have given an overview of the two main ingredients of the calculation. The first is discretisation of
the problem on meshes with hanging nodes, the second is the solver applied to the resulting linear algebraic system. 
In this section, we will further refine these concepts in regard to their combination, and the adaptivity in a massively parallel setting, the ultimate aim of the paper.

\subsection{Parallel mesh handler}
\label{sec:p4est}

Although very sophisticated strategies for mesh partitioning exist, such as equilibrating size of factors of local problems (see e.g. \cite{Jurkova-2010-FBG}), 
these are mostly related to balanced load of the algebraic solver itself. 
In the broader context of a finite element simulation, other parts of the algorithm may dominate the solution. 
For example, integration of element matrices and subassembly of subdomain matrices by far dominated the computational 
time in time-dependent incompressible flow simulations in \cite{Sistek-2015-PIS}.
This trend is also pronounced for recent advanced discretisation techniques such as high order elements with Lagrange or B-spline basis functions.
These parts of the code are embarrassingly parallel, with the amount of work proportional to the local number of elements.
Therefore, keeping the numbers of elements balanced across subdomains is seen as the right requirement for parallel simulations at very large scale.

\begin{figure}[tbh]
\centering
 \includegraphics[height=0.32\textwidth, angle=90]{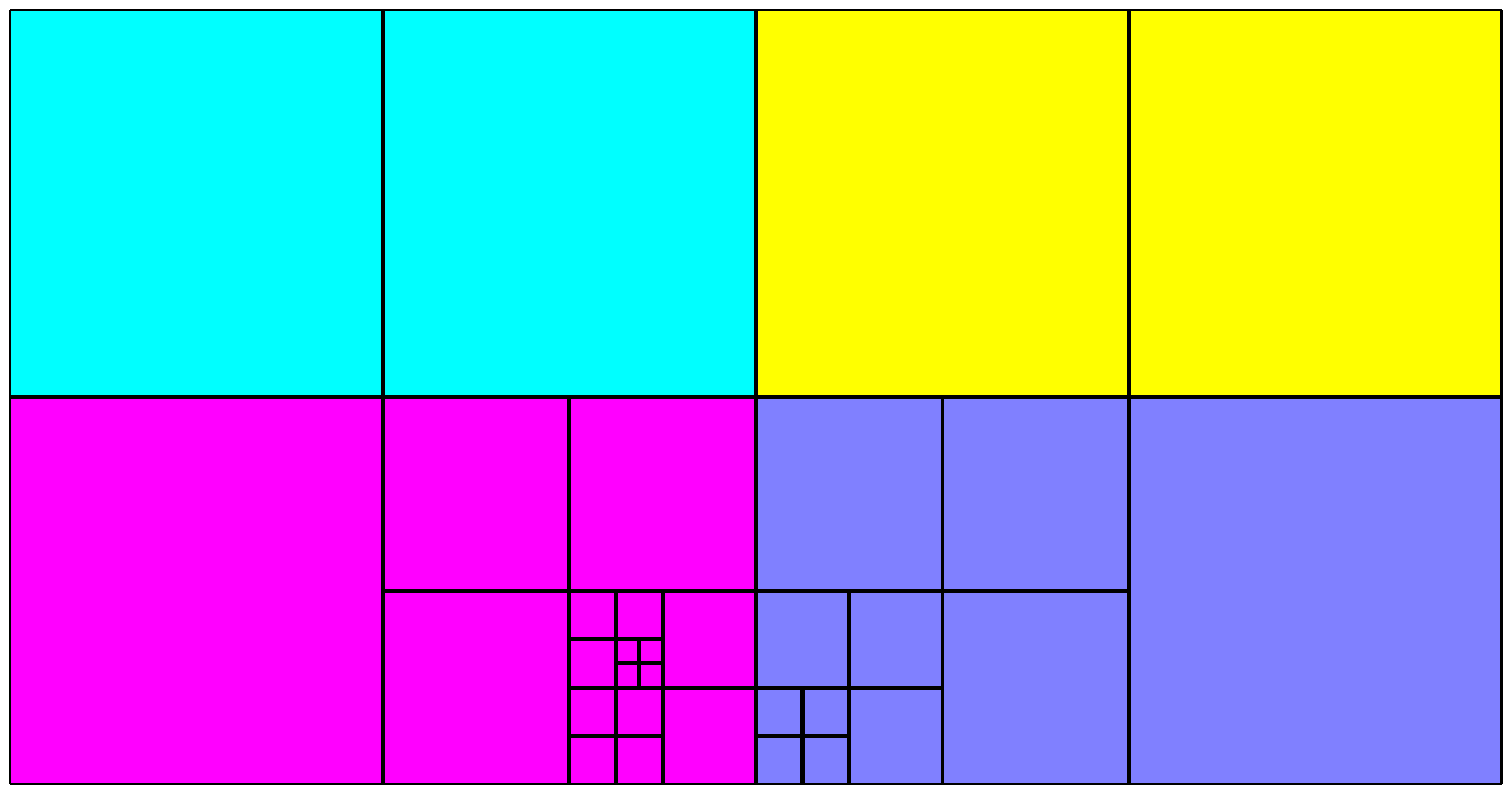} 
 \includegraphics[height=0.32\textwidth, angle=90]{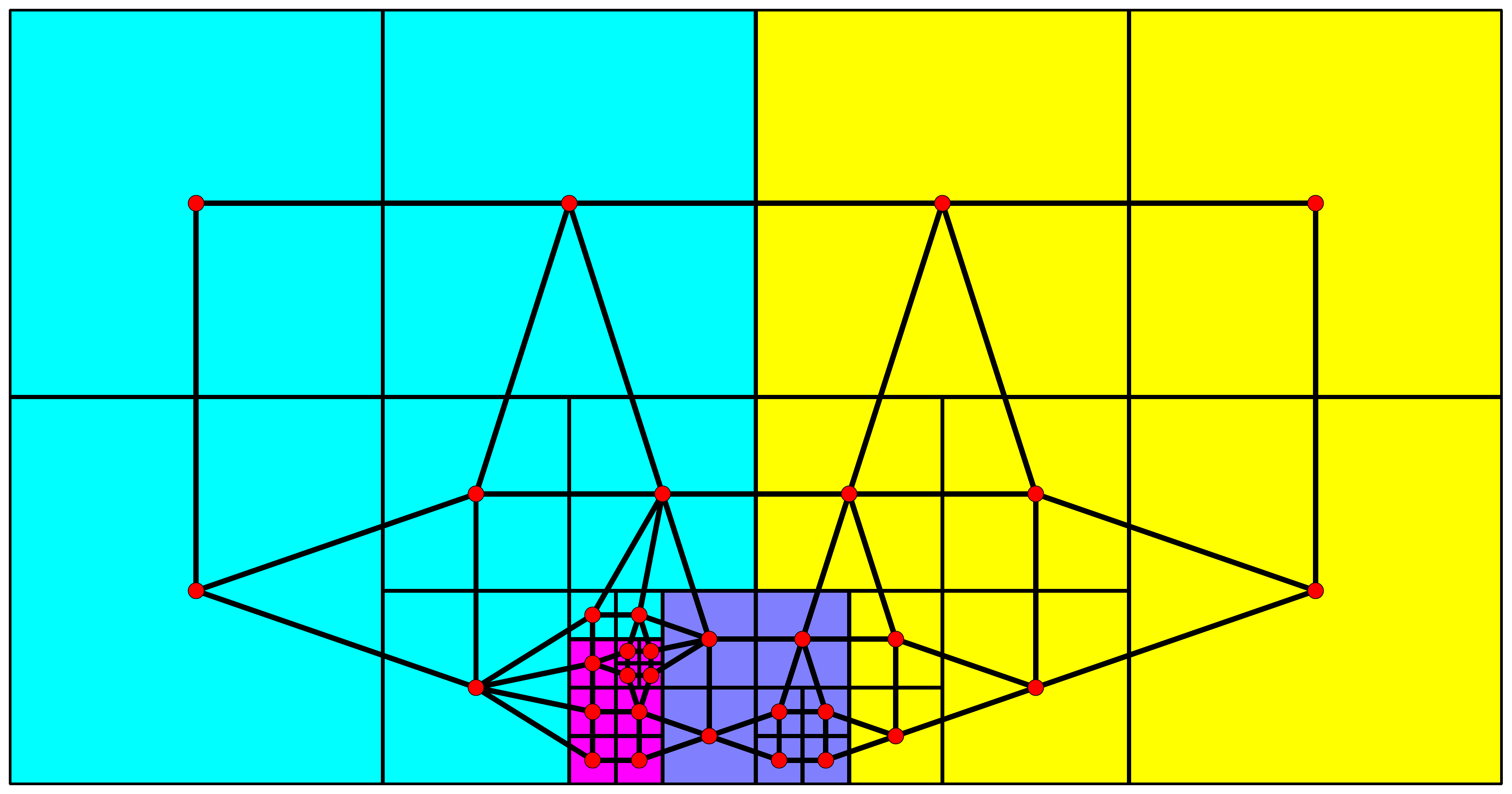} 
 \includegraphics[height=0.32\textwidth, angle=90]{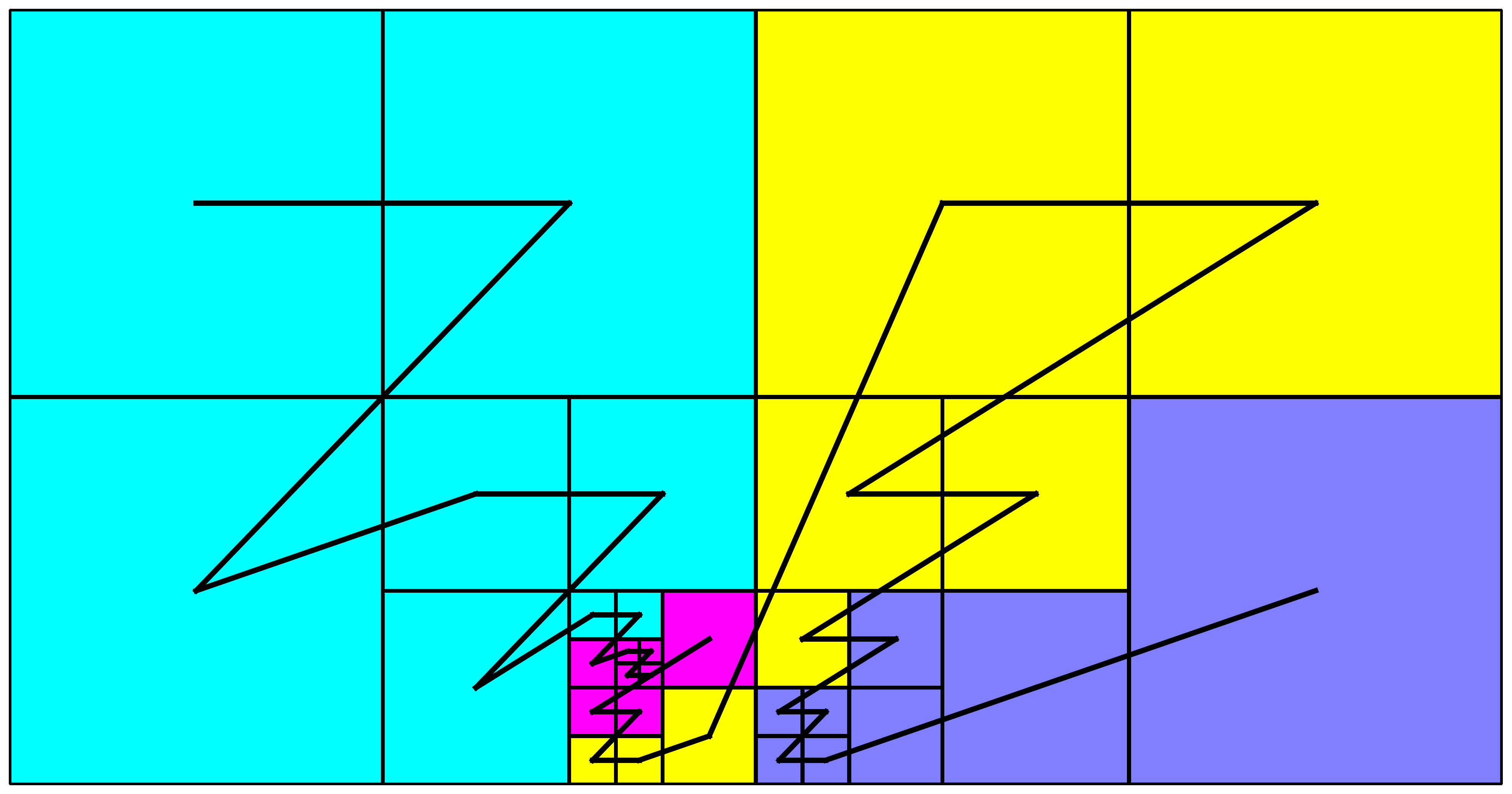} 
 \caption{Unbalanced mesh (left), applying graph partitioner (centre) and applying a space-filling curve (right). Unbalanced mesh is clearly not usable. Graph partitioning
 usually results in better shape of subdomains, but scales only to moderate number of subdomains.}
\label{fig:balancing}
\end{figure}

The main reason for considering adaptive mesh refinement is an effort to invest most of the computational resources into 
troublesome parts of the domain, such as close to singularities or within boundary layers. Since we usually do not know these 
areas a priori, we have to perform several refinement steps, where the mesh is adaptively refined towards such areas. 
It is therefore obvious that in the parallel setting one cannot use adaptivity without some kind of mesh re-balancing between these steps. 
If no balancing was used, adaptivity would lead to meshes such as the example shown in Fig.~\ref{fig:balancing} (left), where 
elements of only few subdomains are refined. 
This leads to complete imbalance of subdomain sizes after several adaptivity steps. 
In such setting, very little benefit would be taken from the employment of domain decomposition. 

A possible solution would be the use of a graph partitioner, such as \codename{METIS}~\cite{Karypis:1998-FHQ}, after each adaptivity step 
in order to create new balanced subdomains. Although one cannot use mesh partitioning directly because of the presence of hanging nodes, 
an incidence graph (such as in the centre of Fig.~\ref{fig:balancing}) can be created and partitioned. This might be a very reasonable approach
for a moderate number of subdomains. 
Although the shape of these subdomains can be far from optimal (see e.g.~\cite{Klawonn-2008-AFA}), 
it is often acceptable. Also, connected subdomains can be enforced.
However, the serial graph partitioning cannot be used for a large number of subdomains within a parallel computation.  
Moreover, parallel versions (such as \codename{ParMETIS}~\cite{Karypis-1999-PMS}) 
also have their own limitations and do not scale to thousands of cores. 

An alternative approach, which we follow in this paper, is the use of space-filling curves, such as a Z-curve, see Fig.~\ref{fig:balancing}, right.
As can be seen, a drawback of Z-curves is that their use for mesh partitioning typically leads to disconnected and loosely coupled subdomains.
However, in connection to the BDDC method, these artifacts can be resolved at the level of the parallel linear solver by selection of appropriate 
constraints for each subdomain component. 
Consequently, the strategy described in Section~\ref{sec:handling_disconnected} can be readily applied to subdomains provided by the \codename{p4est} library.

\subsection{Nonlocal degrees of freedom}

Interface unknowns are located at the geometric boundary of subdomains for many types of finite elements, such as the classical piecewise polynomial Lagrange elements used here. 
However, more general interfaces arise as soon as degrees of freedom shared by adjacent subdomains are not geometrically aligned with element boundaries.
For example, the \emph{fat interface} is naturally encountered for B-spline basis functions~\cite{BeiraoDaVeiga-2013-BPI,BeiraoDaVeiga-2013-IBP}.

\begin{figure}[tbh]
\centering
 \includegraphics[width = 0.5\textwidth]{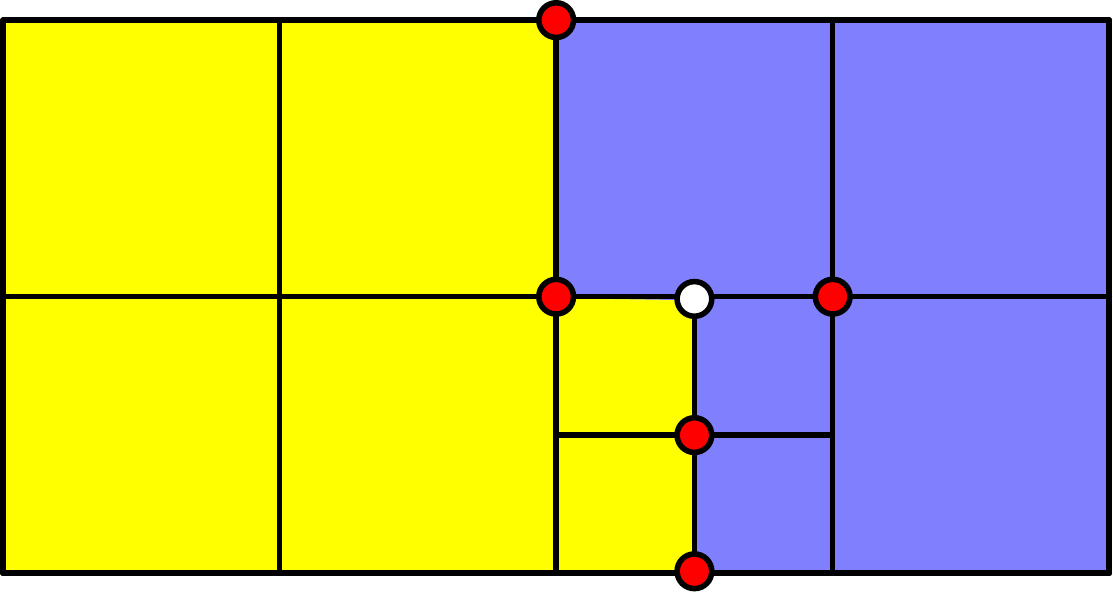}
 \caption{A generalized interface.
The hanging node at the geometric boundary of subdomains is eliminated, resulting in an interface degree of freedom 
proliferating into the interior of the right subdomain.
Interface degrees of freedom are denoted by red circles, the eliminated hanging node is marked by the white circle.
}
\label{fig:generalized_interface}
\end{figure}

A similar generalization is needed also here.
Figure~\ref{fig:generalized_interface} should illustrate the effect.
The situation arises when the division of elements into subdomains suggested by the Z-curves contains a hanging node at the geometric boundary of elements. 
After elimination of this hanging node, described in detail in Section~\ref{sec:global_contraints}, 
a node that is originally interior to one of the subdomains becomes a part of the interface, because its degree of freedom is contained in both subdomains. 
Moreover, the long edge at the hanging node becomes shared by three elements after this elimination.
This complicates the geometry of the interface, but in our experience does not require any modifications of the domain decomposition algorithm itself.

\subsection{Modification of the adaptivity procedure}
\label{sec:modified_refinement}
\hl{
Let us first recall a standard adaptivity algorithm in a parallel setting.
After setting up an initial mesh, the following steps are performed: \\
\vskip 1mm
\noindent Loop until convergence or predefined number of refinements.
\begin{enumerate}
 \item \label{ad:integrate} Construct element and subdomain matrices for a given mesh.
 \item \label{ad:solve} Find the discrete solution.
 \item \label{ad:estimate} Estimate the solution error on each element.
 \item \label{ad:decide} Decide which elements to refine.
 \item \label{ad:refine} Refine the selected elements.
 \item \label{ad:rebalance} Re-balance the mesh. 
\end{enumerate}
End loop.
\vskip 2mm

We have implemented these steps, and the results of the adaptivity runs are presented in Section~\ref{sec:Numerical_adapt}. 
Our detailed performance analysis is, however, restricted to Step~\ref{ad:solve}, which is performed by our BDDC solver.
The other parts of the algorithm benefit from the simplifications of our setting,
and their cost is negligible compared to Step~\ref{ad:solve}.

To be more specific, Step~\ref{ad:integrate} benefits from constant coefficients in the equation and the simple (cubic) shape of the elements.
Steps~\ref{ad:refine} and~\ref{ad:rebalance} are simplified by the fact that there are no specific data (such as solutions from the previous
time levels, etc.) connected to elements, and we do not have to project and transfer any information.
 
Estimating the error in Step~\ref{ad:estimate} is a key part of each adaptivity algorithm.
This estimate can be also used for estimating the global error and terminating the adaptivity procedure.
There are many ways how to estimate the element error. 
In particular, various kinds of a posteriori error estimates are constructed for different problems, 
see e.g. \cite{ainsworth-2000} for an excellent introduction to this topic.
In our calculation, however, we use the knowledge of the exact solution for this step. 

Although Step~\ref{ad:decide} is rather trivial in serial setting, it needs to be elaborated for the case of distributed adaptivity algorithm.}
Let us assume that we are able to assign to each element $\Omega _{K}^{e}$ a real number $\eta_K$, which estimates local error of the solution. 
In each adaptivity step, one should refine some of the elements in such a way, that the global error would decrease as 
fast as possible and converge to zero. Unfortunately, it is rather difficult to prove convergence of adaptive algorithms
in a general way (see e.g. \cite{dorfler_convergent_1996,sebestova-2016} for details). 
In practice, however, 
there are two widely adopted approaches. 
In the simpler one, all elements with an estimated error $\eta_K > \theta\, \eta_{\mbox{max}}$ are refined in each adaptivity step. 
Here $\eta_{\mbox{max}}$ is the maximal element error over the whole mesh and $\theta\in(0,1)$ is a suitable factor. 
If $\theta$ is chosen too small, we might do a lot of unnecessary refinements. 
If it is too large, only few elements would be refined, and we might end up doing prohibitively many
adaptivity steps. Thus, a reasonable choice of $\theta$ is critical. 

\begin{figure}
\centering
 \includegraphics[width = 0.24\textwidth]{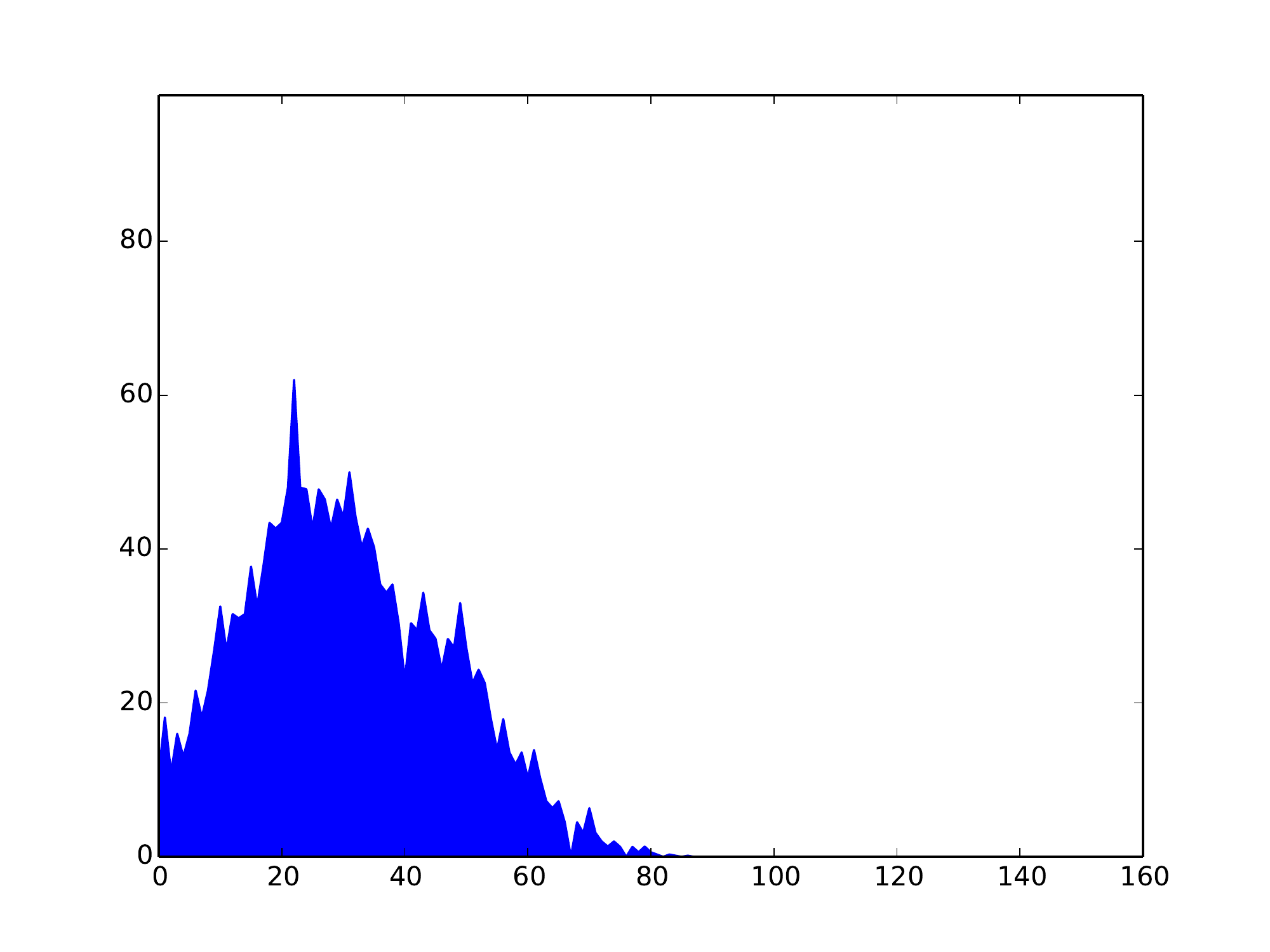}
 \includegraphics[width = 0.24\textwidth]{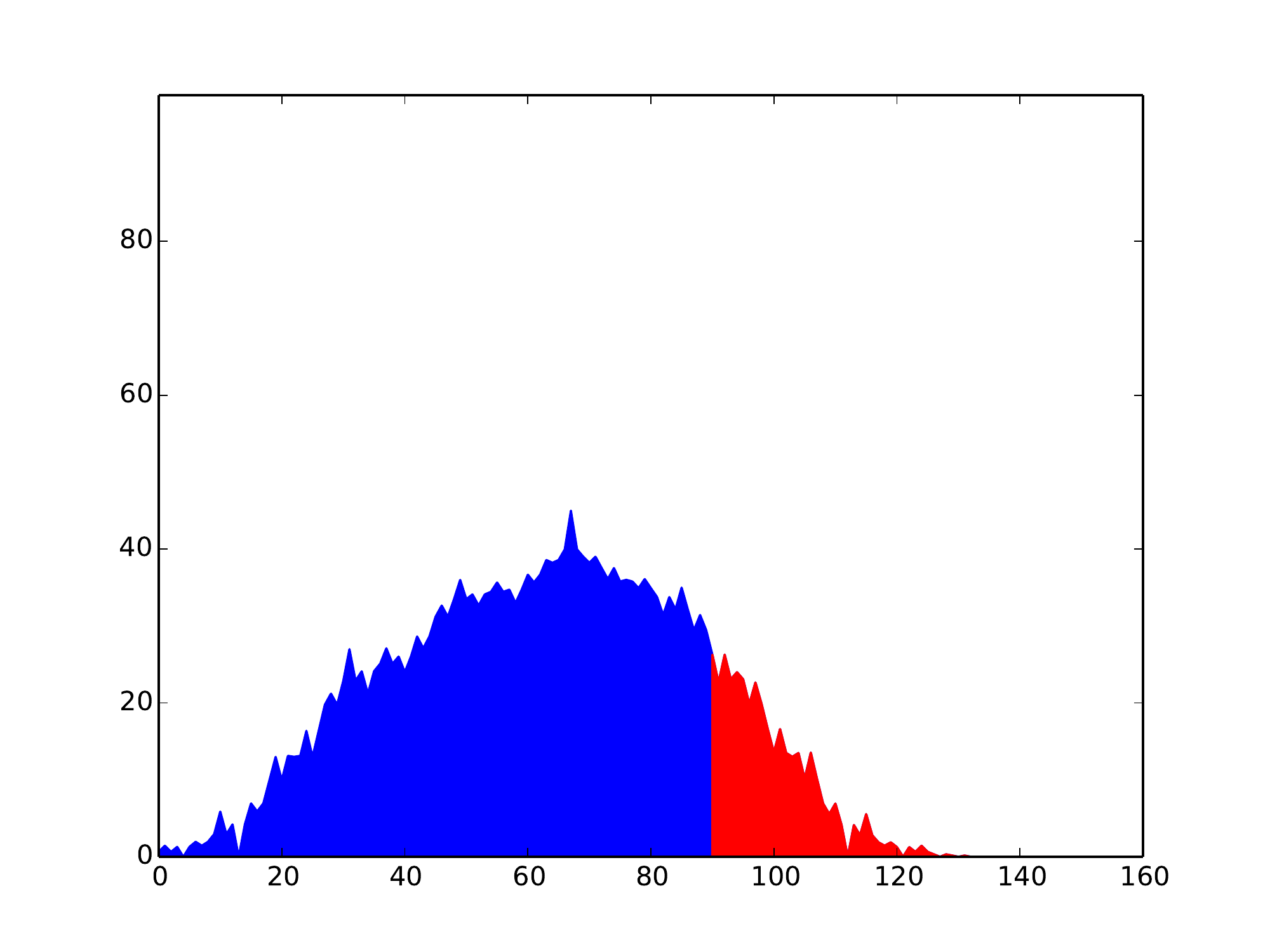}
 \includegraphics[width = 0.24\textwidth]{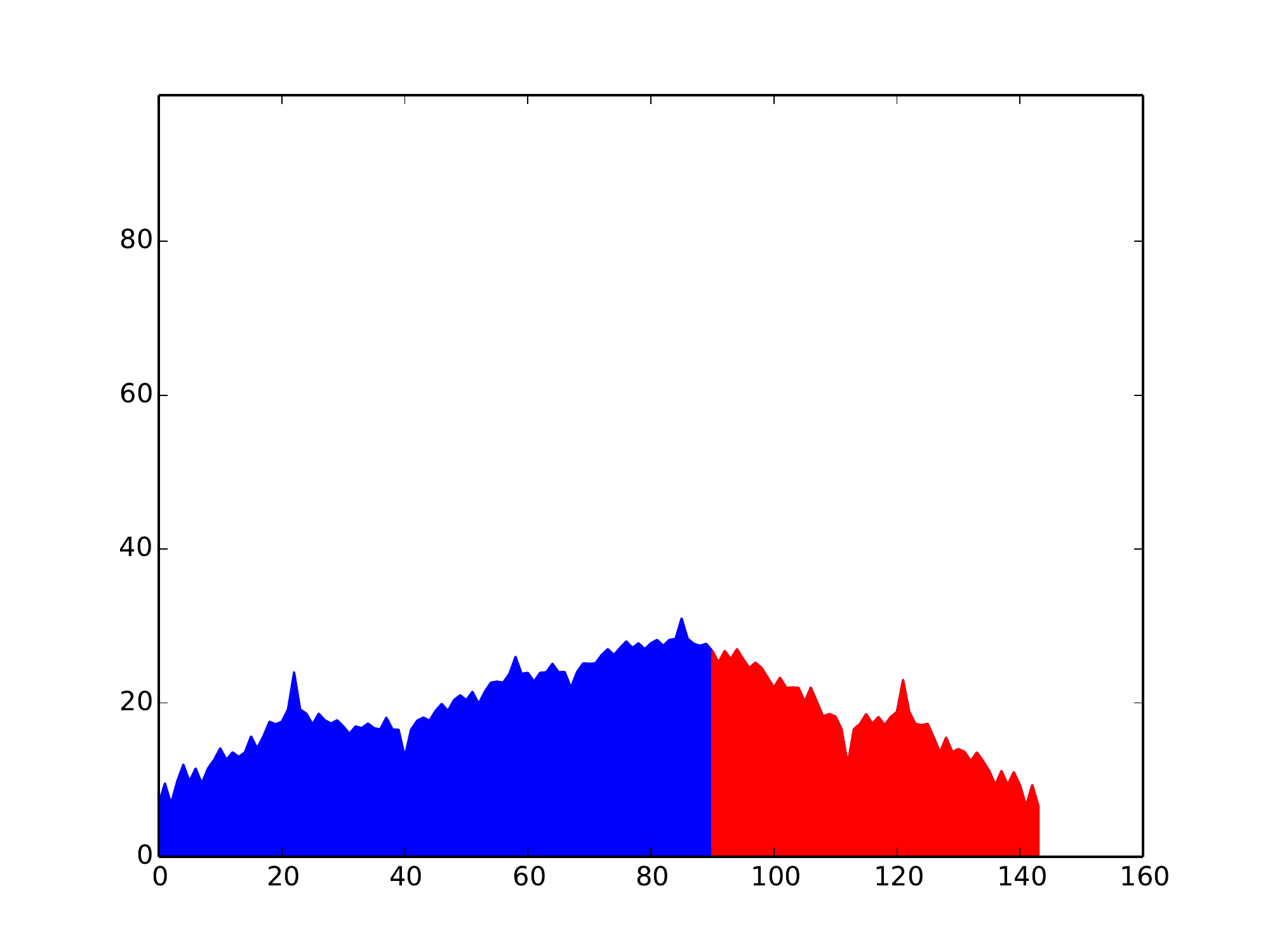}
 \includegraphics[width = 0.24\textwidth]{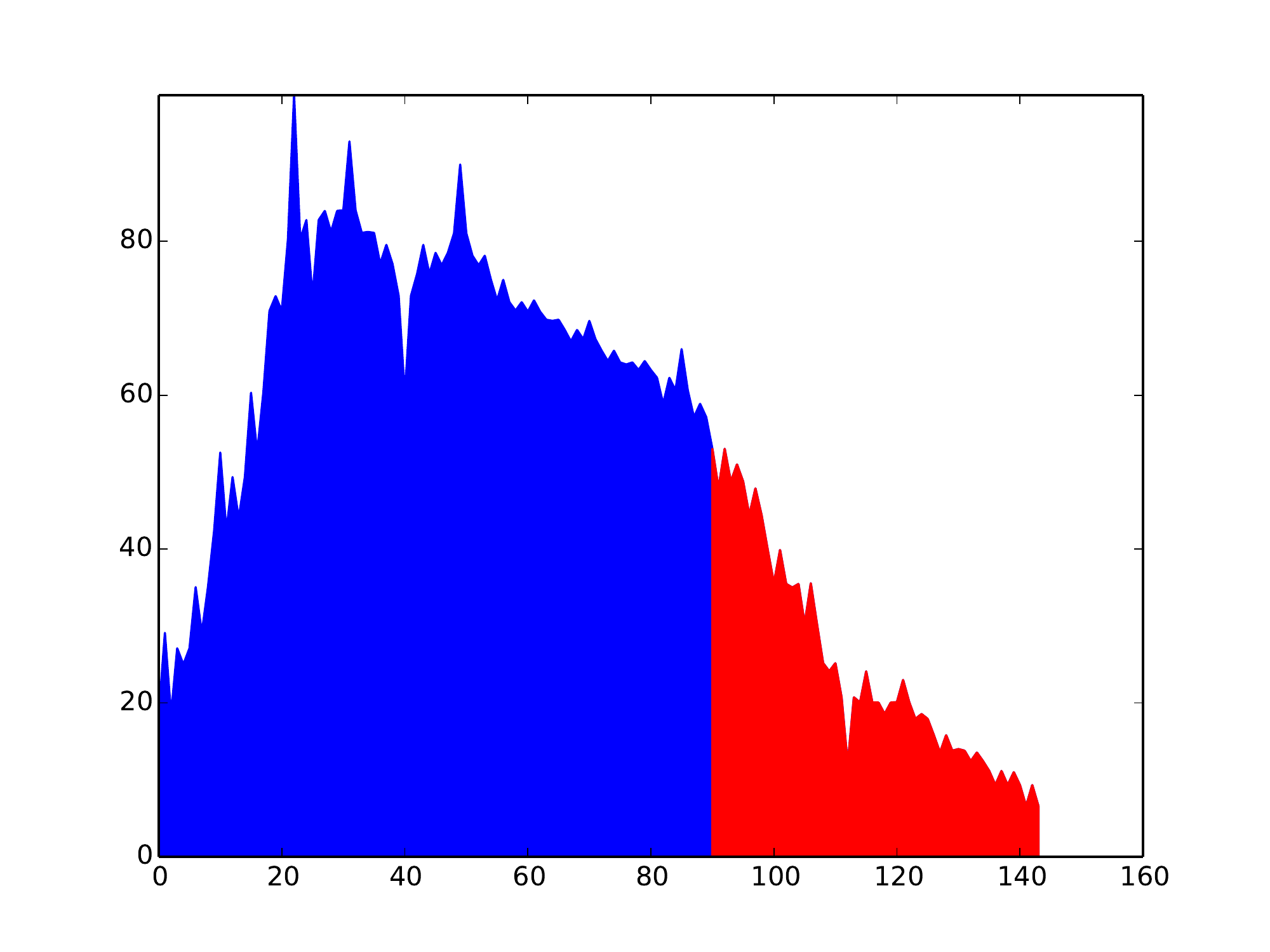}\\
%
%
 \includegraphics[width = 0.24\textwidth]{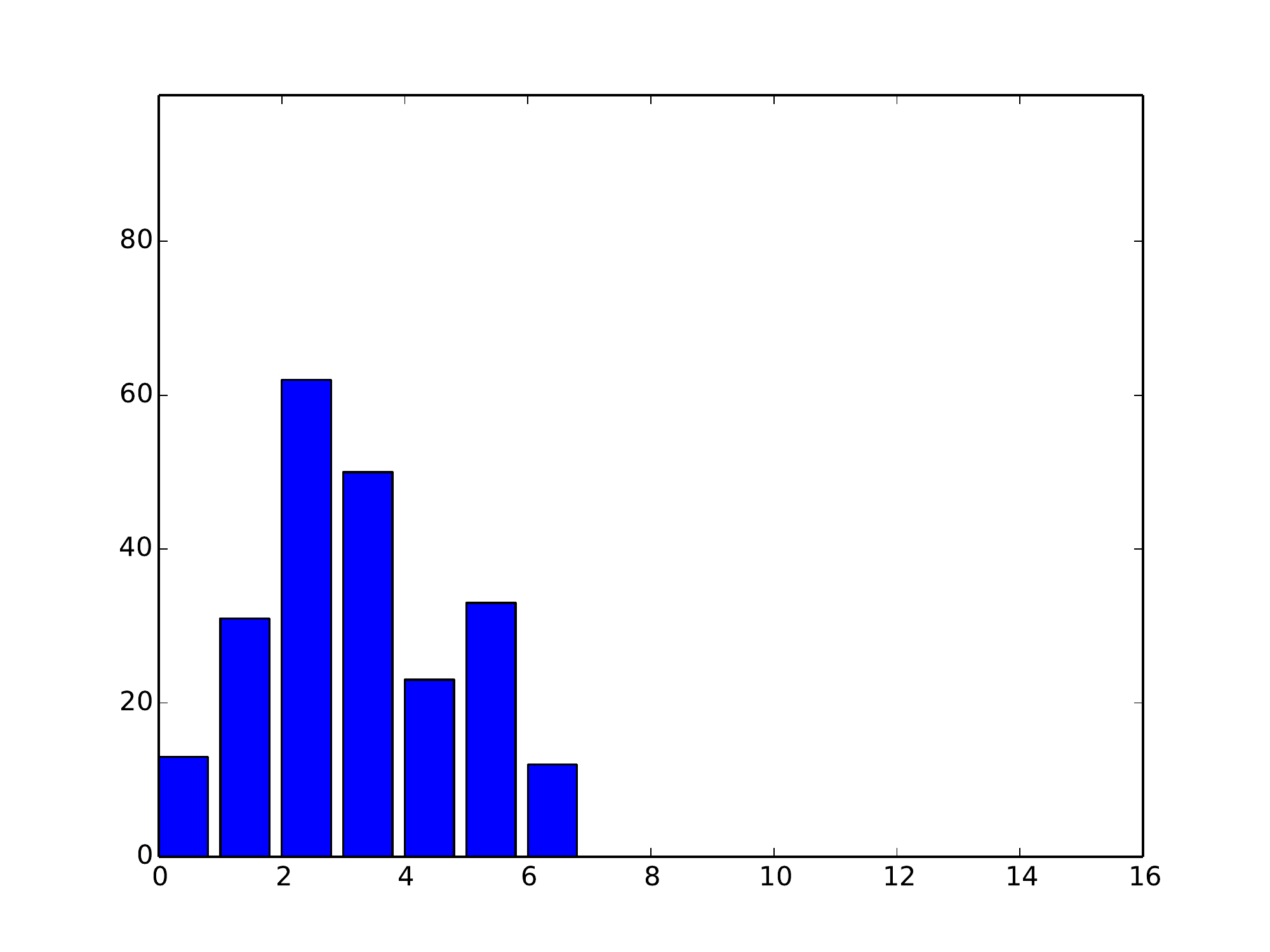}
 \includegraphics[width = 0.24\textwidth]{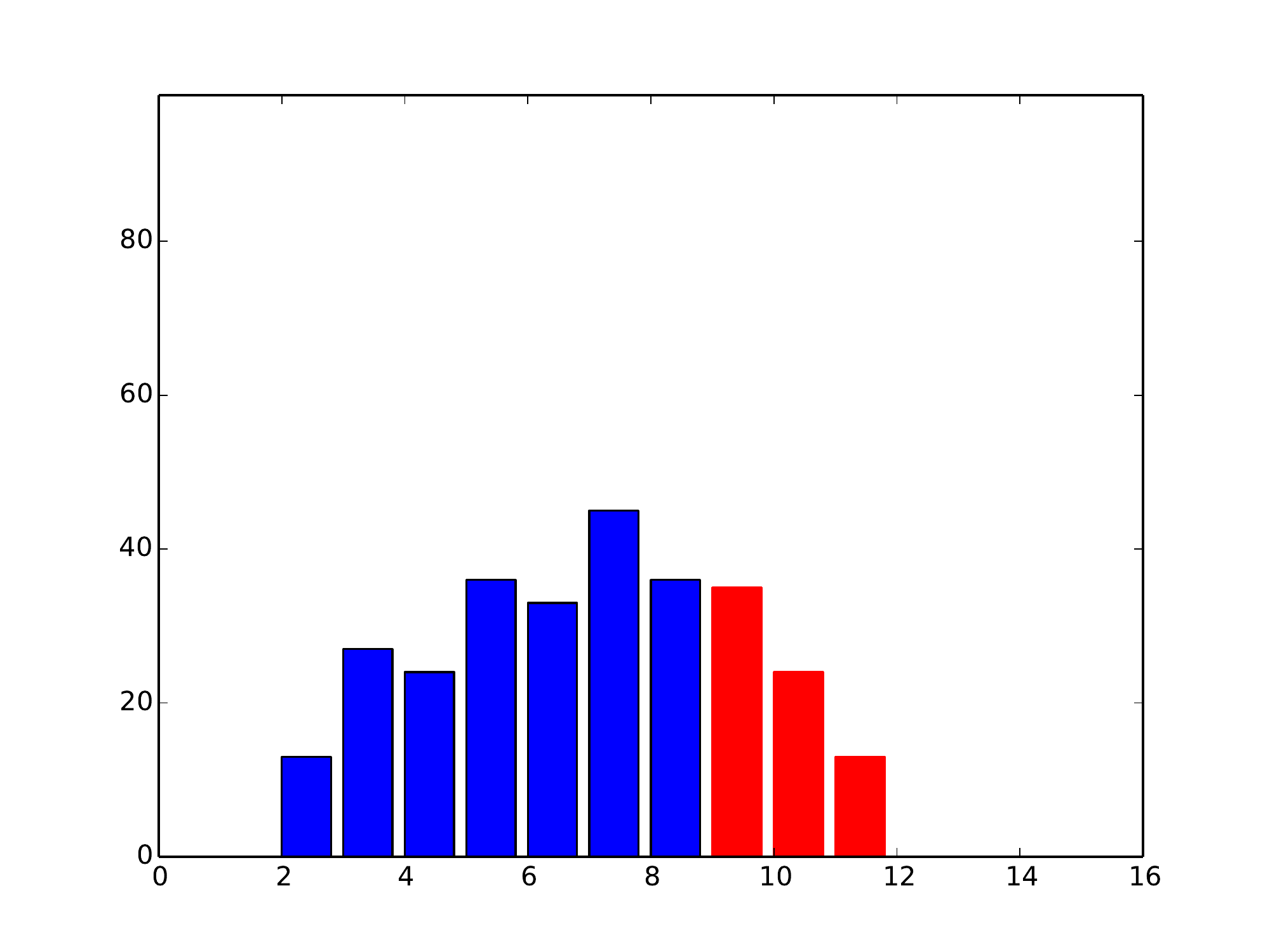}
 \includegraphics[width = 0.24\textwidth]{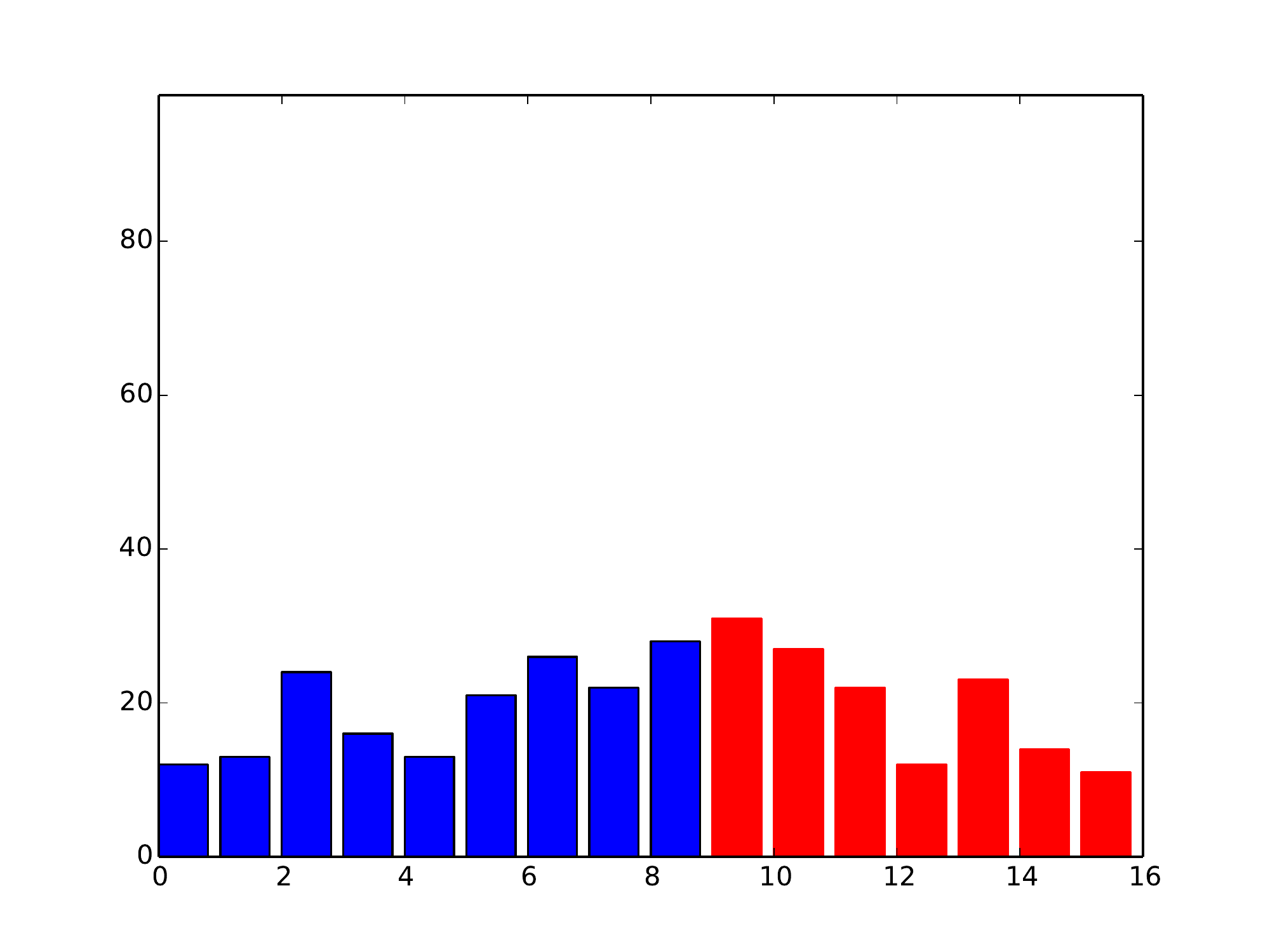}
 \includegraphics[width = 0.24\textwidth]{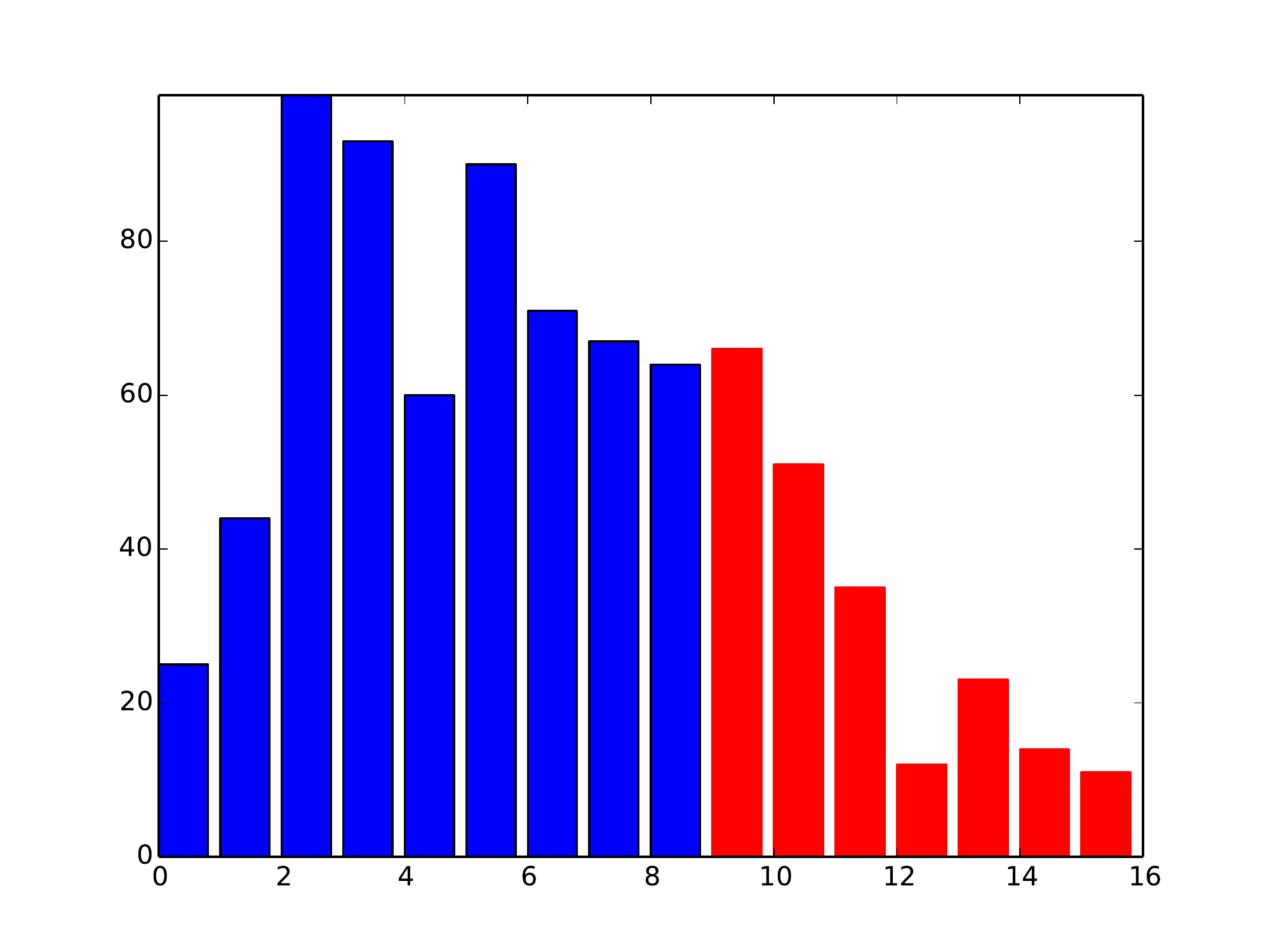}
\caption{Histograms of estimated element error (x axis: element error, y axis: normalised number of elements) on 3 
processors and globally (right). The top line using 160 intervals can be approximated by the bottom line using only 16 intervals.}
 \label{fig:histogram}
\end{figure}

Adapting this algorithm for parallel setting is straightforward. 
We first find a maximal error $\eta_i = \max_{K\in\Omega_i} \eta_K$ for each subdomain $\Omega_i$ and then use a global reduction of a single number
to find the
global maximum $\eta_{\mbox{max}} = \max_{i} \eta_i,\ i = 1,\dots,N_S$.

The disadvantage of this algorithm is the necessity to find a suitable value $\theta$, 
which can be difficult to do in a problem-independent way. 
From the practical point of view, it is often more convenient to refine a given percentage $\zeta$ of 
elements (e.g. 30\%) with the largest estimated error in each adaptivity step. Within a serial calculation, this can be done
easily by sorting elements according to their estimated error, defining $\bar{\theta}$ such that $\zeta N_e$ elements have their error larger than $\bar{\theta}$.
Then, as in the previous approach, all elements with estimated error larger than $\bar{\theta}$ are refined.
Unfortunately, this algorithm is not readily applicable to a parallel computation.

The task of refining a given percentage of elements can be viewed through histograms in Fig.~\ref{fig:histogram}.
Here the $x$-axis corresponds to intervals of the element error $\eta_K$ while the $y$-axis corresponds to numbers of elements in each interval normalized by the interval size.
We can see, that a fine sampling of $\eta_K$ in the top line can be approximated well by a rather coarse sampling at the bottom line. 
This leads us to the following parallel algorithm based on a reasonable number of discrete compartments:
\begin{enumerate}
\item Find the local maximal element error $\eta_i = \max_{K\in\Omega_i} \eta_K$ on each subdomain. 
\item Find the global maximum $\eta_{\mbox{max}}= \max_{i} \eta_i,\ i = 1,\dots,N_S$ by a global reduction of a single scalar. 
\item Define equidistant partitioning of interval $[0, \eta_{\mbox{max}}]$ into moderate number $M$ of subintervals with bounds $0, L, 2L,\ldots, ML = \eta_{\mbox{max}}$, where $L=\eta_{\mbox{max}}/M$. 
\item Define $H_m^i$ as the number of elements with estimated error from interval $\bigl((m-1)L,mL\bigr]$ at each subdomain (i.e.\ process) $\Omega_i$. 
\item Find global sum of number of elements in each interval $H_m = \sum_i H_m^i$. 
We need a global reduction of a vector of length $M$, which is typically not much more expensive than \hl{a} reduction of a single scalar as long as $M$ is small.
\item Find a value $\hat{\theta}$ as the maximal possible value $\bar{m}L$ such that $\sum_{m=\bar{m}}^M H_m \geq \zeta N_e$. 
\end{enumerate}
Now we use $\hat{\theta}$ as the approximation of $\bar{\theta}$.

\section{Numerical results}
\label{sec:Numerical}

We have performed numerical tests for the Poisson problem and linear elasticity in two and three dimensions.
All tests have been run using the \emph{Salomon} supercomputer at the National Supercomputing Centre \emph{IT4Innovations} with up to 4096 cores. 
The largest problem solved (on a highly adapted mesh) had around 1.3$\times$10$^9$ unknowns.

Most of the tests have been performed for the following Poisson problem in 2D and 3D, 
\begin{equation}
-\Delta u = 1\ \mbox{on}\ \Omega,\quad u = 0 \ \mbox{on}\ \partial\Omega,\quad \Omega = [0,1]^d, d\in{2,3}.
\end{equation}

Each node of the \emph{Salomon} supercomputer has 24 cores (two 12-core Intel Xeon E5-1620V3 processors) with 2.5 GHz and 128 GB RAM\@.
The communication network is InfiniBand connected into 7D enhanced hypercube.
In total, the computer has 1008 computing nodes.
\hl{Intel C, C++, and Fortran compilers version 16 and Intel MPI library version 5.1 were consistently used in our study.}

An important ingredient for the performance of the domain decomposition algorithm is a sparse direct solver for the local subdomain problems and the final coarse problem. 
We have used \codename{MUMPS} version 4.10.0 for this purpose.
In particular, the local discrete Dirichlet problems (\ref{eq:discrete_dirichlet}) are solved using the sparse Cholesky ($LL^T$) factorization on each subdomain,
while the discrete Neumann problems (\ref{eq:Neumann-setup}) and (\ref{eq:Neumann-application}) are solved using the $LDL^T$ factorization for symmetric indefinite matrices.
The coarse problem
(\ref{eq:coarse_problem})
is solved on the last level by a parallel Cholesky factorization.

In all presented computations, PCG iterations are stopped based on the relative norm of the residual when ${\|r^{\Gamma}\|_2}/{\|g\|_2} < 10^{-6}$.

\subsection{Weak and strong scalability}
\label{sec:numerical_scaling}

There are many combinations of relevant parameters, for which scaling results have been performed: 
2-D and 3-D problems, 
two and three levels in the BDDC method, 
the first and higher order of FEM approximation, 
selected refinements that determine not only the final size of the problem (number of degrees of freedom)
but also the structure of hanging nodes, and finally the effect of solving the linear elasticity problem instead of the Poisson equation.
Out of these combinations, we have selected those allowing us to best assess the performance of the method. 
We have organized them into four groups, which will be addressed in detail in the following.
A particular interest is devoted to the imbalance in solving the local Neumann problems 
due to the variable number of local constraints for disconnected subdomains.

\subsubsection{Regular subdomains}

\begin{figure}
\begin{center}
\includegraphics[height = 0.33\textwidth]{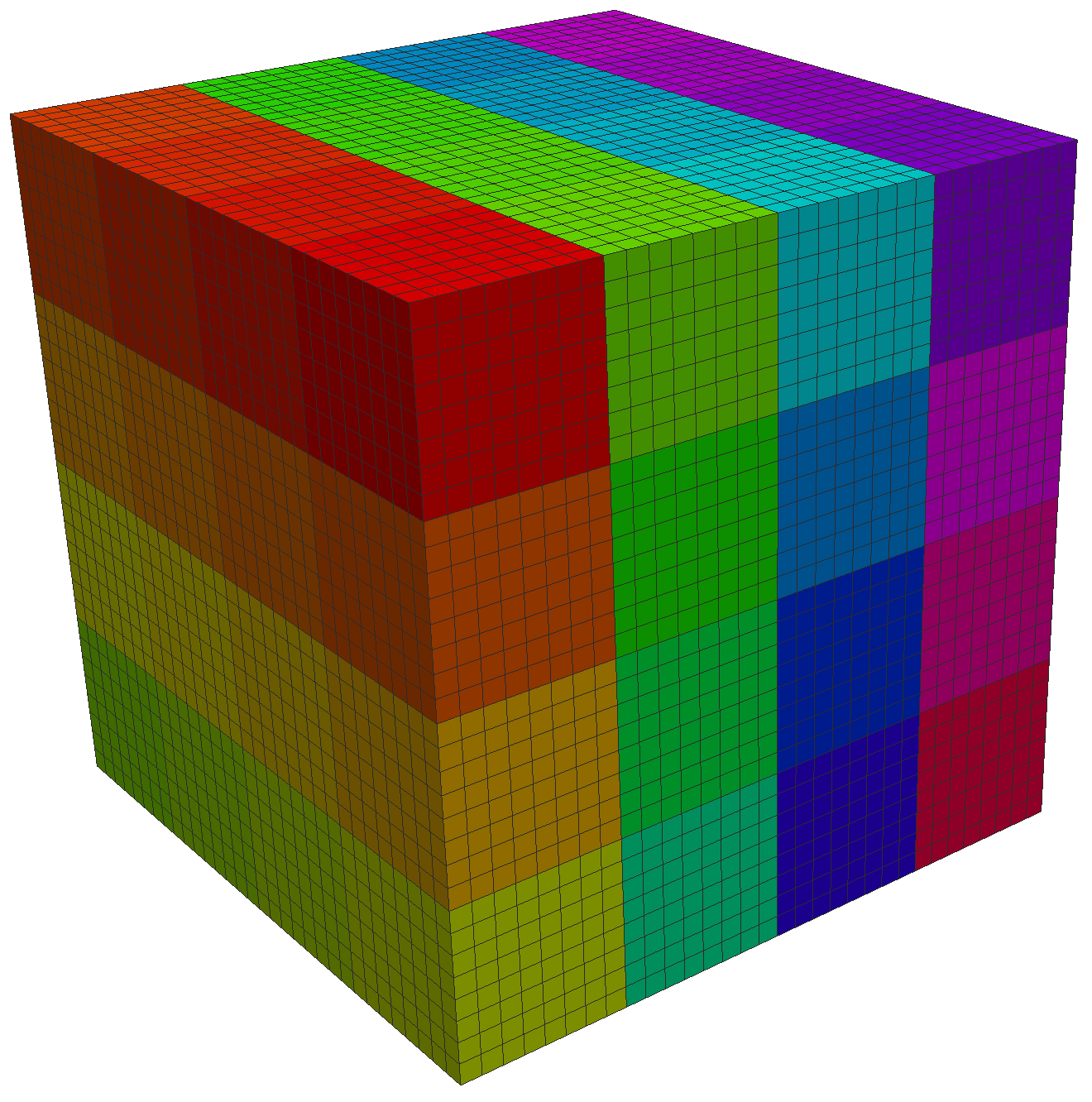}
\hspace{0.1\textwidth}
\includegraphics[height = 0.33\textwidth]{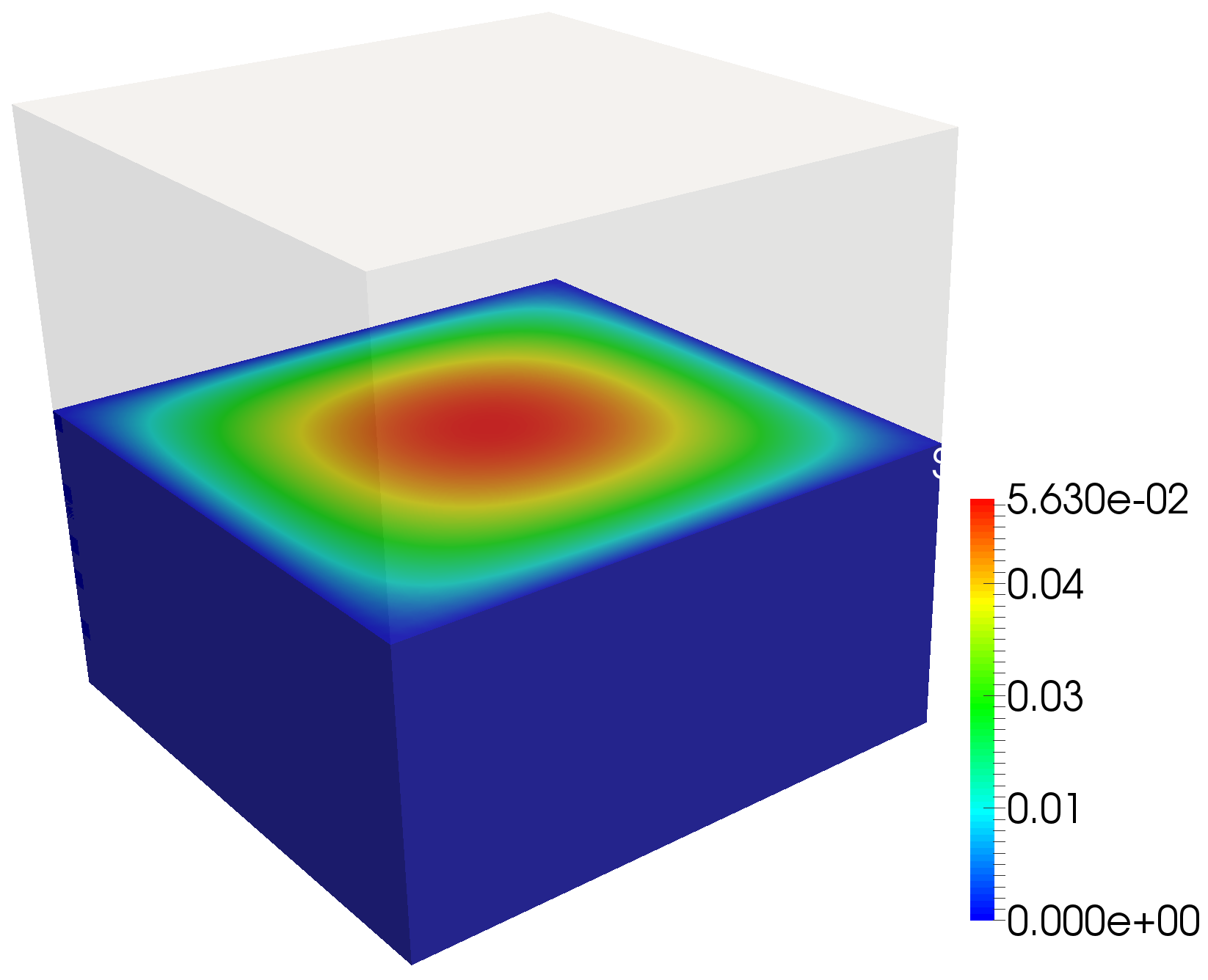}
\end{center}
\caption{\label{fig:poisson_on_cube} An example of the Poisson problem solved with regular subdomains, division into 64 subdomains with $H/h$ = 8 (left) and plot of the solution (right).}
\end{figure}

In order to have a bottom line case for further comparisons, we start by solving the 3D Poisson problem with regular cubic subdomains, see example in Fig.~\ref{fig:poisson_on_cube}. 
This is a widely used benchmark for domain decomposition solvers.
In our case, regular mesh of trilinear finite elements is used. 
First, we perform weak scaling tests for this problem, in which the subdomain size is fixed, 
and the number of subdomains, and correspondingly computing cores, is increased. 
Subdomain size is characterized by the number of elements along an edge of the cubic subdomain, known as the $H/h$ ratio, 
with $H$ being the characteristic subdomain size and $h$ the characteristic element size.
We perform the weak scaling test for $H/h$ = 16, 32, and 64, in each case comparing 2-level and 3-level BDDC method, and using
4$^3$, 5$^3$, 8$^3$, 10$^3$, and 16$^3$ subdomains.
Arithmetic averages on edges and faces of subdomains are used as the coarse degrees of freedom.

Results for the two-level BDDC method are summarized in Table~\ref{table_weak_regular_global_level2},
while Table~\ref{table_weak_regular_global_level3} presents results obtained using three levels of BDDC\@.
As usual, we are interested in the number of PCG iterations for reaching the tolerance. 
Our results confirm the mathematical theory for the BDDC method, with the number of required iterations being independent 
of the number of subdomains and processors, and increasing slightly with the $H/h$ ratio.

When using the three-level BDDC method, the number of iterations slightly increases, while it seems also independent of the number of subdomains.
The slight oscillations in the number of iterations for the 3-level method might be attributed to the fact that the mesh is partitioned into irregular subdomains on the second level.
This is again in accordance to the previously published results. 

\begin{figure}
\begin{center}
\includegraphics[width = 0.49\textwidth]{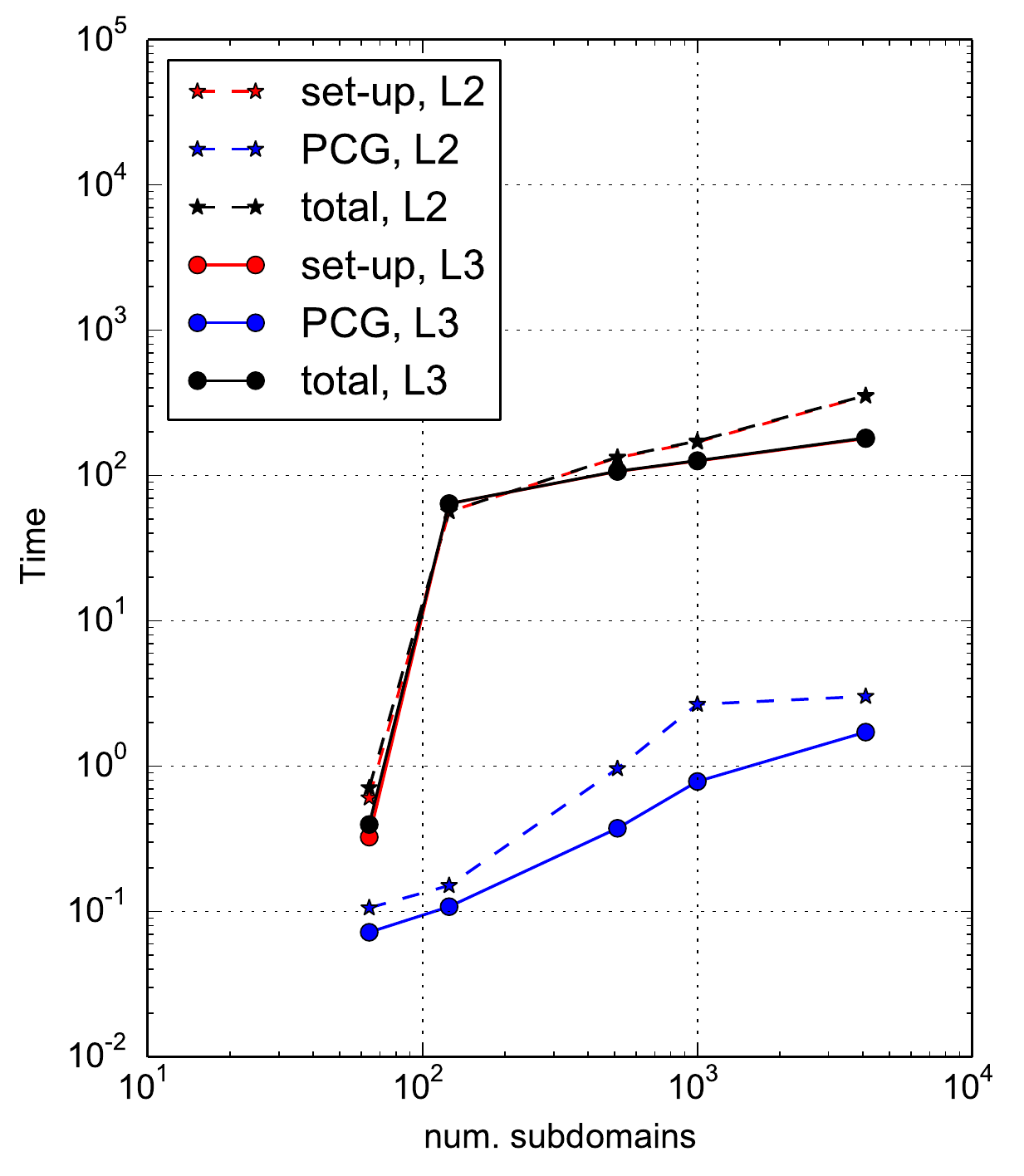}
\includegraphics[width = 0.49\textwidth]{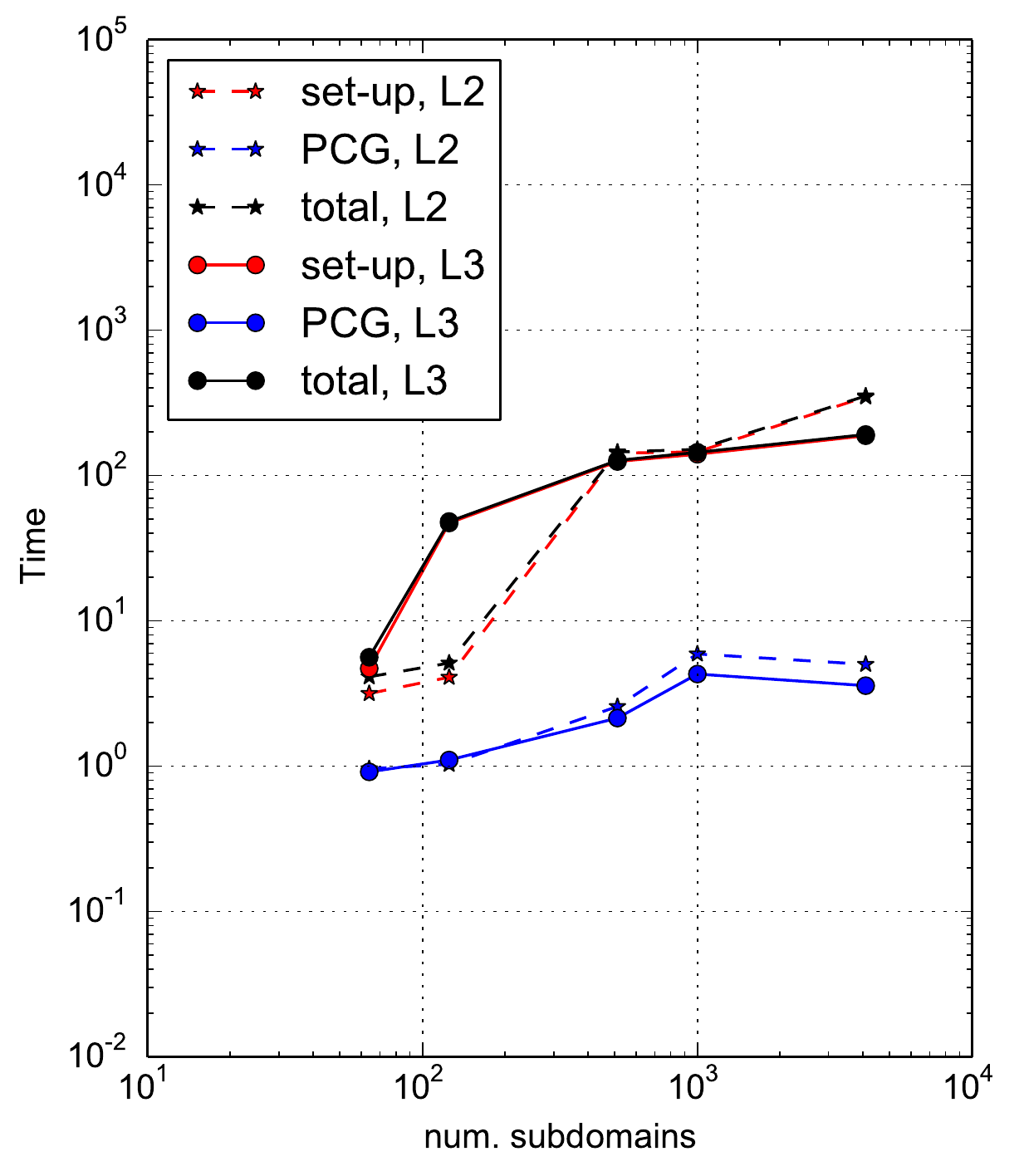} \\
\includegraphics[width = 0.49\textwidth]{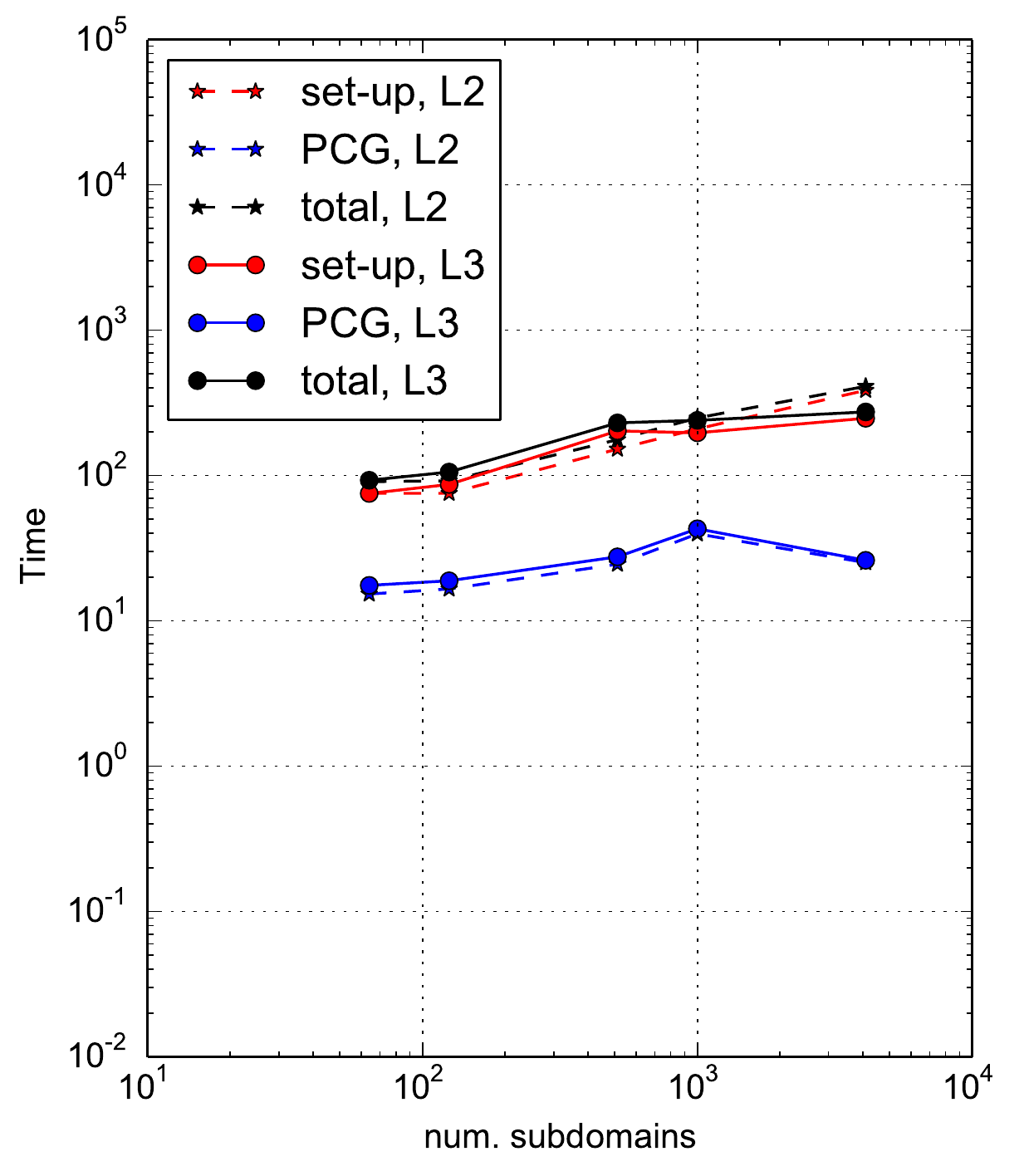}
\includegraphics[width = 0.49\textwidth]{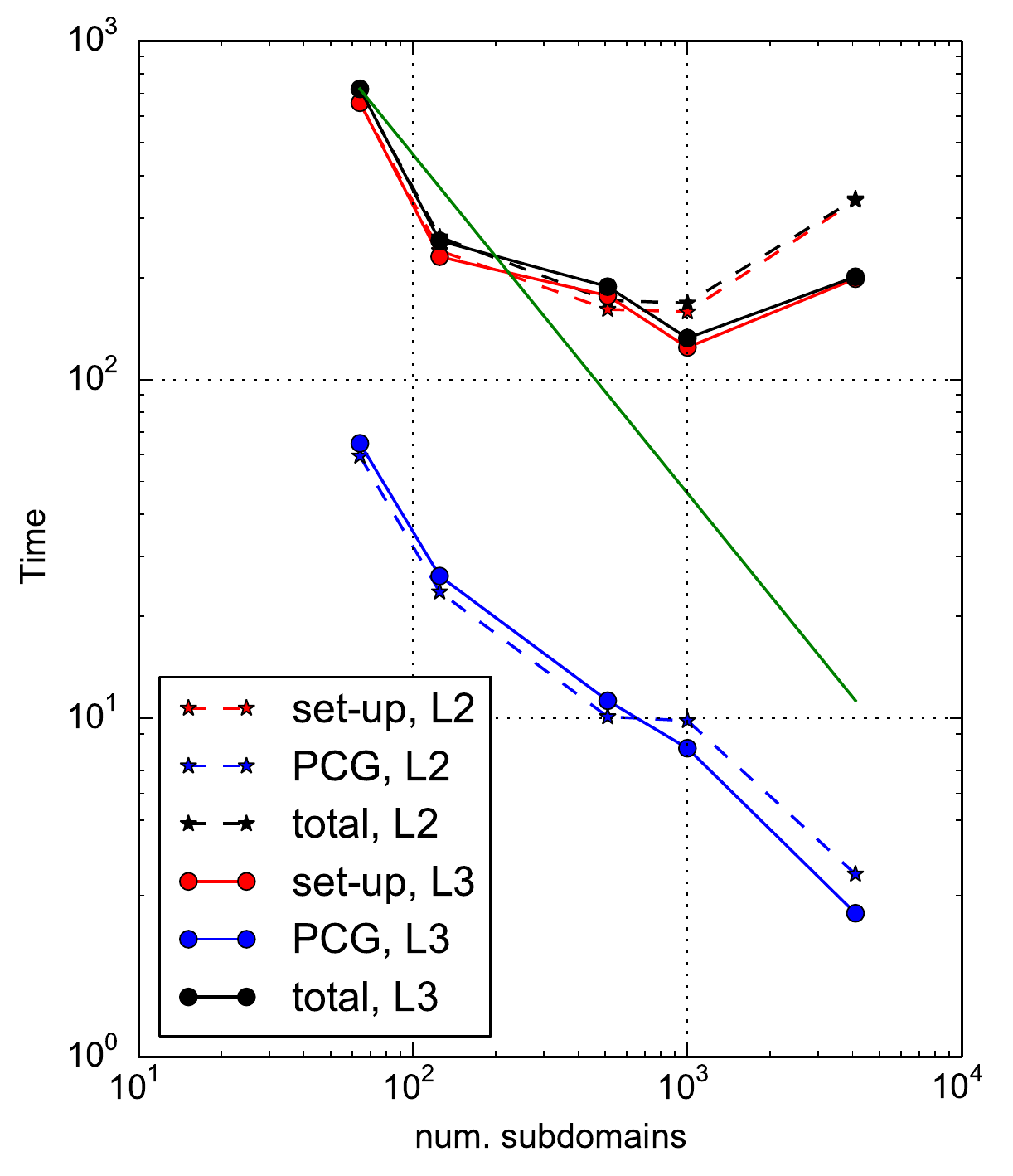} \\
\caption{\label{fig:poisson_regular_scaling1} Poisson equation in 3D, regular cubic subdomains, weak scaling for 2-level and 3-level BDDC, trilinear elements,
$H/h$ = 16 (top left), $H/h$ = 32 (top right), and $H/h$ = 64 (bottom left), and a strong scaling test (bottom right).
}
\end{center}
\end{figure}

This desirable behaviour of the multilevel BDDC method was proven theoretically (see e.g.~\cite{Mandel-2008-MMB,Tu-2007-TBT3D}), 
as well as confirmed by a number of experiments (e.g.~\cite{Badia-2016-MBD,Sousedik-2013-AMB}).
We repeat these results here to evaluate the \codename{BDDCML} solver on the \emph{Salomon} supercomputer, 
and we are interested in the computing times. 
These are also plotted in Fig.~\ref{fig:poisson_regular_scaling1}.
In the plots, we show separately time for preconditioner set-up, time for the PCG iterations, and their sum (total).
As expected, the weak scalability is slightly better for the 3-level method, 
which becomes faster than the 2-level method especially on 4096 cores.


\begin{table}
\begin{center}
\begin{tabular}{c|ccc|cc|c|cc}
$H/h$ & $N_S$ & $n$ & $n/N_S$ & $n^{\Gamma}$ & $n_{C}$ & its. & $t_{set-up}$ & $t_{PCG}$ \\
\hline
\multirow{5}{*}{16} & 
  64 & $2.7\!\cdot\! 10^{5}$ & 4291 & $3.6\!\cdot\! 10^{4}$ & 279 & 9 & 0.6 & 0.11 \\
& 125 & $5.3\!\cdot\! 10^{5}$ & 4251 & $7.5\!\cdot\! 10^{4}$ & 604 & 9 & 56 & 0.15 \\
& 512 & $2.1\!\cdot\! 10^{6}$ & 4192 & $3.3\!\cdot\! 10^{5}$ & 2863 & 9 & 132 & 0.96 \\
& 1000 & $4.2\!\cdot\! 10^{6}$ & 4173 & $6.6\!\cdot\! 10^{5}$ & 5859 & 9 & 169 & 2.7 \\
& 4096 & $1.7\!\cdot\! 10^{7}$ & 4144 & $2.8\!\cdot\! 10^{6}$ & $2.6\!\cdot\! 10^{4}$ & 9 & 352 & 3 \\
\hline
\multirow{5}{*}{32} &  
  64 & $2.1\!\cdot\! 10^{6}$ & $3.4\!\cdot\! 10^{4}$ & $1.5\!\cdot\! 10^{5}$ & 279 & 11 & 3.2 & 0.96 \\
& 125 & $4.2\!\cdot\! 10^{6}$ & $3.3\!\cdot\! 10^{4}$ & $3.0\!\cdot\! 10^{5}$ & 604 & 11 & 4.1 & 1 \\
& 512 & $1.7\!\cdot\! 10^{7}$ & $3.3\!\cdot\! 10^{4}$ & $1.3\!\cdot\! 10^{6}$ & 2863 & 11 & 142 & 2.6 \\
& 1000 & $3.3\!\cdot\! 10^{7}$ & $3.3\!\cdot\! 10^{4}$ & $2.7\!\cdot\! 10^{6}$ & 5859 & 11 & 145 & 5.9 \\
& 4096 & $1.4\!\cdot\! 10^{8}$ & $3.3\!\cdot\! 10^{4}$ & $1.1\!\cdot\! 10^{7}$ & $2.6\!\cdot\! 10^{4}$ & 11 & 347 & 5 \\
\hline
\multirow{5}{*}{64} &
  64 & $1.7\!\cdot\! 10^{7}$ & $2.7\!\cdot\! 10^{5}$ & $5.9\!\cdot\! 10^{5}$ & 279 & 12 & 75 & 15 \\
& 125 & $3.3\!\cdot\! 10^{7}$ & $2.6\!\cdot\! 10^{5}$ & $1.2\!\cdot\! 10^{6}$ & 604 & 13 & 75 & 17 \\
& 512 & $1.4\!\cdot\! 10^{8}$ & $2.6\!\cdot\! 10^{5}$ & $5.5\!\cdot\! 10^{6}$ & 2863 & 14 & 151 & 25 \\
& 1000 & $2.6\!\cdot\! 10^{8}$ & $2.6\!\cdot\! 10^{5}$ & $1.1\!\cdot\! 10^{7}$ & 5859 & 14 & 209 & 40 \\
& 4096 & $1.1\!\cdot\! 10^{9}$ & $2.6\!\cdot\! 10^{5}$ & $4.7\!\cdot\! 10^{7}$ & $2.6\!\cdot\! 10^{4}$ & 14 & 385 & 25 \\
\end{tabular}
\caption{\label{table_weak_regular_global_level2}Weak scaling test, Poisson equation in 3D, regular subdomains, 2-level BDDC;
Columns correspond to number of subdomains $N_S$, size of the global linear system $n$, average subdomain size $n/N_S$, size of the interface problem $n^{\Gamma}$,
global number of coarse degrees of freedom $n_{C}$, number of PCG iterations, 
time of preconditioner set-up $t_{set-up}$ in seconds, and time spent in all PCG iterations $t_{PCG}$ in seconds.}
\end{center}
\end{table}

\begin{table}
\begin{center}
\begin{tabular}{c|ccc|cc|c|cc}
$H/h$ & $N_S$ & $n$ & $n/N_S$ & $n^{\Gamma}$ & $n_{C}$ & its. & $t_{set-up}$ & $t_{PCG}$ \\
\hline
\multirow{5}{*}{16} & 
  64/8 & $2.7\!\cdot\! 10^{5}$ & 4291 & $3.6\!\cdot\! 10^{4}$/118 & 279/25 & 9 & 0.33 & 0.072 \\
& 125/12 & $5.3\!\cdot\! 10^{5}$ & 4251 & $7.5\!\cdot\! 10^{4}$/336 & 604/106 & 11 & 64 & 0.11 \\
& 512/22 & $2.1\!\cdot\! 10^{6}$ & 4192 & $3.3\!\cdot\! 10^{5}$/1544 & 2863/266 & 14 & 106 & 0.37 \\
& 1000/32 & $4.2\!\cdot\! 10^{6}$ & 4173 & $6.6\!\cdot\! 10^{5}$/2069 & 5859/308 & 12 & 125 & 0.79 \\
& 4096/64 & $1.7\!\cdot\! 10^{7}$ & 4144 & $2.8\!\cdot\! 10^{6}$/7452 & $2.6\!\cdot\! 10^{4}$/809 & 13 & 179 & 1.7 \\
\hline
\multirow{5}{*}{32} & 
  64/8 & $2.1\!\cdot\! 10^{6}$ & $3.4\!\cdot\! 10^{4}$ & $1.5\!\cdot\! 10^{5}$/118 & 279/25 & 11 & 4.7 & 0.91 \\
& 125/12 & $4.2\!\cdot\! 10^{6}$ & $3.3\!\cdot\! 10^{4}$ & $3.0\!\cdot\! 10^{5}$/336 & 604/106 & 13 & 47 & 1.1 \\
& 512/22 & $1.7\!\cdot\! 10^{7}$ & $3.3\!\cdot\! 10^{4}$ & $1.3\!\cdot\! 10^{6}$/1544 & 2863/266 & 15 & 124 & 2.1 \\
& 1000/32 & $3.3\!\cdot\! 10^{7}$ & $3.3\!\cdot\! 10^{4}$ & $2.7\!\cdot\! 10^{6}$/2069 & 5859/308 & 14 & 140 & 4.3 \\
& 4096/64 & $1.4\!\cdot\! 10^{8}$ & $3.3\!\cdot\! 10^{4}$ & $1.1\!\cdot\! 10^{7}$/7452 & $2.6\!\cdot\! 10^{4}$/809 & 15 & 187 & 3.6 \\
\hline
\multirow{5}{*}{64} & 
  64/8 & $1.7\!\cdot\! 10^{7}$ & $2.7\!\cdot\! 10^{5}$ & $5.9\!\cdot\! 10^{5}$/118 & 279/25 & 14 & 75 & 18 \\
& 125/12 & $3.3\!\cdot\! 10^{7}$ & $2.6\!\cdot\! 10^{5}$ & $1.2\!\cdot\! 10^{6}$/336 & 604/106 & 15 & 87 & 19 \\
& 512/22 & $1.4\!\cdot\! 10^{8}$ & $2.6\!\cdot\! 10^{5}$ & $5.5\!\cdot\! 10^{6}$/1544 & 2863/266 & 17 & 202 & 28 \\
& 1000/32 & $2.6\!\cdot\! 10^{8}$ & $2.6\!\cdot\! 10^{5}$ & $1.1\!\cdot\! 10^{7}$/2069 & 5859/308 & 16 & 196 & 43 \\
& 4096/64 & $1.1\!\cdot\! 10^{9}$ & $2.6\!\cdot\! 10^{5}$ & $4.7\!\cdot\! 10^{7}$/7452 & $2.6\!\cdot\! 10^{4}$/809 & 17 & 247 & 26 \\
\end{tabular}
\caption{\label{table_weak_regular_global_level3} Weak scaling test, Poisson equation in 3D, regular subdomains, 3-level BDDC;
number of subdomains on the first/second level $N_S$, other columns have the same meaning as in Table~\ref{table_weak_regular_global_level2}.}
\end{center}
\end{table}

As described in Section~\ref{sec:handling_disconnected},
disconnected subdomains by \codename{p4est} are handled by generating constraints 
independently for each component of the substructure.
It is therefore important to evaluate the effect of the varying number of constraints on individual subdomains on the load balance of the solver.
In particular, this mostly affects the factorization of the matrix of the Neumann problem (\ref{eq:Neumann-matrix}) and the 
time of back-substitution in the problem with multiple right-hand sides to find the local coarse basis functions (\ref{eq:Neumann-setup}). 

We are interested in the number of constraints on each subdomain and the time required by the \codename{MUMPS} solver for factorization and back-substitution. 
These \emph{local subdomain properties} are summarized in Table~\ref{table_weak_regular_local}. 
For each case, we present minimum, maximum, and average number of constraints on subdomains and corresponding times for the two operations.
It is interesting to realize that even for this regular case, subdomains inside the cube have 3$\times$ more local constraints than those at the corners
of the domain.

Based on the results, we can conclude that this difference in number of local constraints does not cause significant imbalance of time for factorization of the matrix.
However, the time for the back-substitution part scales with the number of constraints, 
as \codename{MUMPS} does not seem to make a good data reuse of these multiple vectors in version 4.10.0.
We believe that this behaviour might be improved with a newer versions of the \codename{MUMPS} solver.
Fortunately, time for this part of the algorithm is significantly lower than that for the factorization.
We would like to remind that both parts are performed once for each subdomain during the set-up of the BDDC preconditioner. 

\hl{Based on the times for the different $H/h$ ratios}, 
we can also estimate the cost of the sparse direct $LDL^T$ factorization performed 
in \codename{MUMPS} in this case as proportional to $(1/h^{3})^{\alpha}$, with $\alpha \approx 1.7$.

In fact, the local properties could be shown for both 2- and 3-level BDDC method, and for both levels for the latter. 
However, the results shown that properties of the first-level subdomains are critical with respect to time, 
and these are the same for 2- and 3-level method.
Therefore, other tables are omitted in all presented computations.

\begin{table}
\begin{center}
\begin{tabular}{c|c|ccc|ccc|ccc}
\multirow{2}{*}{$H/h$} & \multirow{2}{*}{$N_S$} & \multicolumn{3}{c|}{num.\ coarse dofs} & \multicolumn{3}{c|}{time fact.\ loc.\ (s)} & \multicolumn{3}{c}{time sol.\ loc.\ (s)}\\ 
  & & min & max & avg & min & max & avg & min & max & avg\\
\hline
\multirow{5}{*}{16} & 
  64 & 6 & 18 & 11 & 0.05 & 0.06 & 0.06 & 0.01 & 0.02 & 0.01 \\
& 125 & 6 & 18 & 12 & 0.05 & 0.06 & 0.05 & 0.01 & 0.02 & 0.01 \\
& 512 & 6 & 18 & 14 & 0.05 & 0.07 & 0.05 & 0.01 & 0.02 & 0.01 \\
& 1000 & 6 & 18 & 15 & 0.05 & 0.11 & 0.05 & 0.01 & 0.03 & 0.01 \\
& 4096 & 6 & 18 & 16 & 0.05 & 0.07 & 0.05 & 0.01 & 0.02 & 0.02 \\
\hline
\multirow{5}{*}{32} & 
  64 & 6 & 18 & 11 & 1.25 & 1.34 & 1.28 & 0.07 & 0.19 & 0.12 \\
& 125 & 6 & 18 & 12 & 1.25 & 1.34 & 1.28 & 0.06 & 0.19 & 0.13 \\
& 512 & 6 & 18 & 14 & 1.28 & 1.41 & 1.32 & 0.07 & 0.20 & 0.16 \\
& 1000 & 6 & 18 & 15 & 1.28 & 1.60 & 1.32 & 0.07 & 0.21 & 0.16 \\
& 4096 & 6 & 18 & 16 & 1.28 & 1.43 & 1.33 & 0.07 & 0.21 & 0.18 \\

\hline
\multirow{5}{*}{64} & 
  64 & 6 & 18 & 11 & 34.89 & 38.33 & 36.79 & 0.97 & 2.62 & 1.75 \\
& 125 & 6 & 18 & 12 & 34.55 & 38.19 & 36.38 & 0.93 & 2.63 & 1.88 \\
& 512 & 6 & 18 & 14 & 36.64 & 41.52 & 39.12 & 0.98 & 2.82 & 2.25 \\
& 1000 & 6 & 18 & 15 & 36.89 & 41.66 & 39.80 & 1.03 & 2.88 & 2.38 \\
& 4096 & 6 & 18 & 16 & 36.41 & 42.18 & 39.90 & 1.02 & 3.03 & 2.57 \\
\end{tabular}
\caption{\label{table_weak_regular_local} Weak scaling test, Poisson equation in 3D, regular subdomains, local subdomain properties; 
number of coarse degrees of freedom on subdomain (minimum, maximum and average over all subdomains), 
time for factorization of the local matrix from equation (\ref{eq:Neumann-matrix}), 
time for back-substitution and solution of the local problem in equation (\ref{eq:Neumann-setup}).}
\end{center}
\end{table}

We conclude the study with regular subdomains by performing a strong scaling test. 
In this test, the number of elements is fixed to 400$\times$400$\times$400 = 64$\cdot$10$^{6}$ elements with 64.5$\cdot$10$^{6}$ degrees of freedom.
This problem is solved on 64, 125, 512, 1000, and 4096 regular subdomains, and the results are summarized in Table~\ref{table_regular_global_level2} for 2-level,
and Table~\ref{table_regular_global_level3} for 3-level BDDC\@.
The number of iterations does not grow with number of subdomains, it even slightly decreases, again in an encouraging agreement with the theory for decreasing $H/h$ ratio.
Resulting times are also plotted in Fig.~\ref{fig:poisson_regular_scaling1}.
They show again the superiority of the 3-level method for 1000 and 4096 subdomains. 
However, these scaling results seem to be affected by the observed suboptimal behaviour of the supercomputer 
for larger numbers of cores.

\begin{table}
\begin{center}
\begin{tabular}{ccc|cc|c|cc}
$N_S$ & $n$ & $n/N_S$ & $n^{\Gamma}$ & $n_{C}$ & its. & $t_{set-up}$ & $t_{PCG}$ \\
\hline
64 & $6.4\!\cdot\! 10^{7}$ & $1.0\!\cdot\! 10^{6}$ & $1.4\!\cdot\! 10^{6}$ & 279 & 14 & 656 & 59 \\
125 & $6.4\!\cdot\! 10^{7}$ & $5.2\!\cdot\! 10^{5}$ & $1.9\!\cdot\! 10^{6}$ & 604 & 14 & 240 & 24 \\
512 & $6.4\!\cdot\! 10^{7}$ & $1.3\!\cdot\! 10^{5}$ & $3.3\!\cdot\! 10^{6}$ & 2863 & 13 & 161 & 10 \\
1000 & $6.4\!\cdot\! 10^{7}$ & $6.4\!\cdot\! 10^{4}$ & $4.2\!\cdot\! 10^{6}$ & 5859 & 12 & 158 & 9.8 \\
4096 & $6.4\!\cdot\! 10^{7}$ & $1.6\!\cdot\! 10^{4}$ & $7.0\!\cdot\! 10^{6}$ & $2.6\!\cdot\! 10^{4}$ & 10 & 338 & 3.5 \\
\end{tabular}
\caption{\label{table_regular_global_level2} Strong scaling test, Poisson equation in 3D, regular subdomains, 2-level BDDC.}
\end{center}
\end{table}

\begin{table}
\begin{center}
\begin{tabular}{ccc|cc|c|cc}
$N_S$ & $n$ & $n/N_S$ & $n^{\Gamma}$ & $n_{C}$ & its. & $t_{set-up}$ & $t_{PCG}$ \\
\hline
64/8 & $6.4\!\cdot\! 10^{7}$    & $1.0\!\cdot\! 10^{6}$ & $1.4\!\cdot\! 10^{6}$/118 & 279/25 & 15 & 656 & 65 \\
125/12 & $6.4\!\cdot\! 10^{7}$  & $5.2\!\cdot\! 10^{5}$ & $1.9\!\cdot\! 10^{6}$/336 & 604/106 & 16 & 230 & 26 \\
512/22 & $6.4\!\cdot\! 10^{7}$  & $1.3\!\cdot\! 10^{5}$ & $3.3\!\cdot\! 10^{6}$/1544 & 2863/266 & 17 & 176 & 11 \\
1000/32 & $6.4\!\cdot\! 10^{7}$ & $6.4\!\cdot\! 10^{4}$ & $4.2\!\cdot\! 10^{6}$/2069 & 5859/308 & 14 & 124 & 8.2 \\
4096/64 & $6.4\!\cdot\! 10^{7}$ & $1.6\!\cdot\! 10^{4}$ & $7.0\!\cdot\! 10^{6}$/7452 & $2.6\!\cdot\! 10^{4}$/809 & 14 & 198 & 2.7 \\
\end{tabular}
\caption{\label{table_regular_global_level3} Strong scaling test, Poisson equation in 3D, regular subdomains, 3-level BDDC.}
\end{center}
\end{table}

\begin{table}
\begin{center}
\begin{tabular}{c|c|ccc|ccc|ccc}
\multirow{2}{*}{$N_S$} & \multirow{2}{*}{$n/N_S$} & \multicolumn{3}{c|}{num.\ coarse dofs} & \multicolumn{3}{c|}{time fact.\ loc.\ (s)} & \multicolumn{3}{c}{time sol.\ loc.\ (s)}\\ 
                       &                          & min & max & avg & min & max & avg & min & max & avg \\
\hline
64   & $1.0\!\cdot\! 10^{6}$ & 6 & 18 & 11 & 301.78 & 344.43 & 320.09 & 4.22 & 12.03 & 7.71 \\
125  & $5.2\!\cdot\! 10^{5}$ & 6 & 18 & 12 & 93.03 & 100.83 & 96.19 & 1.81 & 5.16 & 3.68 \\
512  & $1.3\!\cdot\! 10^{5}$ & 6 & 18 & 14 & 10.21 & 11.31 & 10.77 & 0.37 & 1.06 & 0.85 \\
1000 & $6.4\!\cdot\! 10^{4}$ & 6 & 18 & 15 & 3.53 & 4.18 & 3.73 & 0.16 & 0.47 & 0.38 \\
4096 & $1.6\!\cdot\! 10^{4}$ & 6 & 18 & 16 & 0.43 & 0.48 & 0.45 & 0.03 & 0.08 & 0.07 \\
\end{tabular}
\caption{\label{table_regular_local} Strong scaling test, Poisson equation in 3D, regular subdomains, local subdomain properties.}
\end{center}
\end{table}

\begin{figure}
\begin{center}
\includegraphics[width = 0.35\textwidth]{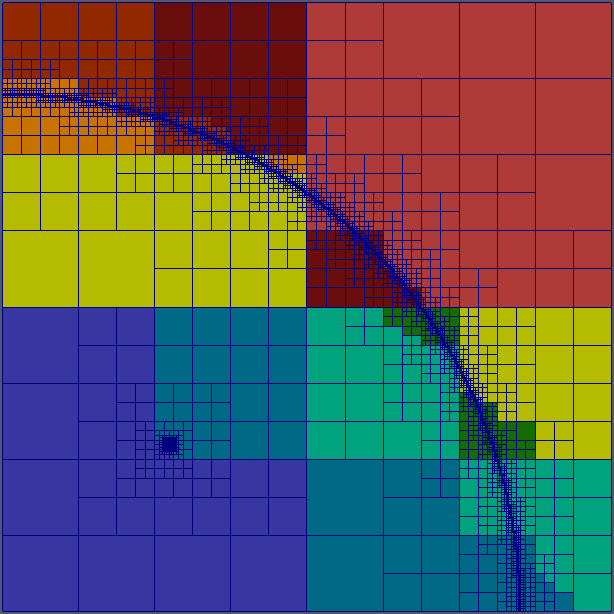}
\hspace{1cm}
\includegraphics[width = 0.4\textwidth]{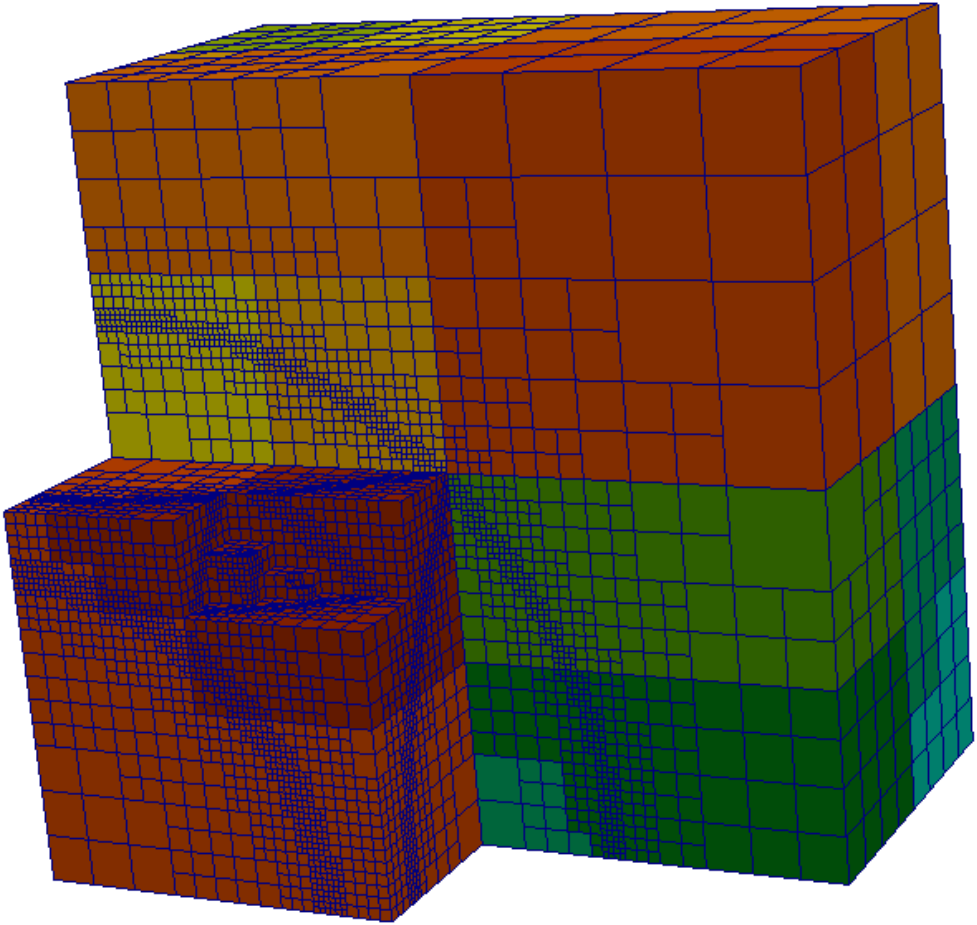}
\end{center}
\caption{\label{fig:scaling_meshes}Example meshes obtained by combination of uniform, square and circle refinements in 2D and 3D. 
Partitioning and preserving 1:2 balance done by the \codename{p4est} library. 
In the 3D case, only subdomains 1-16 of 20 are shown.}
\end{figure}

\subsubsection{Refined meshes in 2D and 3D}
Let us now investigate the behaviour of the parallel BDDC solver for adaptively refined meshes with hanging nodes.
It is not obvious to us how a scaling test could be designed for a parallel adaptive finite element algorithm. 
For this reason, 
we perform strong scaling tests on fixed refined meshes. 
These are created by prescribed refinements of an original one element mesh in such a way that all specific 
phenomena of adapted meshes, such as hanging nodes, different sizes of elements in different parts of the 
mesh, and disconnected subdomains are present. 
Let $R_U$ stand for a uniform refinement, for which all elements in the mesh are subdivided into 4 (in 2D) or 8 (in 3D) new elements, 
let $R_C$ stand for `circular' refinement, for which all elements intersected by a circle or a sphere with centre at the origin and radius 0.85 are refined. 
The last considered refinement $R_S$ (`square') refines all elements which have nonempty intersection with a square or cube $[0.26, 0.28]^d$, $d = 2,3$.
Although these refinements might look somewhat artificial, they produce reasonable testing meshes such as those shown in Fig.~\ref{fig:scaling_meshes}. 
Note that refinements concentrate around the circle and the square, 
which might predict the behaviour of an adaptivity algorithm with a presence of an internal layer or a singularity. 
Arithmetic averages over the edges and faces of subdomains are again used as the coarse degrees of freedom.

Let us start with comparing the solver behaviour for linear finite elements in two and three spatial dimensions.
The 2-D problem is created by 10$\times$ $R_U$, 14$\times$ $R_C$, and 7$\times$ $R_S$, 
while the 3-D problem is obtained by 7$\times$ $R_U$, 5$\times$ $R_C$, and 5$\times$ $R_S$. 
An illustrative example of meshes created by much less refinements of the same type are in Fig.~\ref{fig:scaling_meshes}.

In 2D, the resulting mesh has 150 million elements and 120 million unknowns, 
and in 3D, the considered problem has 
approx.\ 87 million elements translating to about 59 million unknowns. 
Bilinear and trilinear elements are used, and the global number of unknowns is reduced by the fact that degrees of freedom are not assigned to hanging nodes. 
Note that since $R_C$ and $R_S$ are focused in different parts of the domain, the maximal level of refinements can be obtained
as $R_U + \max(R_C, R_S)$, which is 24 for 2D and 12 for 3D. This is still within \codename{p4est} restrictions of maximal 
level of refinements 29 for 2D and 19 for 3D. A relatively high number of refinement levels in our test cases is caused by 
using a single element as an initial mesh and doing a lot of refinements in a relatively small area.

\begin{figure}
\begin{center}
\includegraphics[width = 0.49\textwidth]{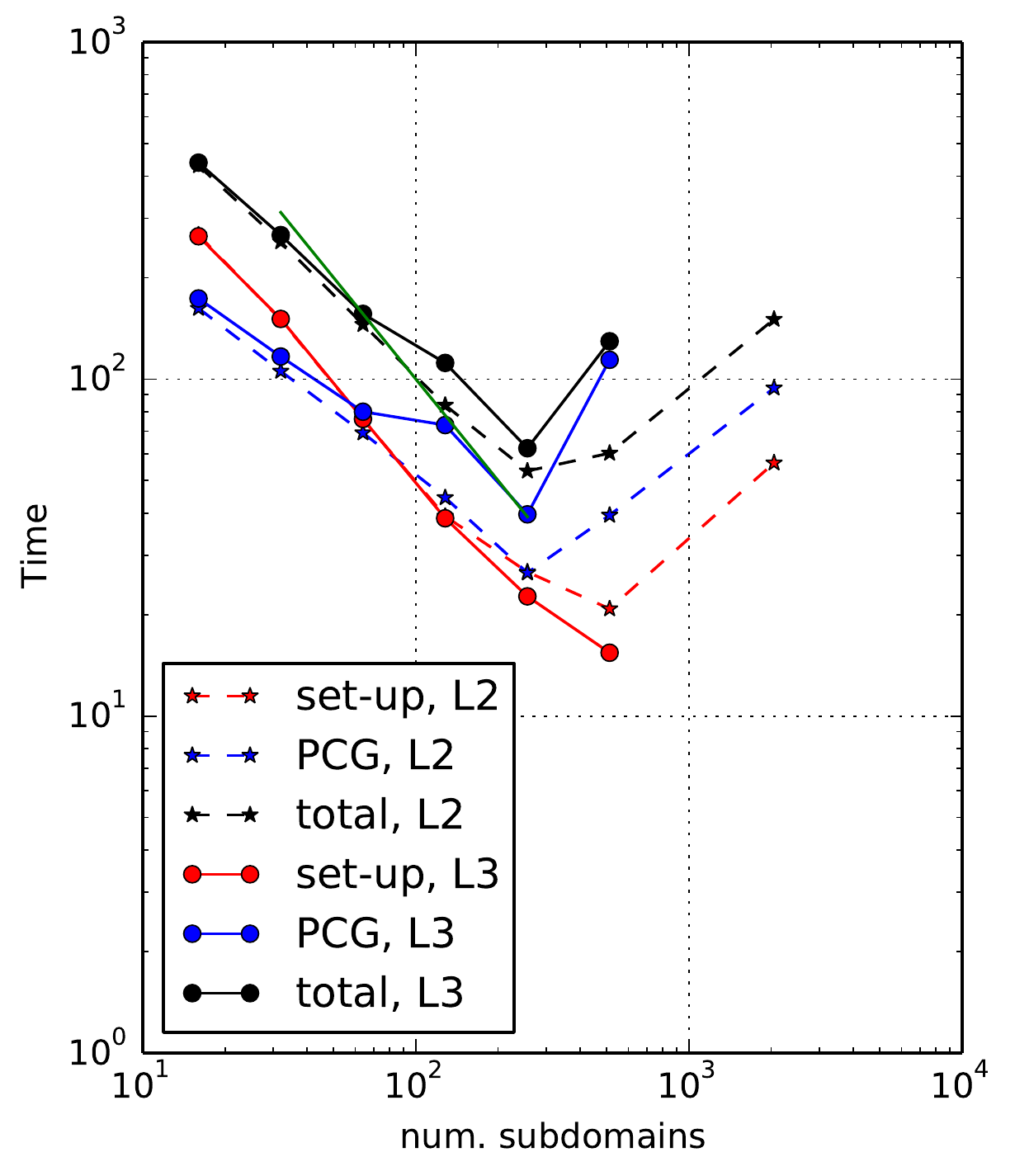}
\includegraphics[width = 0.49\textwidth]{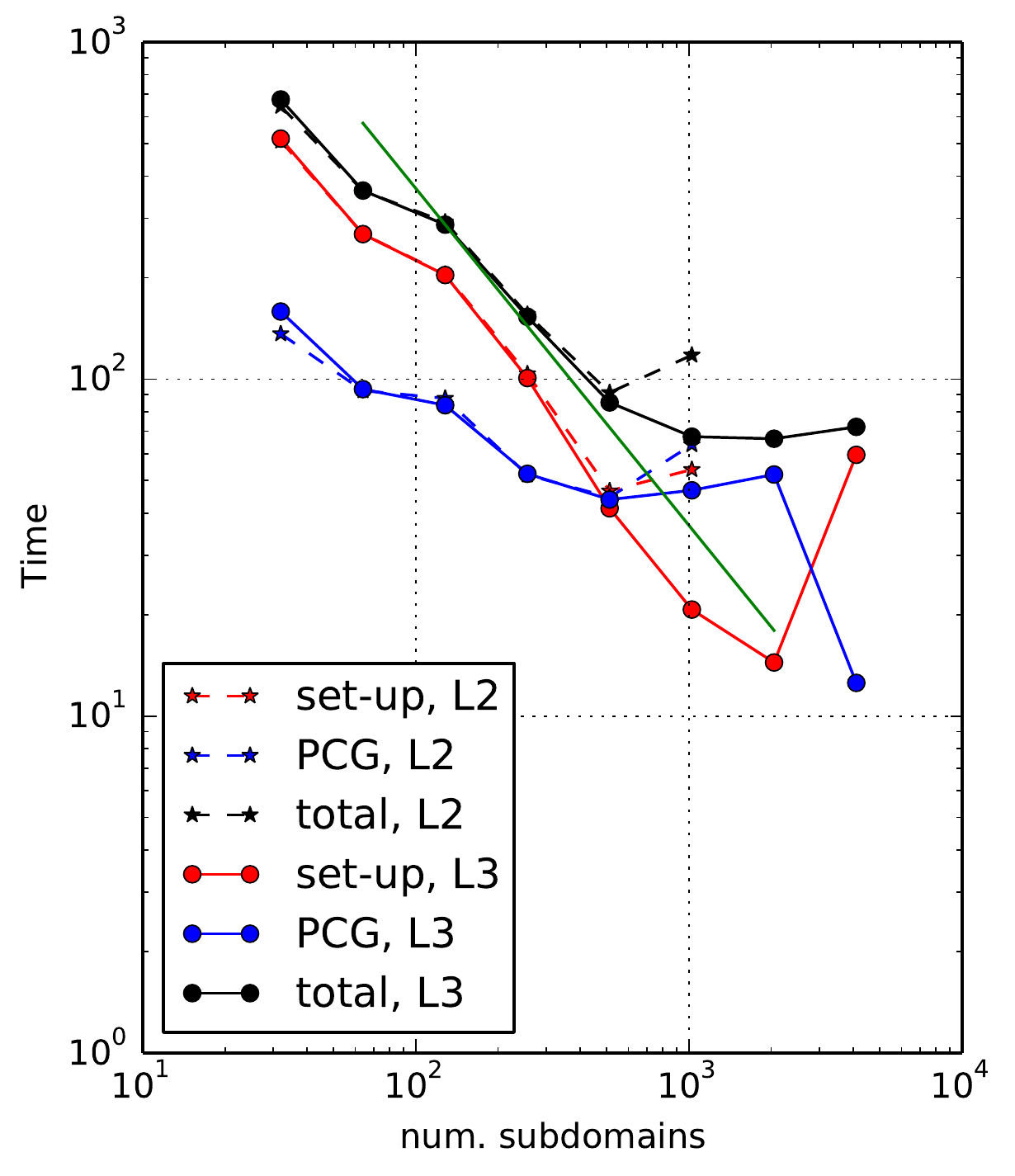} 
\caption{\label{fig:poisson_first_order_scaling}Poisson equation for meshes with a prescribed combination of refinements in 2D (left) and 3D (right). First order (bilinear and trilinear) elements were used. In both cases performance of 2 and 3-level BDDC is compared. 
Data can be found in Tables~\ref{table_glob_2D_level2} and~\ref{table_glob_2D_level3} for 2D, and Tables~\ref{table_glob_3D_level2} and~\ref{table_glob_3D_level3} for 3D.}
\end{center}
\end{figure}

Results of the scaling test with the 2-D problem are presented in Table~\ref{table_glob_2D_level2} for the 2-level method, 
and in Table~\ref{table_glob_2D_level3} for the 3-level method.
Similarly, 3-D results are presented in Tables~\ref{table_glob_3D_level2} and~\ref{table_glob_3D_level3}.
The computational times of the BDDC solver are also shown in Fig.~\ref{fig:poisson_first_order_scaling}. 

For the 2-D problem, 2-level BDDC outperforms the 3-level method not only in the number of iterations (as expected), but also in computational time. 
Neither of these methods has converged for 1024 subdomains, but the 2-level method also provided results for 2048 subdomains.

The situation has changed for the 3-D case, for which the 2-level method did not fit into memory for more than 1024 cores, 
while the 3-level method allowed us to continue the test up to 4096 subdomains.
Even on 1024 cores, the 3-level method has provided a significantly faster solution.
Both 2- and 3-level methods have shown a very acceptable growth in the number of required PCG iterations.
This is seen in the plots by the fact, that the curve for time of iterations has a worse slope than the one for preconditioner set-up.

Let us now have a look at the local subdomain properties in this experiment,
presented in Tables~\ref{table_loc_2D} and~\ref{table_loc_3D} for 2D and 3D, respectively. 
One can see that for these refined meshes, the ratio of the minimum and maximum number of constrains grows to approx.\ 15 in 2D and approx.\ 20 in 3D.
While the ratio of the average value to minimum is still only about 4, there are subdomains with high number of local coarse degrees of freedom. 
Unfortunately, these subdomains dictate the time of the parallel execution of this stage, and we can see some imbalances even in the factorization phase now.
The ratio of maximum and minimum time is roughly 2 in 2D and about 3 in 3D.
These ratios become much higher for the back-substitution in problem (\ref{eq:Neumann-setup}), 
and in the worst cases might take even longer than the factorization.

However, an interesting insight can be obtained by comparing the 3-D results of Table~\ref{table_loc_3D} with Table~\ref{table_regular_local}
obtained for regular subdomains, which are about the same size. 
One can see, that the maximum time for factorization is not worsening a lot, while the minimum is much smaller compared to the regular case. 
This means, that the imbalance is actually caused by the fact that some subdomain problems are getting simpler to be solved for the sparse direct solver,
rather than more difficult by additional constraints. 
It is also worth noting that most of the subdomains have two components (non-contiguous parts of the subdomain). 
The exact percentage varies, but it is between $60$ and $80$ percent of subdomains for meshes presented in this section. 
Interestingly, we have not detected any subdomain with three or more components.


\begin{table}
\begin{center}
\begin{tabular}{ccc|cc|c|cc}
$N_S$ & $n$ & $n/N_S$ & $n^{\Gamma}$ & $n_{C}$ & its. & $t_{set-up}$ & $t_{PCG}$ \\
\hline
16 & $1.2\!\cdot\! 10^{8}$ & $7.5\!\cdot\! 10^{6}$ & $1.6\!\cdot\! 10^{4}$ & 72 & 23 & 267 & 162 \\
32 & $1.2\!\cdot\! 10^{8}$ & $3.7\!\cdot\! 10^{6}$ & $2.5\!\cdot\! 10^{4}$ & 151 & 25 & 150 & 105 \\
64 & $1.2\!\cdot\! 10^{8}$ & $1.9\!\cdot\! 10^{6}$ & $4.4\!\cdot\! 10^{4}$ & 341 & 28 & 76 & 69 \\
128 & $1.2\!\cdot\! 10^{8}$ & $9.3\!\cdot\! 10^{5}$ & $7.0\!\cdot\! 10^{4}$ & 707 & 29 & 39 & 44 \\
256 & $1.2\!\cdot\! 10^{8}$ & $4.7\!\cdot\! 10^{5}$ & $1.2\!\cdot\! 10^{5}$ & 1862 & 21 & 27 & 27 \\
512 & $1.2\!\cdot\! 10^{8}$ & $2.3\!\cdot\! 10^{5}$ & $1.9\!\cdot\! 10^{5}$ & 3713 & 22 & 21 & 39 \\
2048 & $1.2\!\cdot\! 10^{8}$ & $5.8\!\cdot\! 10^{4}$ & $5.3\!\cdot\! 10^{5}$ & $1.5\!\cdot\! 10^{4}$ & 24 & 56 & 94 \\
\end{tabular}
\caption{\label{table_glob_2D_level2} Strong scaling test, Poisson equation in 2D, mesh refined by prescribed refinements, 2-level BDDC.}
\end{center}
\end{table}

\begin{table}
\begin{center}
\begin{tabular}{ccc|cc|c|cc}
$N_S$ & $n$ & $n/N_S$ & $n^{\Gamma}$ & $n_{C}$ & its. & $t_{set-up}$ & $t_{PCG}$ \\
\hline
16/4 & $1.2\!\cdot\! 10^{8}$   & $7.5\!\cdot\! 10^{6}$ & $1.6\!\cdot\! 10^{4}$/25 & 72/10 & 25 & 265 & 173 \\
32/6 & $1.2\!\cdot\! 10^{8}$   & $3.7\!\cdot\! 10^{6}$ & $2.5\!\cdot\! 10^{4}$/46 & 151/20 & 28 & 150 & 116 \\
64/8 & $1.2\!\cdot\! 10^{8}$   & $1.9\!\cdot\! 10^{6}$ & $4.4\!\cdot\! 10^{4}$/81 & 341/26 & 34 & 76 & 80 \\
128/12 & $1.2\!\cdot\! 10^{8}$ & $9.3\!\cdot\! 10^{5}$ & $7.0\!\cdot\! 10^{4}$/202 & 707/40 & 49 & 39 & 73 \\
256/16 & $1.2\!\cdot\! 10^{8}$ & $4.7\!\cdot\! 10^{5}$ & $1.2\!\cdot\! 10^{5}$/368 & 1862/66 & 37 & 23 & 40 \\
512/22 & $1.2\!\cdot\! 10^{8}$ & $2.3\!\cdot\! 10^{5}$ & $1.9\!\cdot\! 10^{5}$/814 & 3713/121 & 62 & 15 & 114 \\
\end{tabular}
\caption{\label{table_glob_2D_level3} Strong scaling test, Poisson equation in 2D, mesh refined by prescribed refinements, 3-level BDDC.}
\end{center}
\end{table}


\begin{table}
\begin{center}
\begin{tabular}{ccc|cc|c|cc}
$N_S$ & $n$ & $n/N_S$ & $n^{\Gamma}$ & $n_{C}$ & its. & $t_{set-up}$ & $t_{PCG}$ \\
\hline
32 & $5.9\!\cdot\! 10^{7}$   & $1.9\!\cdot\! 10^{6}$ & $6.6\!\cdot\! 10^{5}$ & 423 & 33 & 508 & 136 \\
64 & $5.9\!\cdot\! 10^{7}$   & $9.3\!\cdot\! 10^{5}$ & $9.5\!\cdot\! 10^{5}$ & 949 & 41 & 269 & 92 \\
128 & $5.9\!\cdot\! 10^{7}$  & $4.6\!\cdot\! 10^{5}$ & $1.3\!\cdot\! 10^{6}$ & 2383 & 45 & 204 & 88 \\
256 & $5.9\!\cdot\! 10^{7}$  & $2.3\!\cdot\! 10^{5}$ & $1.8\!\cdot\! 10^{6}$ & 5232 & 47 & 103 & 52 \\
512 & $5.9\!\cdot\! 10^{7}$  & $1.2\!\cdot\! 10^{5}$ & $2.5\!\cdot\! 10^{6}$ & $1.1\!\cdot\! 10^{4}$ & 49 & 47 & 45 \\
1024 & $5.9\!\cdot\! 10^{7}$ & $5.8\!\cdot\! 10^{4}$ & $3.3\!\cdot\! 10^{6}$ & $2.3\!\cdot\! 10^{4}$ & 54 & 54 & 64 \\
\end{tabular}
\caption{\label{table_glob_3D_level2} Strong scaling test, Poisson equation in 3D, mesh refined by prescribed refinements, 2-level BDDC.}
\end{center}
\end{table}

\begin{table}
\begin{center}
\begin{tabular}{ccc|cc|c|cc}
$N_S$ & $n$ & $n/N_S$ & $n^{\Gamma}$ & $n_{C}$ & its. & $t_{set-up}$ & $t_{PCG}$ \\
\hline
32/6 & $5.9\!\cdot\! 10^{7}$    & $1.9\!\cdot\! 10^{6}$ & $6.6\!\cdot\! 10^{5}$/280 & 423/30 & 34 & 517 & 158 \\
64/8 & $5.9\!\cdot\! 10^{7}$    & $9.3\!\cdot\! 10^{5}$ & $9.5\!\cdot\! 10^{5}$/507 & 949/63 & 42 & 269 & 93 \\
128/12 & $5.9\!\cdot\! 10^{7}$  & $4.6\!\cdot\! 10^{5}$ & $1.3\!\cdot\! 10^{6}$/1624 & 2383/187 & 49 & 203 & 84 \\
256/16 & $5.9\!\cdot\! 10^{7}$  & $2.3\!\cdot\! 10^{5}$ & $1.8\!\cdot\! 10^{6}$/2514 & 5232/256 & 51 & 100 & 52 \\
512/22 & $5.9\!\cdot\! 10^{7}$  & $1.2\!\cdot\! 10^{5}$ & $2.5\!\cdot\! 10^{6}$/6110 & $1.1\!\cdot\! 10^{4}$/492 & 58 & 41 & 44 \\
1024/32 & $5.9\!\cdot\! 10^{7}$ & $5.8\!\cdot\! 10^{4}$ & $3.3\!\cdot\! 10^{6}$/9261 & $2.3\!\cdot\! 10^{4}$/572 & 67 & 21 & 47 \\
2048/46 & $5.9\!\cdot\! 10^{7}$ & $2.9\!\cdot\! 10^{4}$ & $4.4\!\cdot\! 10^{6}$/$1.3\!\cdot\! 10^{4}$ & $4.5\!\cdot\! 10^{4}$/748 & 71 & 14 & 52 \\
4096/64 & $5.9\!\cdot\! 10^{7}$ & $1.5\!\cdot\! 10^{4}$ & $5.9\!\cdot\! 10^{6}$/$2.0\!\cdot\! 10^{4}$ & $9.0\!\cdot\! 10^{4}$/1267 & 73 & 60 & 13 \\
\end{tabular}
\caption{\label{table_glob_3D_level3} Strong scaling test, Poisson equation in 3D, mesh refined by prescribed refinements, 3-level BDDC.}
\end{center}
\end{table}

\begin{table}
\begin{center}
\begin{tabular}{c|c|ccc|ccc|ccc}
\multirow{2}{*}{$N_S$} & \multirow{2}{*}{$n/N_S$} & \multicolumn{3}{c|}{num.\ coarse dofs} & \multicolumn{3}{c|}{time fact.\ loc.\ (s)} & \multicolumn{3}{c}{time sol.\ loc.\ (s)}\\ 
                       &                          & min & max & avg & min & max & avg & min & max & avg \\
\hline
16  & $7.5\!\cdot\! 10^{6}$ & 3 & 14 & 5 & 80.17 & 103.66 & 84.07 & 4.71 & 12.47 & 6.52 \\
32  & $3.7\!\cdot\! 10^{6}$ & 3 & 19 & 6 & 37.91 & 61.55 & 40.00 & 2.37 & 9.88 & 3.55 \\
64  & $1.9\!\cdot\! 10^{6}$ & 2 & 27 & 6 & 17.82 & 27.97 & 18.96 & 1.16 & 7.12 & 1.94 \\
128 & $9.3\!\cdot\! 10^{5}$ & 2 & 31 & 7 & 8.99 & 13.64 & 9.47 & 0.60 & 3.85 & 1.05 \\
256 & $4.7\!\cdot\! 10^{5}$ & 3 & 43 & 9 & 4.45 & 7.05 & 4.73 & 0.36 & 3.27 & 0.70 \\
512 & $2.3\!\cdot\! 10^{5}$ & 2 & 38 & 9 & 2.10 & 5.17 & 2.24 & 0.17 & 2.22 & 0.36 \\
\end{tabular}
\caption{\label{table_loc_2D} Strong scaling test, Poisson equation in 2D, mesh refined by prescribed refinements, local subdomain properties.}
\end{center}
\end{table}

\begin{table}
\begin{center}
\begin{tabular}{c|c|ccc|ccc|ccc}
\multirow{2}{*}{$N_S$} & \multirow{2}{*}{$n/N_S$} & \multicolumn{3}{c|}{num.\ coarse dofs} & \multicolumn{3}{c|}{time fact.\ loc.\ (s)} & \multicolumn{3}{c}{time sol.\ loc.\ (s)}\\ 
                       &                          & min & max & avg & min & max & avg & min & max & avg \\
\hline
32   & $1.9\!\cdot\! 10^{6}$ & 9 & 61 & 31 & 105.06 & 205.52 & 133.48 & 5.46 & 34.56 & 17.68 \\
64   & $9.3\!\cdot\! 10^{5}$ & 13 & 77 & 34 & 39.47 & 108.08 & 51.94 & 3.13 & 24.94 & 8.65 \\
128  & $4.6\!\cdot\! 10^{5}$ & 13 & 116 & 43 & 15.62 & 85.23 & 20.88 & 1.51 & 24.31 & 5.11 \\
256  & $2.3\!\cdot\! 10^{5}$ & 11 & 162 & 46 & 6.21 & 45.76 & 8.43 & 0.63 & 12.29 & 2.51 \\
512  & $1.2\!\cdot\! 10^{5}$ & 11 & 184 & 47 & 3.06 & 14.16 & 3.83 & 0.40 & 7.22 & 1.28 \\
1024 & $5.8\!\cdot\! 10^{4}$ & 11 & 238 & 48 & 1.29 & 5.00 & 1.59 & 0.17 & 5.84 & 0.61 \\
2048 & $2.9\!\cdot\! 10^{4}$ & 12 & 297 & 46 & 0.56 & 2.55 & 0.64 & 0.07 & 3.35 & 0.26 \\
4096 & $1.5\!\cdot\! 10^{4}$ & 12 & 224 & 45 & 0.25 & 0.60 & 0.29 & 0.03 & 0.69 & 0.12 \\
\end{tabular}
\caption{\label{table_loc_3D} Strong scaling test, Poisson equation in 3D, mesh refined by prescribed refinements, local subdomain properties.}
\end{center}
\end{table}

\subsubsection{Higher order elements}
\begin{figure}
\begin{center}
 \includegraphics[width = 0.49\textwidth]{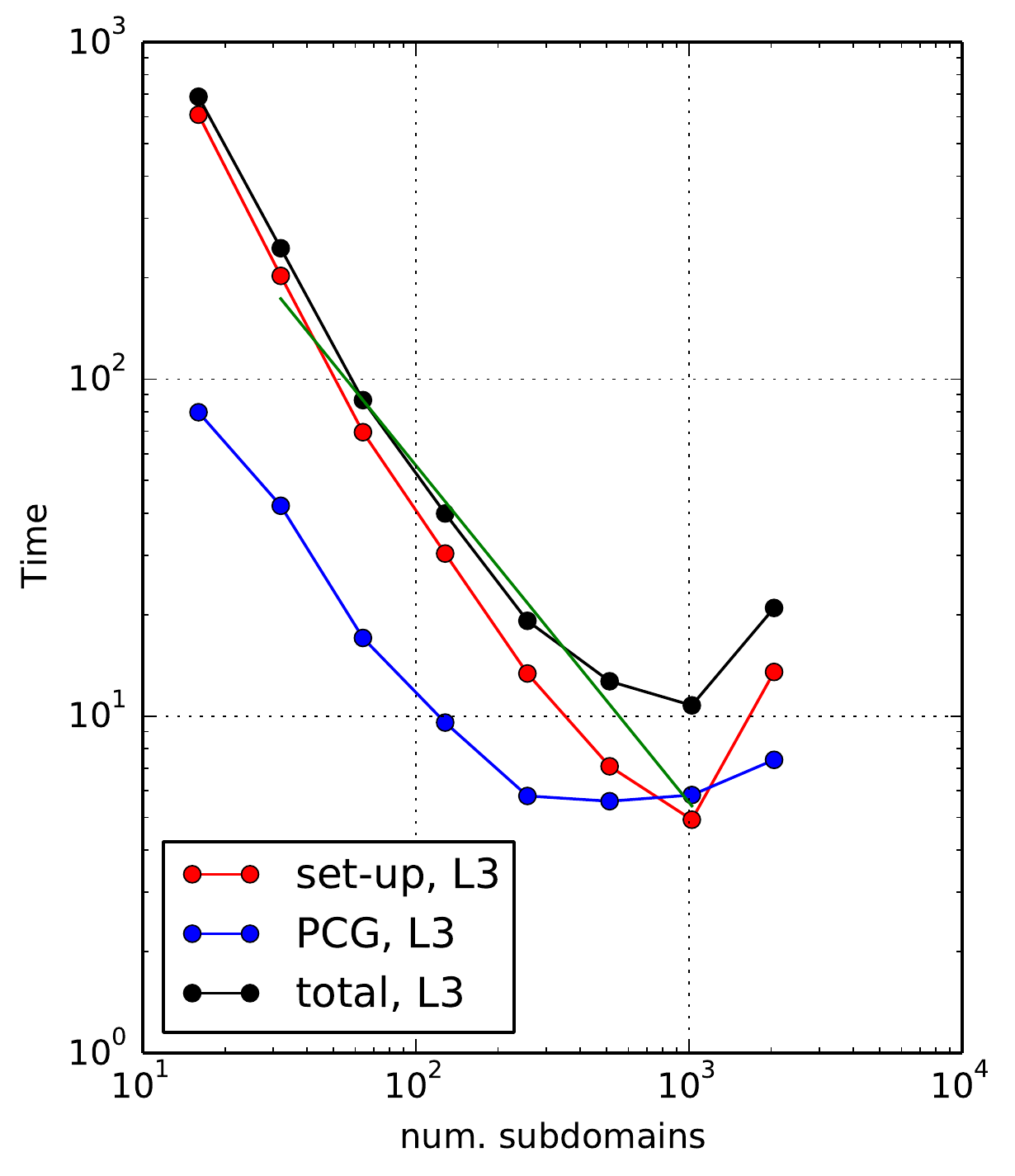}
 \includegraphics[width = 0.49\textwidth]{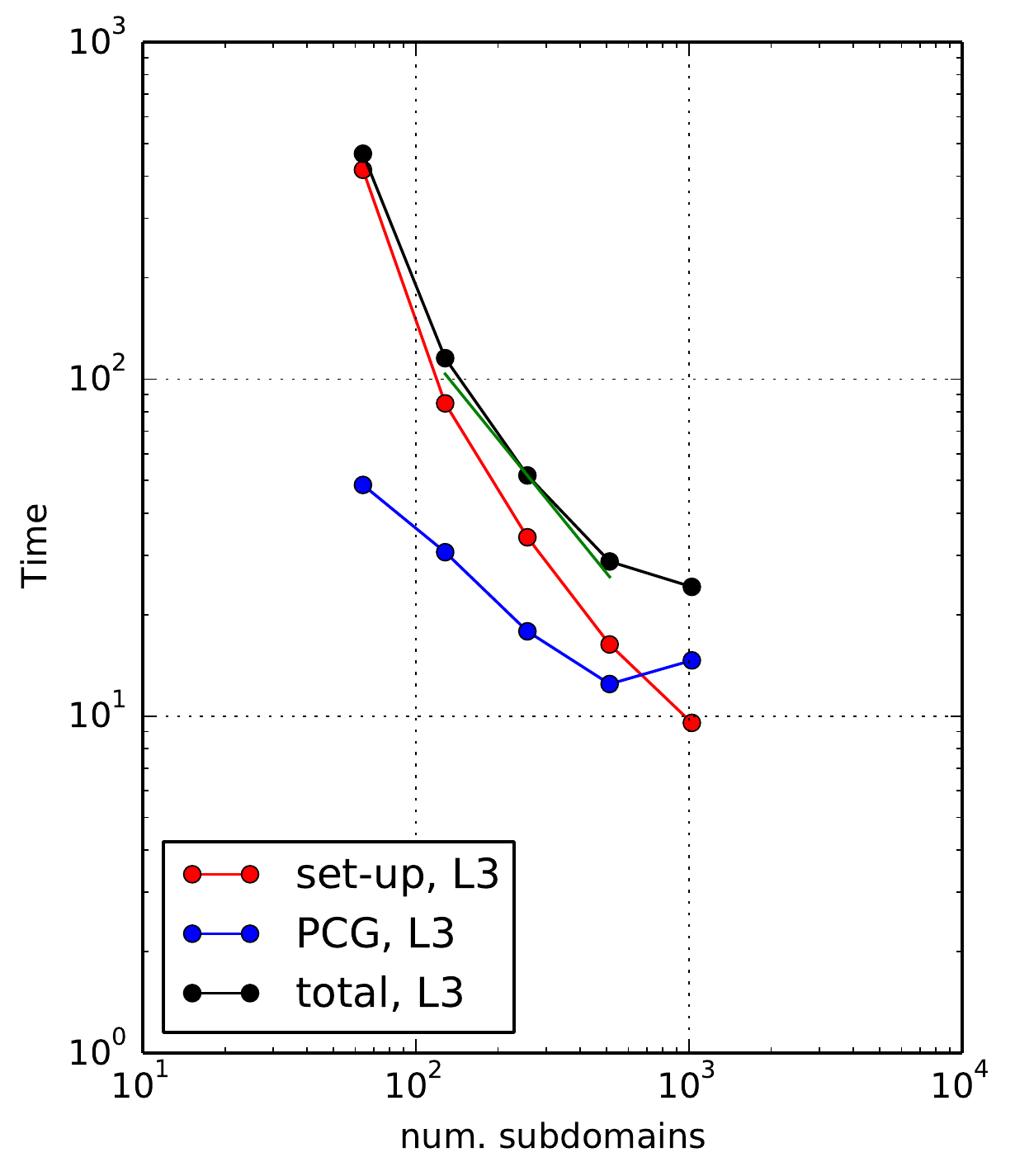}
 \caption{\label{fig:scaling_higher_order_and_elasticity}Poisson in 3D, polynomial order~4 (left) and linear elasticity, polynomial order~1 (right). Global parameters of the 
 runs, including depicted timings, can be found in Tables~\ref{table_glob_3D_order4} and~\ref{table_glob_elasitcity}, respectively. 
 }
\end{center}
\end{figure}

The previous cases of trilinear elements have shown good scaling properties in 3D, especially for the 3-level BDDC method.
In this section, we investigate the performance of the method for higher-order finite elements. 
\hl{Studies of using nonoverlapping DD methods for $p$ and $hp$ finite element methods were performed e.g. in~\cite{Mandel-1990-TDD,Toselli-2004-DDP}. For the BDDC method, convergence analysis with respect to the polynomial order $p$ was provided in~\cite{Klawonn-2008-SEF,Pavarino-2007-BFP} showing only mild worsening of the condition number with increasing $p$.}
We have successfully experimented with second and fourth order finite elements. 
However only results for the latter are presented in Tables~\ref{table_glob_3D_order4} and~\ref{table_loc_3D_order4} 
as the relevant effects are more pronounced for this choice.
Computational times are also plotted in the left part of Fig.~\ref{fig:scaling_higher_order_and_elasticity}.
Only results for the three-level BDDC method are presented.

The selected problem has been obtained by applying  
$R_U = 4$, $R_C = 4$, and $R_S = 4$ refinements to the unit cube. It has 336 thousands of elements and nearly 20 million unknowns. 
We use the averaging scaled by the diagonal stiffness in the weight matrices $W_i$ in the BDDC algorithm (see e.g. \cite{Certikova-2015-DAI} and references therein), 
namely in equations (\ref{eq:making_local_residual}) and (\ref{eq:averaging_subdomain_solves}).
Note that simple arithmetic average has been used in the previous cases.
Nevertheless, this provided us with a significant improvement of convergence for the higher-order elements. 

As we can see from Table~\ref{table_glob_3D_order4}, increasing the polynomial order of finite elements does not have 
any significant effect on the global number of constraints.
More importantly, the convergence of the BDDC method does not seem to deteriorate with the polynomial order. 
In fact, one can see even a slight relative improvement compared to the linear case in Table~\ref{table_glob_3D_level3}.
The number of iterations approximately doubles when going from 32 to 2048 subdomains, which is similar for both cases and seems acceptable.

With respect to local subdomain properties in Table~\ref{table_loc_3D_order4}, the increase in polynomial order approximately doubles the 
times for factorization of the subdomain problem, while the issues with load balance of the coarse basis functions solution remain similar.
Note that this conclusion is based on comparing the case for 1024 subdomains of fourth order elements and 4096 subdomains with the linear elements to get 
comparable subdomain sizes.

\begin{table}
\begin{center}
\begin{tabular}{ccc|cc|c|cc}
$N_S$ & $n$ & $n/N_S$ & $n^{\Gamma}$ & $n_{C}$ & its. & $t_{set-up}$ & $t_{PCG}$ \\
\hline
32/6 & $2.0\!\cdot\! 10^{7}$    & $6.2\!\cdot\! 10^{5}$ & $5.2\!\cdot\! 10^{5}$/308 & 473/37 & 24 & 202 & 42 \\
64/8 & $2.0\!\cdot\! 10^{7}$    & $3.1\!\cdot\! 10^{5}$ & $7.3\!\cdot\! 10^{5}$/531 & 997/53 & 26 & 70 & 17 \\
128/12 & $2.0\!\cdot\! 10^{7}$  & $1.5\!\cdot\! 10^{5}$ & $9.6\!\cdot\! 10^{5}$/1427 & 2106/205 & 31 & 30 & 9.6 \\
256/16 & $2.0\!\cdot\! 10^{7}$  & $7.7\!\cdot\! 10^{4}$ & $1.3\!\cdot\! 10^{6}$/2556 & 4690/206 & 37 & 13 & 5.8 \\
512/22 & $2.0\!\cdot\! 10^{7}$  & $3.9\!\cdot\! 10^{4}$ & $1.7\!\cdot\! 10^{6}$/5899 & 9561/437 & 45 & 7.1 & 5.6 \\
1024/32 & $2.0\!\cdot\! 10^{7}$ & $1.9\!\cdot\! 10^{4}$ & $2.1\!\cdot\! 10^{6}$/9459 & $1.9\!\cdot\! 10^{4}$/708 & 52 & 4.9 & 5.8 \\
2048/46 & $2.0\!\cdot\! 10^{7}$ & 9658                  & $2.7\!\cdot\! 10^{6}$/$1.2\!\cdot\! 10^{4}$ & $3.9\!\cdot\! 10^{4}$/919 & 52 & 14 & 7.4 \\
\end{tabular}
\caption{\label{table_glob_3D_order4} Strong scaling test, Poisson equation in 3D, mesh refined by prescribed refinements, 4th order elements, 3-level BDDC.}
\end{center}
\end{table}

\begin{table}
\begin{center}
\begin{tabular}{c|c|ccc|ccc|ccc}
\multirow{2}{*}{$N_S$} & \multirow{2}{*}{$n/N_S$} & \multicolumn{3}{c|}{num.\ coarse dofs} & \multicolumn{3}{c|}{time fact.\ loc.\ (s)} & \multicolumn{3}{c}{time sol.\ loc.\ (s)}\\ 
                       &                          & min & max & avg & min & max & avg & min & max & avg \\
\hline
32   & $6.2\!\cdot\! 10^{5}$ & 10 & 65 & 34 & 49.51 & 84.43 & 63.92 & 2.91 & 15.25 & 7.83 \\
64   & $3.1\!\cdot\! 10^{5}$ & 13 & 80 & 36 & 16.26 & 26.64 & 21.15 & 1.53 & 7.65 & 3.18 \\
128  & $1.5\!\cdot\! 10^{5}$ & 12 & 110 & 38 & 6.92 & 11.24 & 8.80 & 0.56 & 4.28 & 1.53 \\
256  & $7.7\!\cdot\! 10^{4}$ & 11 & 141 & 42 & 2.49 & 4.36 & 3.50 & 0.23 & 2.15 & 0.75 \\
512  & $3.9\!\cdot\! 10^{4}$ & 15 & 164 & 43 & 1.21 & 2.01 & 1.55 & 0.15 & 1.15 & 0.37 \\
1024 & $1.9\!\cdot\! 10^{4}$ & 9  & 129 & 44 & 0.55 & 0.84 & 0.65 & 0.06 & 0.48 & 0.17 \\
2048 & 9658                  & 12 & 129 & 45 & 0.24 & 0.52 & 0.29 & 0.03 & 0.23 & 0.08 \\
\end{tabular}
\caption{\label{table_loc_3D_order4} Strong scaling test, Poisson equation in 3D, mesh refined by prescribed refinements, 4th order elements, local subdomain properties.}
\end{center}
\end{table}

\subsubsection{Linear elasticity}
In this section, we investigate the performance of the method for a different problem, 
namely problem of the linear elasticity defined as 
\begin{equation}
\label{eq:elasticity}
(\lambda + \mu )\ \mbox{grad}\ \mbox{div}\ {\bf u} + \mu \Delta {\bf u} = {\bf F}\ \mbox{on}\ \Omega, 
\quad  {\bf u} = {\bf 0}\ \mbox{on}\ \partial \Omega,\quad \Omega = [0,1]^3,
\end{equation}
where $\lambda$ and $\mu$ are the Lam{\' e} constants, ${\bf u}$ is the displacement vector, 
and ${\bf F}$ represent body forces. The Lam{\' e} constants are related to 
the Young's modulus $E$ and the Poisson's ratio $\nu$ as
\begin{equation}
 \lambda = \frac{\nu E}{(1+\nu)(1-2\nu)},\qquad \mu = \frac{E}{2(1+\nu)}.
\end{equation}
For our calculations, we use $E = 10^{10}$, $\nu=1/3$, ${\bf F} = (0, 0, -10^5)$ together with zero Dirichlet boundary 
condition on the whole $\partial\Omega$. 
Solution to this benchmark problem corresponds to displacement caused by 
application of a constant body force to a cube supported at all faces.

We can see two interesting differences from the Poisson problem that are of interest for finite element solution and a domain decomposition method.
Firstly, the structure of the matrices is more complicated, with more non-zero entries on each row. 
This makes the problem more complicated for a sparse direct solver used on subdomains.

Nevertheless, the main difference is the more complicated structure of the nullspace of local subdomain matrices. 
While it is simply the one-dimensional space of constants for the Poisson problem, it is a six-dimensional space for a floating subdomain for the elasticity problem.
This space is spanned by the rigid body motions, which are three translations and three rotations around the given coordinate axes.
However, if the subdomain consists of several components that are completely disconnected or connected only through a single finite element node, 
the nullspace structure becomes more complicated. 
In general, the dimension of the local nullspace can be anywhere between six times the number of subdomain components and zero, 
depending also on the fact whether the subdomain is in touch with boundary conditions.

Three options seem important in the set-up of the BDDC preconditioner in this case. 
First, we consider dual graph of subdomain mesh (elements are graph vertices) to detect components. 
Since a graph edge is created only when two elements share a face, two components will be detected even if the elements are connected via 
a single node or a single edge, and coarse degrees of freedom are generated independently for each of them.
Next, we consider selection of corners for this case. 
This was shown to play an important role for complex elasticity problems solved on unstructured meshes in \cite{Sistek-2012-FSC},
and the algorithm from this reference is used in \codename{BDDCML}.
Finally, we use again the averaging scaled by the diagonal stiffness in the weight matrices $W_i$. 
All these three modifications had a positive effect on convergence for the elasticity problem.

The numerical tests presented in Table~{\ref{table_glob_elasitcity}} and in the right part of Fig.~\ref{fig:scaling_higher_order_and_elasticity} are obtained for problem (\ref{eq:elasticity}) solved on a unit cube
after application of $R_U = 7$, $R_C = 3$, and $R_S = 3$ refinements, leading to 7.3 million elements and 16.7 million unknowns.
%
The results of a strong scaling test are presented for the 3-level BDDC method in Table~{\ref{table_glob_elasitcity}}, 
and the solution times are also plotted in the right part of Fig.~\ref{fig:scaling_higher_order_and_elasticity} along with those for a 2-level method.
We can see again a certain increase in the number of iterations during the strong scaling test, 
starting at 32 for 64 subdomains and reaching 79 iterations for 1024 subdomains.
We consider this again to be rather acceptable for the given highly refined unstructured meshes. 
Despite this growth of the number of iterations, the times in the plot show a rather nice strong scaling properties. 
With the 2-level BDDC method, we were not able to obtain results with more than 128 subdomains.

Let us conclude this study by looking at the local subdomain properties presented in Table~\ref{table_loc_elasticity},
and compare them to Table~\ref{table_loc_3D}.
To look at comparable subdomain sizes, one can consider the case of 1024 subdomains for the elasticity problem and the case with 
4096 subdomains for the Poisson problem.
We can see that the number of local constraints approximately triples for the elasticity problem compared to the Poisson problem, however, 
the ratio of maximum and minimum number of local constraints is smaller, only around 10 for linear elasticity. 
The factorization time increases slightly more than twice, and the ratio also slightly worsens. 
As we have already discovered, the worst part in terms of scaling with number of constraints is the back-substitution part 
which, due to the high number of local constrains, is now more than twice more expensive than the factorization.

\begin{table}
\begin{center}
\begin{tabular}{ccc|cc|c|cc}
$N_S$ & $n$ & $n/N_S$ & $n^{\Gamma}$ & $n_{C}$ & its. & $t_{set-up}$ & $t_{PCG}$ \\
\hline
64/8 & $1.7\!\cdot\! 10^{7}$    & $2.6\!\cdot\! 10^{5}$ & $9.1\!\cdot\! 10^{5}$/2526 & 1282/108 & 32 & 417 & 49 \\
128/12 & $1.7\!\cdot\! 10^{7}$  & $1.3\!\cdot\! 10^{5}$ & $1.2\!\cdot\! 10^{6}$/6759 & 2815/274 & 49 & 85 & 31 \\
256/16 & $1.7\!\cdot\! 10^{7}$  & $6.5\!\cdot\! 10^{4}$ & $1.6\!\cdot\! 10^{6}$/$1.2\!\cdot\! 10^{4}$ & 5777/369 & 61 & 34 & 18 \\
512/22 & $1.7\!\cdot\! 10^{7}$  & $3.3\!\cdot\! 10^{4}$ & $2.1\!\cdot\! 10^{6}$/$2.6\!\cdot\! 10^{4}$ & $1.2\!\cdot\! 10^{4}$/728 & 66 & 16 & 12 \\
1024/32 & $1.7\!\cdot\! 10^{7}$ & $1.6\!\cdot\! 10^{4}$ & $2.7\!\cdot\! 10^{6}$/$3.8\!\cdot\! 10^{4}$ & $2.4\!\cdot\! 10^{4}$/919 & 79 & 9.5 & 15 \\
\end{tabular}
\caption{\label{table_glob_elasitcity} Strong scaling test, linear elasticity in 3D, mesh refined by prescribed refinements, trilinear elements, 3-level BDDC.}
\end{center}
\end{table}

\begin{table}
\begin{center}
\begin{tabular}{c|c|ccc|ccc|ccc}
\multirow{2}{*}{$N_S$} & \multirow{2}{*}{$n/N_S$} & \multicolumn{3}{c|}{num.\ coarse dofs} & \multicolumn{3}{c|}{time fact.\ loc.\ (s)} & \multicolumn{3}{c}{time sol.\ loc.\ (s)}\\ 
                       &                          & min & max & avg & min & max & avg & min & max & avg \\
\hline
64   & $2.6\!\cdot\! 10^{5}$ & 75 & 342 & 191 & 13.68 & 186.26 & 31.21 & 6.32 & 88.89 & 21.36 \\
128  & $1.3\!\cdot\! 10^{5}$ & 81 & 441 & 211 & 5.55 & 31.63 & 11.47 & 3.74 & 29.86 & 10.66 \\
256  & $6.5\!\cdot\! 10^{4}$ & 75 & 537 & 218 & 2.36 & 11.80 & 4.27 & 1.25 & 13.72 & 4.78 \\
512  & $3.3\!\cdot\! 10^{4}$ & 81 & 621 & 225 & 0.96 & 3.80 & 1.70 & 0.72 & 8.00 & 2.33 \\
1024 & $1.6\!\cdot\! 10^{4}$ & 69 & 633 & 229 & 0.38 & 1.47 & 0.65 & 0.25 & 3.77 & 1.05 \\
\end{tabular}
\caption{\label{table_loc_elasticity} Strong scaling test, linear elasticity in 3D, mesh refined by prescribed refinements, trilinear elements, local subdomain properties.}
\end{center}
\end{table}

\subsection{Impact of non-uniform mesh}
\begin{figure}
\centering
 \includegraphics[height=0.32\textwidth]{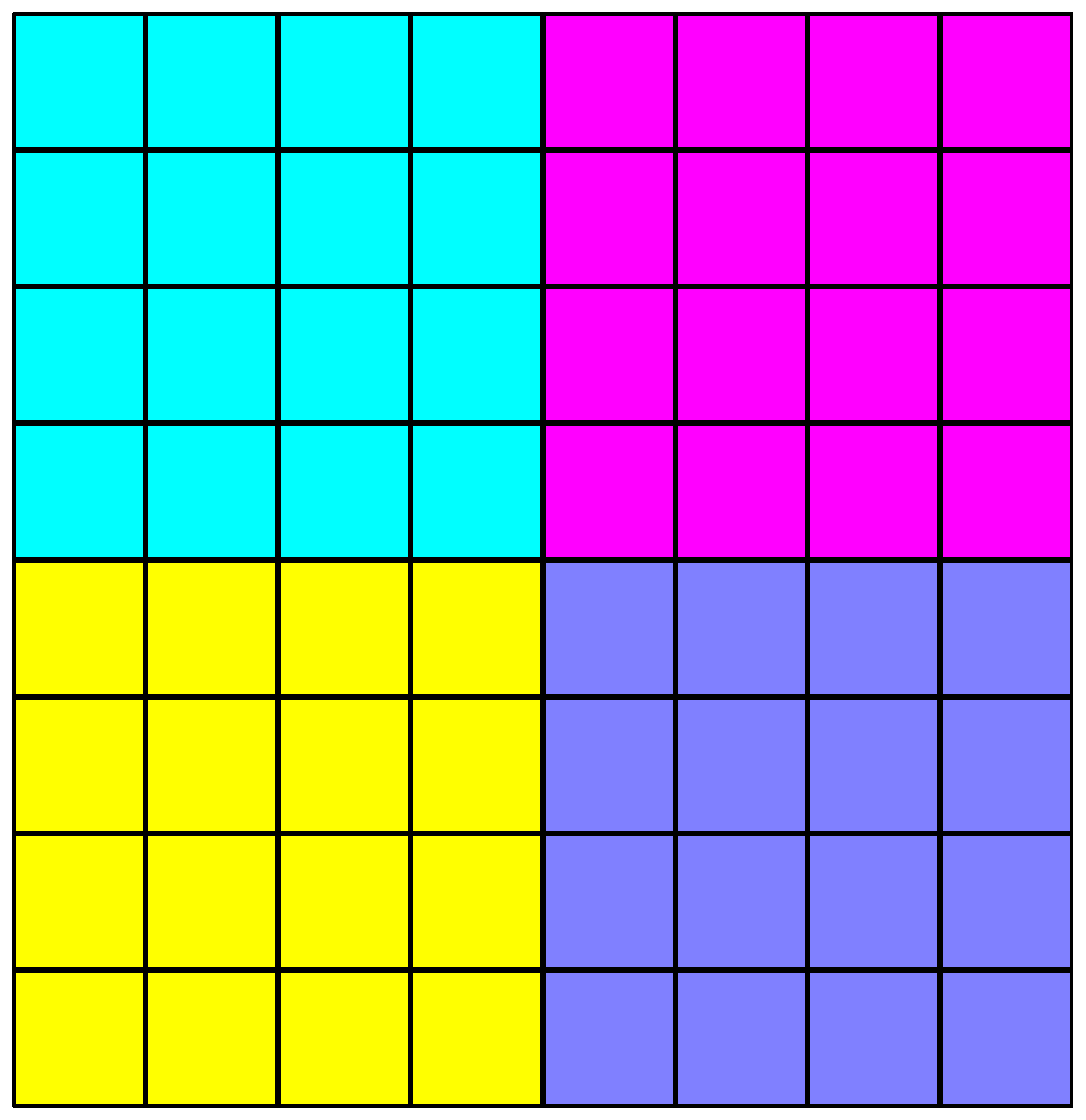} 
 \includegraphics[height=0.32\textwidth]{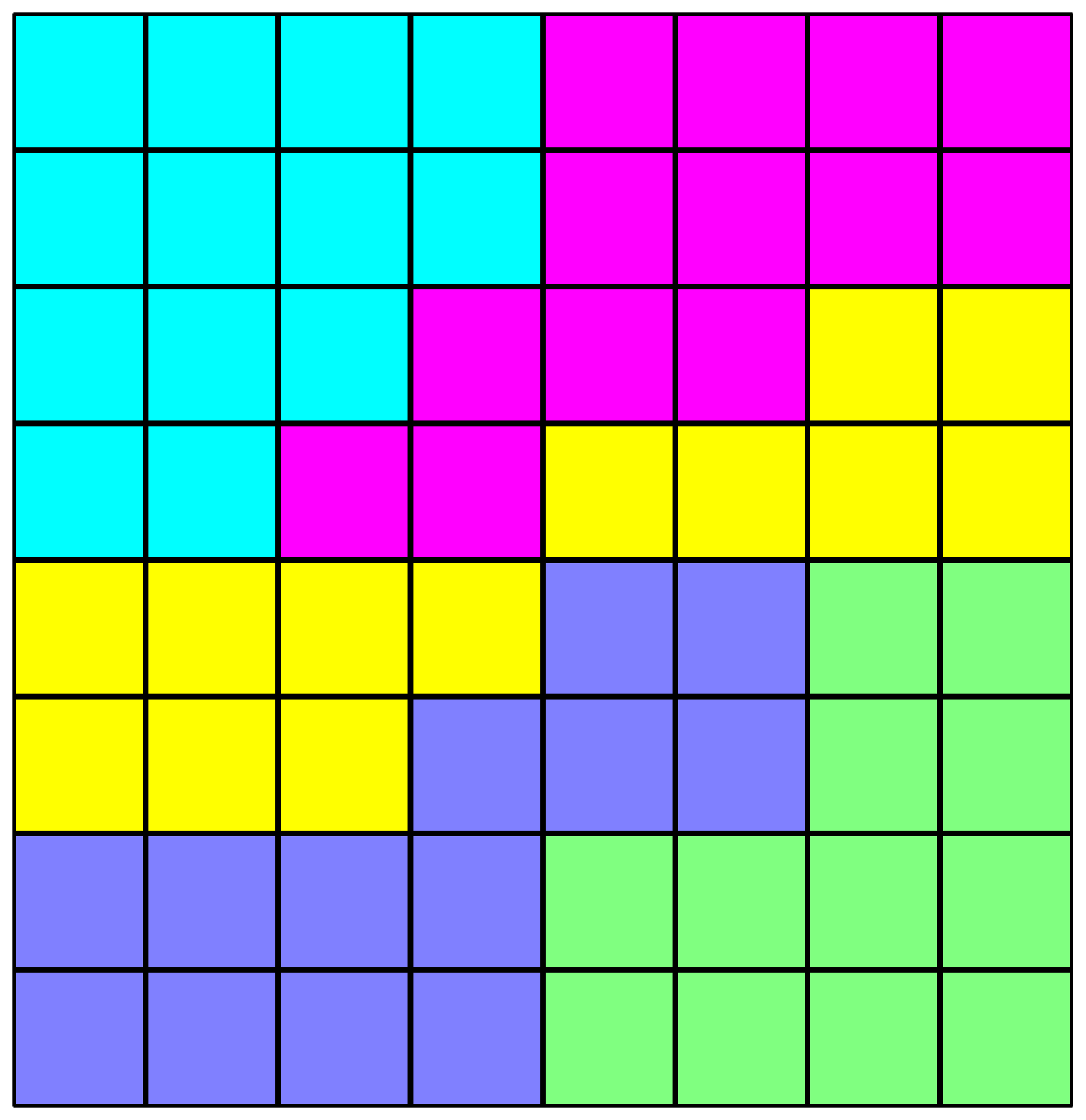} 
 \includegraphics[height=0.32\textwidth]{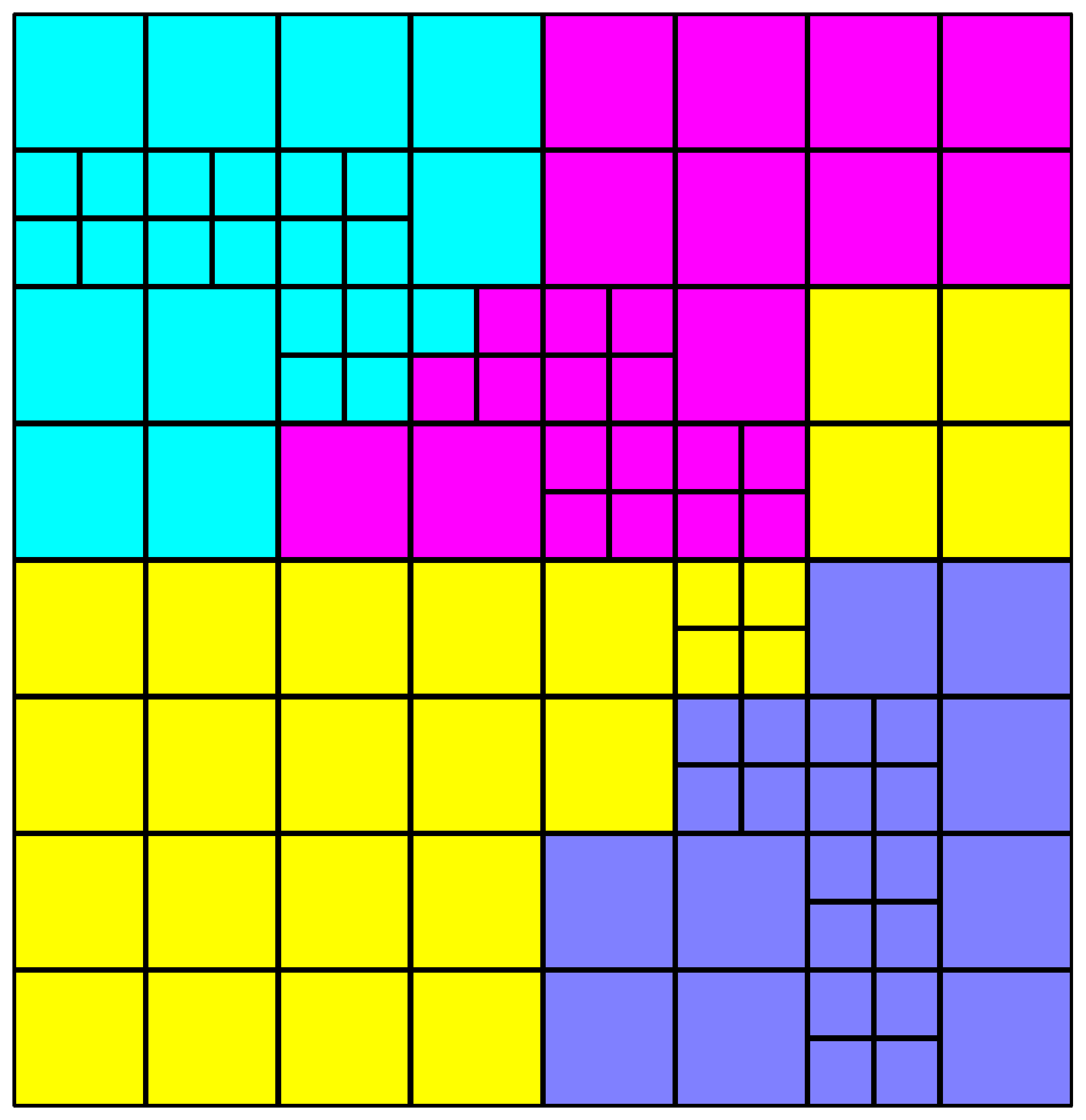} 
 \caption{Illustrative examples (using few elements in 2D) of the difference between uniform mesh with uniform subdomains (left), the same mesh with nonaligned subdomains (centre) and adapted mesh (right). Note that the difference between the first and the second case is only in the number of subdomains, which is not aligned with the mesh in the second case. Figures correspond to the three cases shown in Tables~\ref{table_glob_structure} and~\ref{table_loc_structure}.}
\label{fig:structure}
\end{figure}

\hl{In the previous section,} we have shown good scaling properties of the method for different 
types of elements and different configurations. 
In this section, we want to summarize previous results and perform a simple experiment to assess the influence of partitioning by 
the \codename{p4est} library on the performance of the domain decomposition solver. For this, we 
present properties of three calculations. 
The first one is obtained by 8 consecutive applications of $R_U$ refinement. 
It is thus a cube uniformly refined into $(2^8)^3$ elements. When we use \codename{p4est} to partition this mesh into 512 processors, 
we obtain regular subdomains, as it is shown for a smaller problem in 2D at the left panel of Fig.~\ref{fig:structure}. 
It is due to the fact that the number of elements is divisible by the number of subdomains. 
The second configuration we consider is the same mesh divided into 513 processors. 
Since now the number of elements is not divisible by number of processors, we obtain more complicated substructures, as it is suggested in the centre of Fig.~\ref{fig:structure}. 
The last considered configuration is for the same mesh with one additional $R_C$ and one $R_S$ refinement. 
We partition this mesh again into 512 processors, but this time the mesh is not uniform anymore, and it does not lead to a partition with regular subdomains.

The data obtained from these three calculations are shown in Tables~\ref{table_glob_structure} and~\ref{table_loc_structure}. 
We can observe that the irregularity of subdomains and adaptively refined meshes have similar role in increasing the complexity of the problem and the number of iterations,
and the number of iterations approximately doubles from the case of uniform mesh with regular subdomains.

\begin{table}
\begin{center}
\begin{tabular}{ccc|cc|c|cc}
$N_S$ & $n$ & $n/N_S$ & $n^{\Gamma}$ & $n_{C}$ & its. & $t_{set-up}$ & $t_{PCG}$ \\
\hline
512/22 & $1.70\!\cdot\! 10^{7}$ & $3.3\!\cdot\! 10^{4}$ & $1.3\!\cdot\! 10^{6}$/1586 & 2863/252 & 16 & 3.4 & 2.1 \\
513/22 & $1.70\!\cdot\! 10^{7}$ & $3.3\!\cdot\! 10^{4}$ & $1.8\!\cdot\! 10^{6}$/3530 & 8197/303 & 22 & 4.3 & 3 \\
512/22 & $1.74\!\cdot\! 10^{7}$ & $3.4\!\cdot\! 10^{4}$ & $1.8\!\cdot\! 10^{6}$/5028 & 7745/530 & 30 & 4.7 & 3.9 \\
\end{tabular}
\caption{\label{table_glob_structure}Impact of nonuniformity. Similar cases (but much larger and 3D) as shown in Fig.~\ref{fig:structure}, in the same order. 
Both using unaligned subdomains and unstructured meshes 
contribute to the higher complexity and the corresponding increase of the
number of iterations and both set-up and solution times.}
\end{center}
\end{table}

\begin{table}
\begin{center}
\begin{tabular}{c|c|ccc|ccc|ccc}
\multirow{2}{*}{$N_S$} & \multirow{2}{*}{$n/N_S$} & \multicolumn{3}{c|}{num.\ coarse dofs} & \multicolumn{3}{c|}{time fact.\ loc.\ (s)} & \multicolumn{3}{c}{time sol.\ loc.\ (s)}\\ 
                       &                          & min & max & avg & min & max & avg & min & max & avg \\
\hline
512/22 & $3.3\!\cdot\! 10^{4}$ & 6 & 18 & 14 & 1.23 & 1.36 & 1.28 & 0.06 & 0.19 & 0.15 \\
513/22 & $3.3\!\cdot\! 10^{4}$ & 8 & 62 & 39 & 0.88 & 1.32 & 1.03 & 0.07 & 0.48 & 0.32 \\
512/22 & $3.4\!\cdot\! 10^{4}$ & 12 & 63 & 36 & 0.78 & 1.43 & 1.08 & 0.13 & 0.48 & 0.31 \\
\end{tabular}
\caption{\label{table_loc_structure}Impact of nonuniformity, local subdomain properties, the same three meshes as in Table~\ref{table_glob_structure}.}
\end{center}
\end{table}

\subsection{Parallel adaptivity}
\label{sec:Numerical_adapt}

We have shown so far that the BDDC method preserves good scaling properties for significantly refined meshes with hanging nodes.
Those meshes were created artificially, by applying prescribed refinements, 
which were designed to mimic the expected behaviour of an adaptive algorithm for certain type of problems.

\begin{figure}[tbh]
\centering
 \includegraphics[width = 0.33\textwidth]{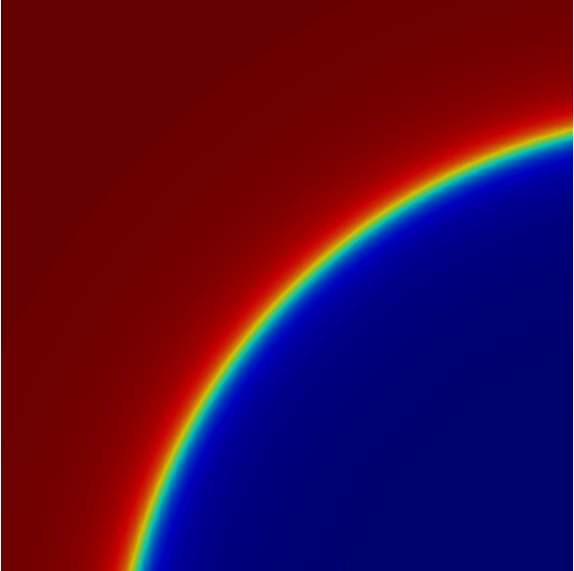}
 \hskip 10mm
 \includegraphics[width = 0.33\textwidth]{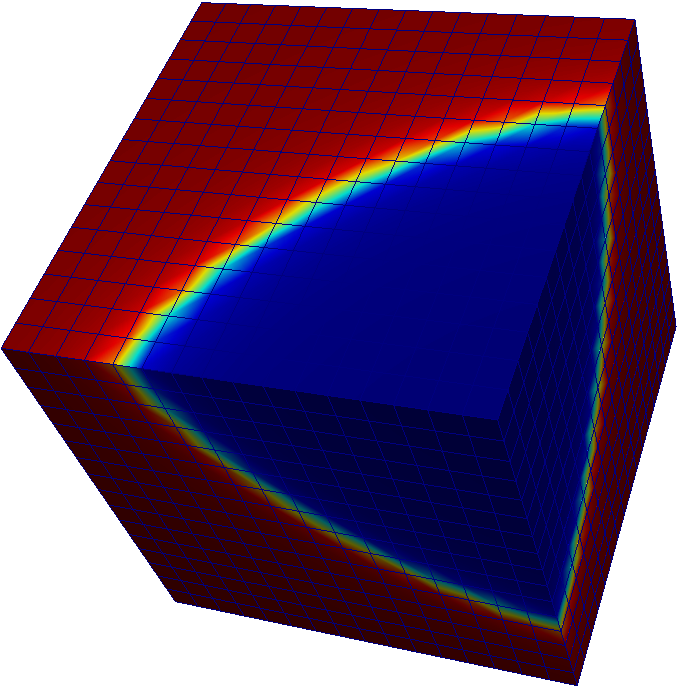} 
 \caption{\label{fig:exact_solution}Exact solution of the test problem in 2D and 3D. } 
\end{figure}

In this section, we show results of adaptive runs of the described algorithm, with no refinements prescribed a priori. 
Instead, the adaptivity algorithm described in Section~\ref{sec:modified_refinement} is used. 
We use again the Poisson problem defined in Section~\ref{sec:Numerical}. 
However, now we prescribe the exact solution as
\begin{equation}
   u^{*} = \arctan \Bigl(s \cdot \bigl(r-\frac{\pi}{3}\bigr)\Bigr),
\end{equation}
where $r$ is the distance from the point $[1.25, -0.25]$ in 2D or $[1.25, -0.25, -0.25]$ in 3D. 
The right-hand side $f$ and Dirichlet boundary condition $g$ are defined to recover the prescribed solution, 
\begin{equation}
   f = - \Delta u^{*},\quad g = u^{*}|_{\partial\Omega}.
\end{equation}
Function $u^{*}$ exhibits an internal layer\hl{, i.e. it has large gradients on a narow region inside the domain,} as can be seen in Fig.~\ref{fig:exact_solution} for 2D and 3D. 
The steepness of the layer is controlled by parameter $s$, which is set to $s=60$ in presented calculations.
The knowledge of the exact solution makes it easy to estimate the error of the current solution for each element and to select
elements for refinements accordingly. 
This is not possible for real-life applications though, where a different strategy has to be employed.
For example, one can use one of the existing a posteriori error estimators.

The parameters of the algorithm described in detail in Section~\ref{sec:modified_refinement}
are set as follows. In each adaptivity step we refine $\zeta = 15 \%$ of elements in linear case and $\zeta = 12 \%$ of elements otherwise. In 3D, when each refined element is replaced by eight new elements, this leads to roughly doubling the number of elements after each step. The change of the number of degrees of freedom is similar. 
We have also experimented with several values of the number of compartments $M$ from the algorithm of Section~\ref{sec:modified_refinement}
and found that for our setting with 2048 subdomains and meshes with more than 10$^9$ elements, the value $M=100$ is sufficient and its increase does not significantly change the behaviour of the algorithm. 

\begin{figure}
\begin{center}
 \includegraphics[width = 1\textwidth]{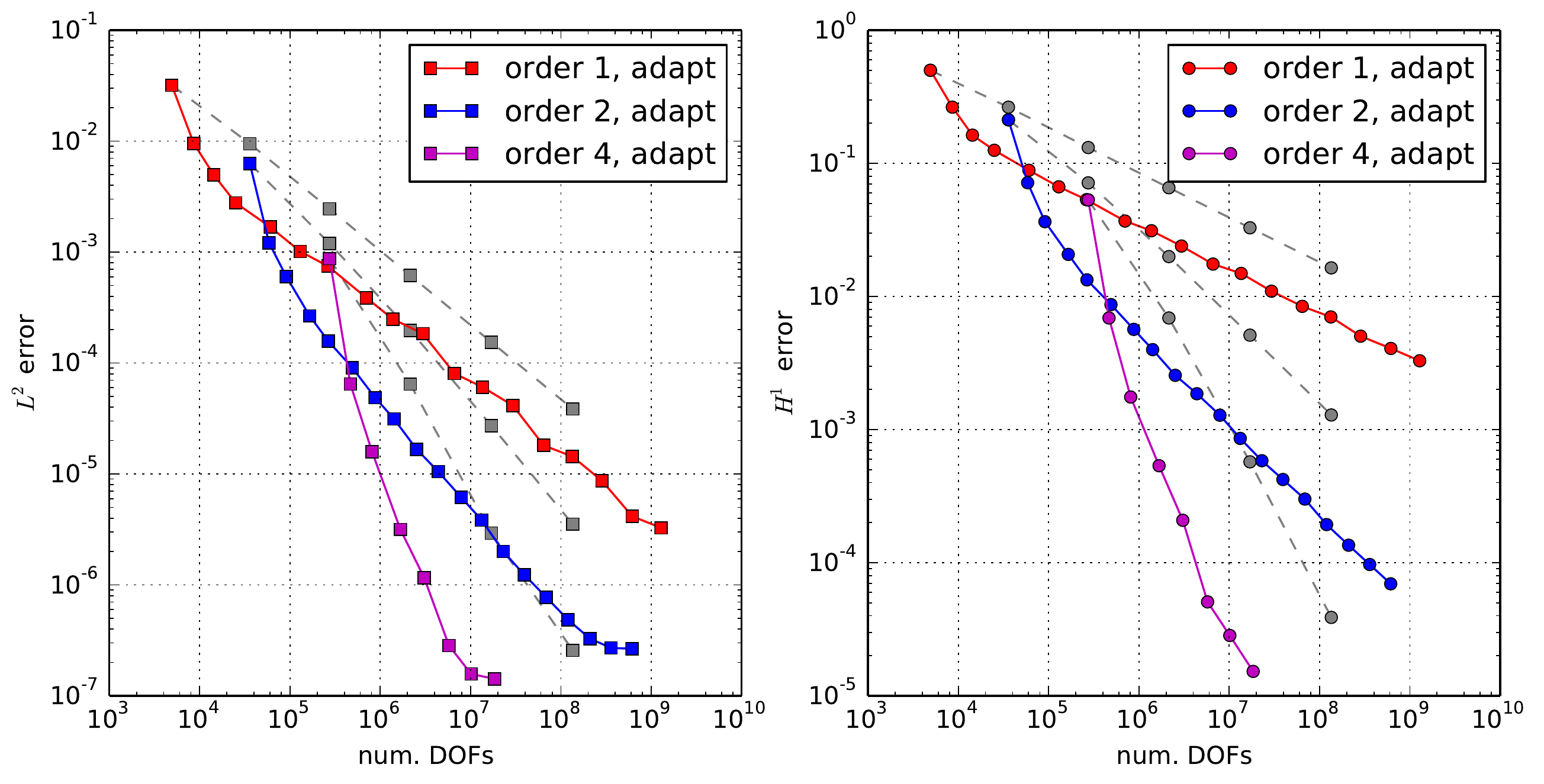}
  \caption{\label{fig:convergence}Convergence of the $L^2$ and $H^1$ norms of the error, respectively, on 2048 subdomains in 3D, 3-level BDDC method.
Different colours indicate adaptive method with different polynomial orders.
Grey colour indicates convergence obtained by corresponding uniform refinements. In both cases, polynomial orders 1, 2, and 4 are used.} 
\end{center}
\end{figure}

The error is measured as the norm of the difference between the approximate and the exact solution. 
Size of the error in each adaptivity step is depicted in Fig.~\ref{fig:convergence} together with the error for uniform refinements (all elements refined in each step), \hl{leading to an $8\times$ increase of the number of elements}.

\begin{figure}
 \includegraphics[width = 0.32\textwidth]{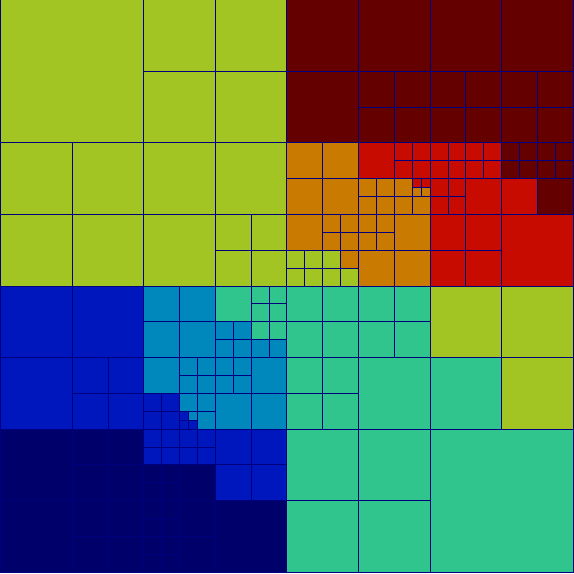}
 \includegraphics[width = 0.32\textwidth]{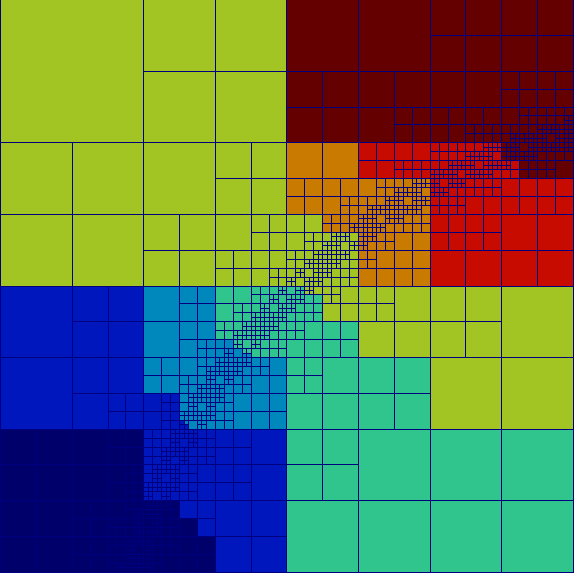}
 \includegraphics[width = 0.32\textwidth]{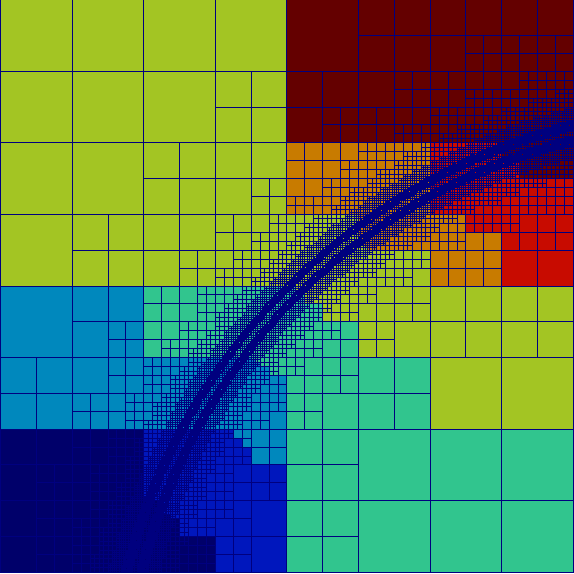}
 \caption{\label{fig:adapt_meshes_2D}Three illustrative adaptivity steps in 2D.}
\end{figure}

\begin{figure}
 \includegraphics[width = 0.32\textwidth]{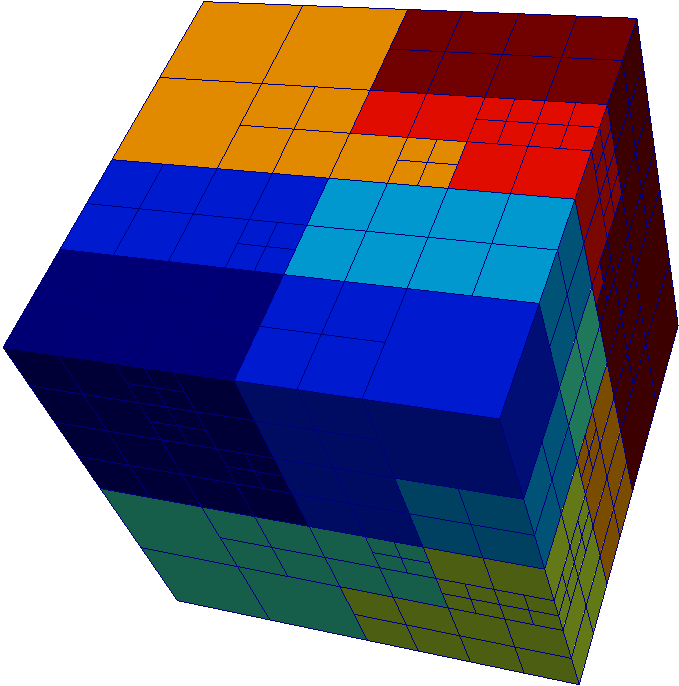}
 \includegraphics[width = 0.32\textwidth]{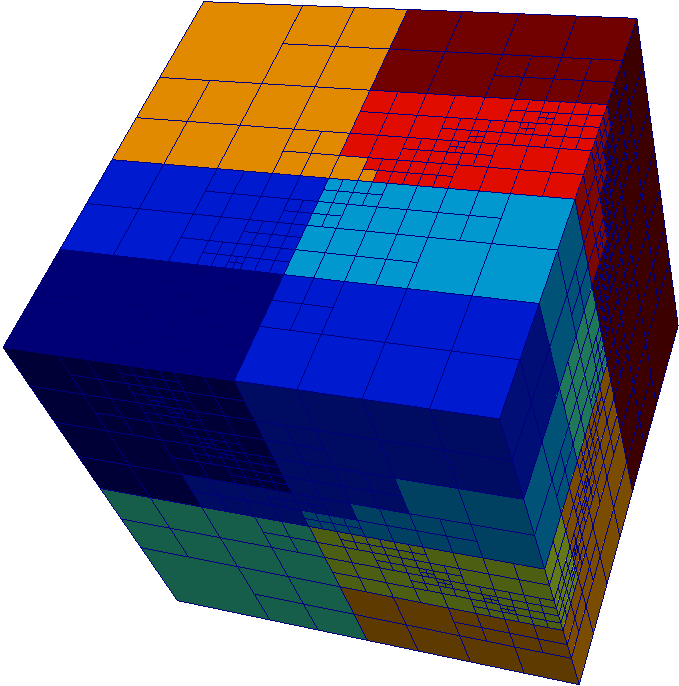}
 \includegraphics[width = 0.32\textwidth]{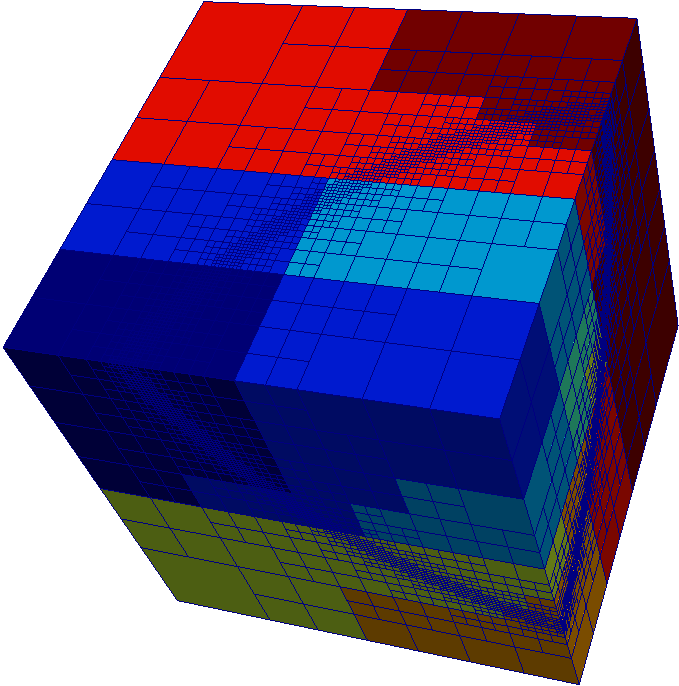}
 \caption{\label{fig:adapt_meshes_3D}Three illustrative adaptivity steps in 3D.}
\end{figure}

The adaptivity algorithm behaves according to the expectations and refines elements close to the internal layer, 
as can be seen in Fig.~\ref{fig:adapt_meshes_2D} for the 2-D and in Fig.~\ref{fig:adapt_meshes_3D} for the 3-D problem.
Different colours indicate different subdomains (processors). 
The figures are only for illustration and were obtained from runs on 8 processors. 
Notice how some areas of the computational domains change colours, indicating rebalancing of the mesh and moving some elements from one processor to another. 

Convergence of the global error of the solution is shown in Fig.~\ref{fig:convergence}, 
in $L^2$ and $H^1$ norms in 3D and for orders 1, 2, and 4. 
We show convergence of the adaptive algorithm (individual colours) in comparison with
the corresponding uniform refinements.
As can be observed, the use of adaptive algorithm achieves higher precision with much fewer degrees of freedom than the uniform refinements. 
To be more precise, let us suppose that we want to achieve $H^1$ error around $10^{-3}$. Using linear elements (no matter if uniform or adaptive), this precision has not been reached even for the largest computed problems. Using quadratic elements with uniform refinements, we need more than $2\cdot 10^8$ degrees of freedom. 
If we use quadratic elements with adaptivity (or fourth order elements with uniform refinements), approx.\ $10^7$ degrees of freedom are sufficient. 
Finally, using adaptivity with fourth order elements, the same error is achieved using only approx.\ $10^6$ degrees of freedom. In the same way, one can fix the number of degrees of freedom and find reduction of error when higher-order and/or adaptive algorithm is used. For example, using $10^7$ degrees of freedom, one can reduce the $L^2$ error 
from more than $2\cdot 10^{-4}$ for linear uniform refinements to less then $2\cdot 10^{-7}$ for adaptivity with order 4. The reduction of error is thus more then $1000\times$ by using a better method, while keeping the final linear system of the same size (although the matrix might be significantly denser).

Let us now analyze the difference between adaptive and corresponding uniform convergence curves (grey with the same starting 
point in Fig.~\ref{fig:convergence}).
At the beginning, the adaptive algorithm reduces the error very quickly compared to 
uniform refinements. After several steps, which refine elements close to the internal layer, the error is approximately 
uniformly distributed
among elements. From this point on, the adaptive algorithm will not work significantly differently from uniform refinements -- 
there is no area where elements should be refined with higher priority. This corresponds to the observed straightening of the convergence curves for the adaptive algorithm, which, after initial phase, have the same slope as uniform refinements of the corresponding order. Adaptivity algorithm in this case thus does not improve the rate of convergence. 
Note that in order to achieve an exponential convergence, as e.g.\ in~\cite{hermes_2014}, 
one would need to allow also increasing the polynomial degree during the adaptive process.

One can also observe that the reason for spreading of refinements away from the internal layer is not due to avoiding 
hanging nodes of a higher level. 
Although these could be restrictive in the case of a singularity, in our case, the solution is smooth (although very steep), and the refinements are needed for equilibrating the error.

In Table~\ref{table_adapt_global}, global properties of problems solved on individual adaptivity levels are presented.
Note that this can be seen as yet another type of scaling test. Unlike in the case of strong or weak scaling, 
we use the same number of processors (2048 in our case), while both global and local problem sizes increase. 

We can see from Table~\ref{table_adapt_global} that the performance of the method is good for the selected number of processors. The largest solved problem has more than 10$^9$ unknowns. Although the solution time gets rather large for the last adaptivity steps, the number of iterations remains similar.

\begin{table}
\begin{center}
\begin{tabular}{ccc|cc|c|cc}
$N_S$ & $n$ & $n/N_S$ & $n^{\Gamma}$ & $n_{C}$ & its. & $t_{set-up}$ & $t_{PCG}$ \\
\hline
2048/46 & 4913 & 2 & 4873/2144 & 4055/552 & 9  & 2.8 & 1.1 \\
2048/46 & 8594 & 4 & 8488/4134 & 7756/942 & 29  & 0.55 & 1.6 \\
2048/46 & $2.5\!\cdot\! 10^{4}$ & 12 & $2.1\!\cdot\! 10^{4}$/7978 & $1.7\!\cdot\! 10^{4}$/1005 & 53  & 0.6 & 3 \\
2048/46 & $1.3\!\cdot\! 10^{5}$ & 63 & $8.3\!\cdot\! 10^{4}$/$1.6\!\cdot\! 10^{4}$ & $3.6\!\cdot\! 10^{4}$/988 & 60  & 0.67 & 3.5 \\
2048/46 & $7.0\!\cdot\! 10^{5}$ & 342 & $3.1\!\cdot\! 10^{5}$/$2.0\!\cdot\! 10^{4}$ & $4.5\!\cdot\! 10^{4}$/1109 & 54 &  0.89 & 3.5 \\
2048/46 & $3.0\!\cdot\! 10^{6}$ & 1445 & $8.2\!\cdot\! 10^{5}$/$2.1\!\cdot\! 10^{4}$ & $4.7\!\cdot\! 10^{4}$/1152 & 56  & 1.6 & 4.8 \\
2048/46 & $1.4\!\cdot\! 10^{7}$ & 6623 & $2.4\!\cdot\! 10^{6}$/$2.0\!\cdot\! 10^{4}$ & $4.7\!\cdot\! 10^{4}$/1057 & 55  & 2.9 & 10 \\
2048/46 & $6.4\!\cdot\! 10^{7}$ & $3.1\!\cdot\! 10^{4}$ & $7.2\!\cdot\! 10^{6}$/$2.1\!\cdot\! 10^{4}$ & $4.5\!\cdot\! 10^{4}$/1011 & 55 &  10 & 33 \\
2048/46 & $2.9\!\cdot\! 10^{8}$ & $1.4\!\cdot\! 10^{5}$ & $2.0\!\cdot\! 10^{7}$/$2.2\!\cdot\! 10^{4}$ & $4.8\!\cdot\! 10^{4}$/1222 & 56  & 61 & 130 \\
2048/46 & $1.3\!\cdot\! 10^{9}$ & $6.3\!\cdot\! 10^{5}$ & $5.2\!\cdot\! 10^{7}$/$2.1\!\cdot\! 10^{4}$ & $4.7\!\cdot\! 10^{4}$/1088 & 51  & 565 & 521 \\
\end{tabular}
\caption{\label{table_adapt_global} Adaptivity in 3D, Poisson problem, 3-level BDDC, trilinear elements. For brevity, only initial and even steps 1, 2, 4, 6, 8, \ldots, 18 are shown.}
\end{center}
\end{table}

\section{Conclusions}
\label{sec:Conclusions}

We have presented an approach to employ adaptive mesh refinement within a massively parallel finite element software.
The necessary equilibrating of subdomain sizes during the adaptivity process is based on repartitioning a space-filling Z-curve, 
and the solver relies on the \codename{p4est} library in this regard. 
The goal of the partitioning strategy is to obtain subparts with equal number of elements.
While more advanced partitioning strategies exist, this set-up is, to the best of our knowledge, the only one currently available for rebalancing at thousands of computer cores.
Moreover, for real-world problems, which might be time-dependent, nonlinear or with other complexities, and can involve high-order elements, 
the assembly time might actually dominate the whole computation.
This also supports the use of this simple partitioning strategy, 
as the matrix integration and assembly algorithms are typically embarrassingly parallel with time proportional to the number of elements in the given subdomain.

For solving the arising systems of linear equations, the multilevel BDDC method has been employed. 
A parallel implementation of this method is freely available in our \codename{BDDCML} library.
The standard (two-level) and the three-level versions of the method have been used.

A relatively simple strategy for dealing with hanging nodes has been recalled. 
The assumption on 2:1 mesh regularity has been advocated as a reasonable trade-off between unnecessary enforced refinements and algorithmic complexity.
Under this assumption, hanging nodes can be eliminated on the element level while the number of local degrees of freedom of an element and, consequently, a major part of the assembly algorithm remain unchanged.
Meshes satisfying this requirement are automatically provided by the \codename{p4est} library without any extra coding effort.

We have demonstrated that this set-up fits naturally into the framework of nonoverlapping domain decomposition methods,
and it is very straightforward to pass the subdomain-wise assembled algebraic system to the \codename{BDDCML} solver for the subsequent solution.
One of the challenges of this combination seems to be the fact that disconnected subdomains composed of several components as well as loosely coupled subdomains 
are often encountered. 
The BDDC solver detects these situations and is able to accommodate them by an adjusted selection of constraints.

We have investigated parallel performance of the \codename{BDDCML} solver with respect to adaptively refined meshes. 
In particular, we have studied the performance in two and three dimensions, for Poisson and elasticity problems, using linear, quadratic and fourth-order finite elements, 
and employing two and three levels within the BDDC method.
We have presented tests of parallel strong scaling on adapted unstructured meshes as well as studies of the effect of irregular subdomains.
The solver has been also verified using convergence tests.
Finally, we have presented an adaptive finite element calculation in 3D, with the finest mesh resulting in more than 10$^9$ unknowns.

Presented results show a very good performance of the developed solver in most circumstances. 
In particular, the number of PCG iterations required for solving the system remains almost constant for all the studied cases.

The only issue seems to be the potential load imbalance caused by significantly different number of constraints on different subdomains. 
However, the number of components of a subdomain is limited, and it has reached at most two in our computations. 
Consequently, we have shown on a number of experiments that the differences from a simple case with regular subdomains remain rather small and acceptable.

\section*{Acknowledgements}
We are grateful to Santiago Badia, Javier Principe and Fehmi Cirak for valuable discussions on related topics.

This research was supported by Czech Science Foundation through grant 14--02067S, 
and by Czech Academy of Sciences through RVO:67985840.
Computational time on the \emph{Salomon} supercomputer has been provided by the IT4Innovations Centre of Excellence project (CZ.1.05/1.1.00/02.0070), 
funded by the European Regional Development Fund and the national budget of the Czech Republic via the Research and Development for Innovations Operational Programme, 
as well as Czech Ministry of Education, Youth and Sports via the project Large Research, Development and Innovations Infrastructures (LM2011033).

\bibliographystyle{elsarticle-harv}
\bibliography{literatura}

\end{document}